\documentclass[11pt]{amsart}
\usepackage{graphicx}
\usepackage{latexsym,amsmath,amsfonts,amscd, amsthm}
\usepackage{changebar}
\usepackage{color}
\usepackage{bm}
\usepackage{tikz}
\usepackage{multirow}
\usetikzlibrary{arrows,backgrounds,snakes}
\usepackage{epsfig,subfigure}
\usepackage{ulem}

\newtheoremstyle{plainNoItalics}{}{}{\normalfont}{}{\bfseries}{.}{ }{}
\theoremstyle{plain}
\newtheorem{thm}{Theorem}[section]
\theoremstyle{plainNoItalics}

\newtheorem{defn}[thm]{Definition}
\newtheorem{rem}[thm]{Remark}
\newtheorem{prop}[thm]{Proposition}

\newcommand{\f}{\frac}
\newcommand{\ra}{\rightarrow}
\newcommand{\beq}{\begin{equation}}
	\newcommand{\eeq}{\end{equation}}
\newcommand{\beqa}{\begin{eqnarray}}
	\newcommand{\eeqa}{\end{eqnarray}}
\newcommand{\bit}{\begin{itemize}}
	\newcommand{\eit}{\end{itemize}}
\newcommand{\bedef}{\begin{defn}}
	\newcommand{\edefn}{\end{defn}}
\newcommand{\bpro}{\begin{prop}}
	\newcommand{\epro}{\end{prop}}

\newcommand{\Dt}{\Delta t}

\newcommand{\mI}{{\mathbb  I}}

\newcommand\vave[1]{{\langle{#1}\rangle}}
\newcommand{\eps}{\varepsilon}

\newcommand{\bx}{{\bf x}}
\newcommand{\bv}{{\bf v}}

\newcommand{\bX}{{\bf X}}
\newcommand{\bY}{{\bf Y}}

\title[]{Asymptotic preserving and uniformly unconditionally stable finite difference schemes for kinetic transport equations} 

\begin{document}

\maketitle

\centerline{\scshape Guoliang Zhang}
\medskip
{\footnotesize
\centerline{Institute of Natural Sciences, Shanghai Jiao Tong University, Shanghai 200240, PR China }
\centerline{zglmath@sjtu.edu.cn}
}
	
\medskip

\centerline{\scshape Hongqiang Zhu}
\medskip
{\footnotesize
	\centerline{School of Natural Sciences, Nanjing University of Posts and Telecommunications, Nanjing, Jiangsu 210023, PR China  }
	\centerline{zhuhq@njupt.edu.cn}
}

\medskip

\centerline{\scshape Tao Xiong\footnote{Corresponding author.}}
\medskip
{\footnotesize
	\centerline{School of Mathematical Sciences, Fujian Provincial Key Laboratory of Mathematical Modeling} 
	\centerline{and High-Performance Scientific Computing, Xiamen University}
	\centerline{Xiamen, Fujian 361005, PR China}
	\centerline{txiong@xmu.edu.cn}
}

	
	
	\begin{abstract}
		In this paper, uniformly unconditionally stable first and second order finite difference schemes are developed for kinetic transport equations in the diffusive scaling. We first derive an approximate evolution equation for the macroscopic density, from the formal solution of the distribution function, which is then discretized by following characteristics for the transport part with a backward finite difference semi-Lagrangian approach, while the diffusive part is discretized implicitly. After the macroscopic density is available, the distribution function can be efficiently solved even with a fully implicit time discretization, since all discrete velocities are decoupled, resulting in a low-dimensional linear system from spatial discretizations at each discrete velocity. Both first and second order discretizations in space and in time are considered. The resulting schemes can be shown to be asymptotic preserving (AP) in the diffusive limit.  Uniformly unconditional stabilities are verified from a Fourier analysis based on eigenvalues of corresponding amplification matrices. Numerical experiments, including high dimensional problems, have demonstrated the corresponding orders of accuracy both in space and in time, uniform stability, AP property, and good performances of our proposed approach.
	\end{abstract}
	
\vspace{0.1cm}

	

	\section{Introduction}
	\label{sec_intro}
	\setcounter{equation}{0}
	\setcounter{figure}{0}
	\setcounter{table}{0}
	
	Kinetic transport equations have been widely used in rarefied gas dynamics \cite{Bird1994}, neutron transport \cite{case1967} and radiative transfer \cite{Chandrasekhar1960,pomraning2005equations}, etc. In this work, we are interested in the nonstationary kinetic transport equation in the diffusive scaling:
	\begin{equation}
		\label{model}
		\eps\,f_t + v\,f_x = \frac{1}{\eps}\mathcal{C}(f),
	\end{equation}
	where $f=f(x,v,t)$ is the probability distribution function of particles. For simplicity, here we take one-dimension (1D) in space with $x\in \Omega_x \subset \mathbb{R}$, and 1D in velocity where $v\in \Omega_v \subset \mathbb{R}$ and time $t\geq 0$. For the velocity, we may either have $v\in \{-1, 1\}$ for a discrete-velocity model, or $v\in [-1,1]$ for a one-group velocity model. $\mathcal{C}(f)$ is a collision operator which describes the interaction between particles, or particles and their outer media. $\varepsilon$ is the dimensionless Knudsen number, which is defined as the ratio of the mean free path over the characteristic length of the system.
	Multidimensions of \eqref{model} can be formulated similarly.
	
	Owing to the multi-scale property of \eqref{model}, from both time stiffness and collision operator, it has great numerical challenges to develop effective and efficient algorithms. First of all, a standard explicit numerical scheme would suffer from a strict time step condition $\Delta t\le\mathcal{O}(\eps h)$, where $h$ is the mesh size, which is prohibitive as $\eps\ra 0$; Secondly, direct implicit time discretization may avoid the stability issue, however, it may fail to capture the correct asymptotic limit as $\eps\ra 0$ on underresolved meshes \cite{caflisch1997uniformly,naldi1998numerical}. Besides, fully implicit schemes may couple the space $x$ and the velocity $v$ together, leading to a large size matrix, which requires advanced iterative methods especially for multidimensions \cite{klar1998asymptotic,li2017implicit}.
	
	
	To address such issues, many different methods have been developed in literature, among them asymptotic preserving (AP) schemes \cite{jin1999efficient} have offered a great framework. The main idea of an AP scheme is to
	preserve the asymptotic behavior of the physical model on the discrete level \cite{jin2010asymptotic,jin2022asymptotic}. For the problem \eqref{model}, which leads to a diffusion equation as $\eps\ra0$, an AP scheme would converge to a scheme for the limiting diffusion equation.
	There are several different ways to build an AP scheme for \eqref{model}, such as relaxation schemes \cite{jin1998diffusive,jin2000uniformly,boscarino2013implicit,albi2020implicit}, schemes based on micro-macro reformulation \cite{lemou2008new,jang2015high,crestetto2019asymptotically,peng2020stability,peng2021asymptotic}, unified gas-kinetic (UGK) type schemes \cite{mieussens2013asymptotic,sun2015,shi2020,li2020unified}, etc. For these AP schemes, the time stability condition relying on the parameter $\eps$ can be avoided, especially as $\eps\rightarrow 0$ in the diffusive regime. However, due to the limiting diffusion equation, some of these schemes, such as \cite{jin1998diffusive,jin2000uniformly,lemou2008new,jang2015high}, will subject to a parabolic time step condition, namely $\Delta t=\mathcal{O}(h^2)$ in the diffusive regime $\eps\ll1$ (see the corresponding analyses given in \cite{klar2002uniform,liu2010analysis,jang2014analysis}). We name these schemes AP explicit schemes. 
	
	On the other hand, to achieve unconditional stability in the diffusive regime, in \cite{boscarino2013implicit,peng2020stability,peng2021stability}, a weighted diffusive term is added to both sides of the original equation, one is discretized explicitly and the other implicitly, resulting a limiting scheme being implicit to the diffusion equation. Some other schemes propose to modify their AP explicit schemes by applying an implicit treatment to diffusion-linked terms, see e.g. \cite{mieussens2013asymptotic,crestetto2019asymptotically,albi2020implicit,peng2021asymptotic}. We name these schemes AP implicit schemes.
	Although such AP implicit schemes are unconditionally stable in the diffusive limit, it is still conditionally stable away from the diffusive regime. For these schemes, especially in the intermediate regime when $\eps$ is comparable with the mesh size, the time step restriction is still $\mathcal{O}(\eps h)$ (the CFL number relies on $\eps$) \cite{peng2021stability,peng2021asymptotic}, which is small and not very pleasant for some real applications, especially involving long time simulations. 
	
	{In this work, we would like to propose a uniformly unconditionally stable AP scheme for all range of $\eps$. Noticing that the transport velocity $v/\eps$ in \eqref{model} is independent of $x$, we may overcome the time stability restriction from this advective term by following characteristics. 
		Our main novel idea is to first derive an approximation model for the macroscopic density, whose flux consists of a part from transport by following the characteristics, and another part coming from diffusion. The main advantage of this model equation is that it is built based on two main physical procesess during kinetic transport: transportation and diffusion. Numerically we can efficiently update the density by using a semi-Lagrangian finite difference scheme for the transport part, and a linearly implicit finte difference discretization for the diffusion part.
		The resulting scheme for the macroscopic density is uniformly unconditionally stable, and only a low-dimensional linear system in space needs to be solved. After the density is available, the microscopic equation for the distribution function, which is high-dimensional, is decoupled for all discrete velocities if a discrete-ordinate method (DOM) is used. It is very convenient for parallel computing, so that this approach is very friendly for high-dimensional problems. Besides, we utilize the original microscopic equation, as it has a constant transport speed for each fixed velocity. Numerically we can either use a semi-Lagrangian finite difference scheme, or a linear implicit finite difference scheme. Both are also uniformly unconditionally stable and only a low-dimensional linear system in space needs to be solved. The overall scheme for the macro-micro system is then expected to be uniformly unconditionally stable. We would emphasize that this framework is very different from the micro-macro reformulation \cite{lemou2008new,jang2015high,crestetto2019asymptotically,peng2020stability,peng2021asymptotic}. As the latter is based on a modeling by splitting the distribution function into equilibrium and non-equilibrium parts, an orthogonal projection coupling the whole velocity space is used, which makes it not easy to design uniformly unconditionally stable schemes. Our approach is similar to UGK type schemes \cite{sun2015,shi2020,li2020unified}, but the main difference is that our new model includes a long characteristic tracking for transport. We propose both first and second order discretizations in space and in time. Uniformly unconditional stability is verified by a Fourier analysis for the two-discrete-velocity telegraph equation, with asympotic preserving (AP) property formally proved, for both fully discrete schemes. We note that this stability analysis is based on the eigenvalues of the corresponding amplification matrix. Numerical experiments, including high-dimensional problems, will validate the designed orders of accuracy, AP, large time step stablity, higher efficiency, and good performances of our proposed approach.}
	
	The rest of the paper is organized as follows. In Section 2, the model equation is reviewed. The approximation model is derived in Section 3. In Section 4, the first and second order schemes are proposed, and their uniformly unconditional stability is analyzed by a Fourier analysis. In Section 5, numerical experiments are performed to verify the orders of accuracy, the effectiveness and good performance of the proposed methods, including high-dimensional problems. Concluding remarks are made in Section 6.
	
	\section{Model equation}
	\label{sec_model}
	\setcounter{equation}{0}
	\setcounter{figure}{0}
	\setcounter{table}{0}
	
	We rewrite the linear transport kinetic equation \eqref{model} as
	\begin{equation}
		\label{eq1}
		f_{t}+\frac{1}{\varepsilon} v f_{x} =\frac{1}{\varepsilon^{2}}\,\mathcal{C}(f).
	\end{equation}
	We consider the following four collision operators for $\mathcal{C}(f)$ as in \cite{jang2015high}
	\begin{subequations}
		\label{Qf}
		\begin{align}
			\mathcal{C}(f)\,=\, &\langle f\rangle -f,  \label{Q1}\\
			\mathcal{C}(f)\,=\, & K\langle f\rangle ^m\,(\langle f\rangle -f),    \label{Q2}\\
			\mathcal{C}(f)\,=\,& \langle f\rangle -f+A\varepsilon v \langle f\rangle ,  \quad |A\varepsilon < 1|,   \label{Q3}\\
			\mathcal{C}(f)\,=\,& \langle f\rangle -f+C\varepsilon v [\langle f\rangle ^2-(\langle f\rangle -f)^2], \quad C> 0,   \label{Q4}
		\end{align}
	\end{subequations}
	where we denote $\langle f\rangle \,:=\,\int f d \mu $ and $d \mu$ is either a discrete Lebesgue measure on the discrete velocity space $\{-1,1\}$, namely
	\begin{equation}
		\label{vave_disv}
		\langle f\rangle =\frac{f(x, v=1, t)+f(x, v=-1, t)}{2},
	\end{equation}
	or $d\mu=\f12\,dv$ for the one-group velocity of $v\in[-1,1]$
	\begin{equation}
		\label{vave_groupv}
		\langle f\rangle =\f12\int^1_{-1}f(x,v,t)dv.
	\end{equation}
	The system \eqref{eq1} with the collision operator \eqref{Q1} is known as the Goldstein-Taylor model or the telegraph equation.
	
	Now if we define the macroscopic density $\rho$ and current $g$ for the probability density function $f$ as
	\begin{equation}
		\rho=\langle f\rangle, \quad f=\rho+\varepsilon g,
	\end{equation}
	where $\langle\cdot\rangle$ is either from \eqref{vave_disv} or \eqref{vave_groupv}, we may derive their corresponding diffusive limiting equations for the macroscopic density $\rho$ from \eqref{eq1} with the collision operators \eqref{Qf}. Taking \eqref{Q1} as an example, we introduce two orthogonal projections $\Pi=\langle\cdot \rangle$ and $\mathbb{I}- \Pi$, which are defined on $L^2(\Omega_v)$. Applying them to \eqref{eq1}, we obtain
	\begin{subequations}
		\begin{align}
			&\rho_t+\langle vg\rangle_x=0,  \label{eq2_1} \\
			&\varepsilon^2 g_t+\varepsilon(\mI-\Pi)(vg)_x+ v\,\rho_x=-g.\label{eq2_2}
		\end{align}
	\end{subequations}
	Substituting $g=-v \rho_x-\eps(\mI-\Pi)(vg)_x-\eps^2 g_t$ from \eqref{eq2_2} into \eqref{eq2_1}, as $\eps \rightarrow 0$, it yields the following heat equation for the density $\rho$
	\begin{equation}
		\rho_t-\vave{v^2}\,\rho_{xx}=0,
		\label{lQ1}
	\end{equation}
	where $\vave{v^2}=1$ with \eqref{vave_disv} and $\vave{v^2}=\f13$ from \eqref{vave_groupv}.
	Similarly, for other collision operators, corresponding diffusive limiting equations for the density $\rho$ can be obtained, see \cite{jin1998diffusive,jang2015high} for more details. For the collision operators in \eqref{Qf}, the corresponding diffusive limiting equations, which we accumulate together, are given as
	\begin{subequations}
		\label{limiteq}
		\begin{align}
			\label{heat}
			&\rho_t=\vave{v^2}\,\rho_{xx},   \quad \text{heat equation}; \\
			\label{porous}
			&\rho_t=\langle v^2\rangle \frac{1}{K(1-m)}(\rho^{1-m})_{xx},   \quad \text{porous media equation}; \\
			\label{adv-diff}
			&\rho_t+A\langle v^2\rangle\rho_x=\langle v^2\rangle\rho_{xx}, \quad \text{linear advection-diffusion equation};\\
			\label{visburg}
			&\rho_t+C\langle v^2\rangle(\rho^2)_x=\langle v^2\rangle\rho_{xx}, \quad \text{viscous nonlinear Burgers' equation}.
		\end{align}
	\end{subequations}
	
	\section{Model approximation}
	\label{sec_model_apprx}
	\setcounter{equation}{0}
	\setcounter{figure}{0}
	\setcounter{table}{0}
	
	In this section, we will derive an approximation model for the macroscopic density $\rho$, from which we can update $\rho$ directly. The idea is to substitute a formal solution of $f$ in terms of $\rho$ into the flux function of \eqref{eq1}, and then average it in the velocity space. This technique has been used in developing the UGK scheme, e.g. in \cite{mieussens2013asymptotic,li2020unified}, in which the resulting flux is then approximated by numerical flux reconstructions. The main difference is that we further integrate by parts and approximate the integrals along the $v$-direction in the macroscopic flux to get a new approximation model. The approximation errors can be accurately analyzed in a formal form. The new model includes a long characteristic tracking for transport, but still easily leads to a correct asymptotic diffusive limit.
	
	We first define the characteristics for \eqref{eq1} as
	\beq
	\label{eq31}
	\frac{d \mathbf{X}(t)}{d t}=\frac{v}{\varepsilon}.
	\eeq
	{Starting from a given point $(x_*,t_*)$ with a velocity $v$, the characteristic curve can be written as a four-variable function, which is traced back via
		\beq
		\label{eq32}
		\mathbf{X}(t, x_*, t_*, v)=x_*-\frac{v}{\varepsilon} (t_*-t).
		\eeq
		Along the characteristic curve, the main unknown $f(x,v,t)$ now can be written as $f(\mathbf{X}(t,x_*,t_*,v),v,t)$, and its corresponding material derivative is
		\begin{equation}
			\label{eq33}
			\begin{aligned}
				\frac{d f(\mathbf{X}(t,x_*,t_*,v), v, t)}{d t} =f_x \frac{d \mathbf{X}(t,x_*,t_*,v)}{d t} + f_t =f_t+\frac{v}{\varepsilon} \, f_x.
			\end{aligned}
		\end{equation}
		With this, \eqref{eq1} is equivalent to
		\begin{equation}
			\label{eq34}
			\frac{d f(\mathbf{X}(t,x_*,t_*,v), v, t)}{d t} =\frac{1}{\varepsilon^{2}}\,Q(f).
		\end{equation}
		Taking the telegraph equation with the collision operator \eqref{Q1} as an example, \eqref{eq34} becomes
		\begin{equation}
			\label{eq35}
			\frac{d f(\mathbf{X}(t,x_*,t_*,v), v, t)}{d t} =\frac{1}{\varepsilon^{2}} \left(\rho(\mathbf{X}(t,x_*,t_*,v),t) -f(\mathbf{X}(t,x_*,t_*,v),v,t)\right),
		\end{equation}
		or equivalently
		\begin{equation}
			\label{eq36}
			\frac{d f(\mathbf{X}(t,x_*,t_*,v), v, t)}{d t} +\frac{1}{\varepsilon^{2}}f(\mathbf{X}(t,x_*,t_*,v), v, t)=\frac{1}{\varepsilon^{2}}\rho(\mathbf{X}(t,x_*,t_*,v), t).
		\end{equation}
		Now denoting $\mu=\frac{1}{\varepsilon ^2}$ and integrating \eqref{eq36} over $(t_n, t)$ after multiplying the integrating factor $e^{\mu t}$, we obtain a formal
		solution for $f$ in terms of $\rho$
		\begin{equation}
			\label{eq37}
			f(\mathbf{X}(t,x_*,t_*,v), v, t)=e^{-\mu \left(t-t_{n}\right)} f\left(\mathbf{X}(t_{n},x_*,t_*,v), v, t_{n}\right)+ \int_{t_n}^{t} \mu e^{-\mu (t-s)} \rho\left(\mathbf{X}(s,x_*,t_*,v), s\right) ds.
		\end{equation}
		When $t=t_*$, $\bX(t_*,x_*,t_*,v)=x_*$, the unknown $f$ at the point $(x_*,t_*)$ is given by
		\begin{equation}
			\label{eq38}
			f(x_*,v,t_*)=e^{-\mu \left(t_*-t_{n}\right)} f\left(\mathbf{X}(t_n,x_*,t_*,v), v, t_{n}\right)+ \int_{t_n}^{t_*} \mu e^{-\mu (t_*-s)} \rho(\mathbf{X}(s,x_*,t_*,v), s) ds.
		\end{equation}
	}
	
	In what follows, we will derive an equation for the macroscopic density $\rho$, based on the formal solution of $f$ in \eqref{eq38}. Without confusion, we now drop the subindex $*$ for clarity. We start with integrating \eqref{eq1} over $v$, which yields
	\begin{equation}
		\label{eq310}
		\rho_{t}+\frac{1}{\varepsilon}\langle v f \rangle_x  \,= \,0.
	\end{equation}
	Substituting the formal solution \eqref{eq38} into \eqref{eq310}, we have
	\begin{equation}
		\label{eq311}
		\rho_t+\frac{e^{-\mu \left(t-t_{n}\right)}}{\varepsilon} \left\langle vf_x\left(\bX(t_n,x,t,v),  v,  t_{n}\right)\right\rangle
		+ \frac{1}{\varepsilon} \left\langle v\int_{t_n}^{t} \mu e^{-\mu (t-s)} \rho_x(\bX(s,x,t,v), s)\, ds \right\rangle =0.
	\end{equation}
	\eqref{eq311} does not explicitly guarantee a right diffusive limit. In \cite{mieussens2013asymptotic}, the authors proposed to approximate $\rho_x(\bX(s,x,t,v), s)$ by a linear reconstruction from a Taylor expansion at the point $(x,t_n)$, which would lead to time stability restrictions. Our idea is to further perform an integration by parts for the third term in \eqref{eq311}, 
	\begin{align}
		 \int_{t_n}^{t} \mu e^{-\mu (t-s)} \rho_x(\bX(s,x,t,v), s)\, ds=&\int_{t_n}^{t}  \rho_x(\bX(s,x,t,v), s)\, d\, e^{-\mu (t-s)} \notag \\
		=&\rho_x(x, t)-\rho_x(\bX(t_n,x,t,v), s)e^{-\mu (t-t_n)} \notag \\
		&-\int_{t_n}^{t}  e^{-\mu (t-s)}\left(\frac{v}{\eps}\,\rho_{xx}(\bX(s,x,t,v), s)+\rho_{xt}(\bX(s,x,t,v), s)\right)d\,s , \notag
	\end{align}
	here $d\bX(s,x,t,v)/ds=v/\eps$ is used. Due to $\vave{v}=0$ and $\mu=1/\eps^2$, \eqref{eq311} gives
	\begin{align}
		\label{eq312} \rho_t&+\frac{e^{-\mu \left(t-t_{n}\right)}}{\varepsilon} \langle v\,(f_x(\bX(t_n,x,t,v),  v,  t_{n})-\rho_x(\bX(t_n,x,t,v), t_{n}))\rangle \\ &-\left\langle  \int_{t_n}^{t}  v^2 \, \rho_{xx}(\bX(s,x,t,v), s)\,\mu\,e^{-\mu (t-s)}ds \right\rangle
		-\frac{1}{\varepsilon}\left\langle  \int_{t_n}^{t} v \rho_{xt}(\bX(s,x,t,v), s)) e^{-\mu (t-s)} ds\right\rangle=0, \notag
	\end{align}
	where $\rho_{xx}$ is the second derivative with respect to its first argument and $\rho_{xt}$ is the second order mixed derivatives for the first and second arguments respectively. Up to \eqref{eq312}, there is no approximation error yet.
	
	To make \eqref{eq312} applicable as a numerical scheme, we take the following two approximations
	\begin{subequations}
		\label{eq313}
		\begin{align}
			\label{eq313a}
			&\int_{t_{n}}^{t} \mu e^{-\mu\left(t-s\right)} \rho_{xx}(\bX(s,x,t,v), s) d s
			\approx \rho_{xx}(x,t)\int_{t_{n}}^{t} \mu e^{-\mu\left(t-s\right)}  d s = (1-e^{-\mu(t-t_n)})\rho_{xx}(x,t), \\
			\label{eq313b}
			&\int_{t_{n}}^{t} \frac{1}{\varepsilon} e^{-\mu\left(t-s\right)} \rho_{xt}(\bX(s,x,t,v), s) d s
			\approx \rho_{xt}(x,t)\int_{t_{n}}^{t} \frac{1}{\varepsilon} e^{-\mu\left(t-s\right)}  d s =\varepsilon (1-e^{-\mu(t-t_n)})\rho_{xt}(x,t).
		\end{align}
	\end{subequations}
	That is, we extract $\rho_{xx}(\bX(s,x,t,v),s)$ and $\rho_{xt}(\bX(s,x,t,v),s)$ by taking $s=t$ so that $\bX(t,x,t,v)=x$, while the remaining integration can be done explicitly. With these two approximations, since $\vave{v}=0$ in the last term on the left hand side, \eqref{eq312} becomes
	\begin{equation}
		\label{eq314}
		\rho_t+\frac{e^{-\mu \left(t-t_{n}\right)}}{\varepsilon} \langle v(f-\rho)_x(\bX(t_n,x,t,v),  v,  t_{n})\rangle
		-\langle v^2\rangle (1-e^{-\mu(t-t_n)})\rho_{xx}\left(x,t\right)=0.
	\end{equation}
	Here in the second term on the left, we abuse the notation by using $\rho_x(\bX(t_n,x,t,v),  v,  t_{n})$ with an extra $v$, in order to make \eqref{eq314} to be more concise, similarly in the following. We refer to \eqref{eq312} for a complete expression.
	We can easily see that \eqref{eq314} formally converges to the heat equation
	\[
	\rho_t - \vave{v^2}\,\rho_{xx} = 0,
	\]
	for $t>t_n$, since
	\[
	\frac{1}{\varepsilon}e^{-\mu(t-t_n)} \rightarrow 0, \quad e^{-\mu(t-t_n)} \rightarrow 0,  \text{ as } \eps \ra 0.
	\]
	
	For the model approximation used in \eqref{eq313}, we have the following theorem for the approximation errors.
	\begin{thm}
		\label{thm1}
		For $t\in(t_n, t_{n+1}]$ and denoting $\Delta t=t_{n+1}-t_n$, the approximation errors in \eqref{eq313} {can be formally} given as
		\begin{subequations}
			\label{modelappr2}
			\begin{align}
				\label{eq315a}
				&\int_{t_{n}}^{t} \mu e^{-\mu\left(t-s\right)} \rho_{xx}(\bX(s,x,t,v), s)\, d\,s
				= (1-e^{-\mu(t-t_n)})\rho_{xx}(x,t) + \mathcal{O}\left(\frac{\eps\,\Dt^2}{2\eps^4+2\eps^2\Dt+\Dt^2}\right); \\
				\label{eq315b}
				&\int_{t_{n}}^{t} \frac{1}{\varepsilon} e^{-\mu\left(t-s\right)} \rho_{xt}(\bX(s,x,t,v), s)\,d\,s
				=\varepsilon (1-e^{-\mu(t-t_n)})\rho_{xt}(x,t)+ \mathcal{O}\left(\frac{\eps^2\,\Dt^2}{2\eps^4+2\eps^2\Dt+\Dt^2}\right).
			\end{align}
		\end{subequations}
	\end{thm}
	
	\begin{proof}
		By a Taylor expansion, from \eqref{eq313a} we have
		\begin{equation}
			\label{err1}
			\begin{aligned}
				\int_{t_{n}}^{t}& \left(\rho_{xx}(x-\frac{v}{\varepsilon}(t-s), s)-\rho_{xx}\left(x, t\right) \right) \mu e^{-\mu\left(t-s\right)}ds \\
				&=\int_{t_{n}}^{t}
				\sum_{\ell=1}^{\infty}\frac{(s-t)^\ell}{\ell !}(\frac{v}{\varepsilon}\partial_{x} +\partial_{t})^\ell\rho_{xx}(x, t)\, \mu e^{-\mu\left(t-s\right)}ds \\
				&=\sum_{\ell=1}^{\infty}\,(v\,\partial_{x} +\eps\,\partial_{t})^\ell\rho_{xx}(x, t)\,\int_{t_{n}}^{t}
				\frac{(s-t)^\ell}{\eps^\ell\,\ell !} \mu e^{-\mu\left(t-s\right)}ds.
			\end{aligned}
		\end{equation}
		If we assume $(v\,\partial_{x} +\eps\,\partial_{t})^\ell\rho_{xx}(x, t)=\mathcal{O}(1)$ for any $\ell$, we only need to estimate the integral term $\int_{t_{n}}^{t}\frac{1}{\varepsilon^\ell\,\ell!}\,(s-t)^{\ell}\, \mu e^{-\mu\left(t-s\right)}ds$. Notice that $\mu=1/\eps^2$, for any fixed $\eps>0$, we have
		\begin{equation}
			\label{eq316}
			\begin{aligned}
				\int_{t_{n}}^{t}(s-t)^\ell\, \mu e^{-\mu\left(t-s\right)}ds&=(s-t)^\ell\,e^{-\mu\left(t-s\right)}\Bigg|^{s=t}_{s=t_n}-\ell\int_{t_{n}}^{t}(s-t)^{\ell-1}\, e^{-\mu\left(t-s\right)}ds \\
				&=-(t_n-t)^\ell\,e^{-\mu(t-t_n)}-\ell \int_{t_{n}}^{t}\,(s-t)^{\ell-1}\, e^{-\mu\left(t-s\right)}ds \\
				&=(-1)^\ell\, \ell!\, \eps^{2\ell}e^{-\mu(t-t_n)}\left(e^{\mu(t-t_n)}-\sum_{j=0}^\ell\frac{\big(\mu\,(t-t_n)\big)^j}{j!}\right),
			\end{aligned}
		\end{equation}
		so that \eqref{err1} gives
		\begin{equation}
			\label{err12}
			\begin{aligned}
				\int_{t_{n}}^{t}& \left(\rho_{xx}(x-\frac{v}{\varepsilon}(t-s), s)-\rho_{xx}\left(x, t\right) \right) \mu e^{-\mu\left(t-s\right)}ds \\
				&=\sum_{\ell=1}^{\infty}(-1)^\ell(v\,\partial_{x} +\eps\,\partial_{t})^\ell\rho_{xx}(x, t)\,\eps^{\ell}e^{-\mu(t-t_n)}\left(e^{\mu(t-t_n)}-\sum_{j=0}^\ell\frac{\big(\mu\,(t-t_n)\big)^j}{j!}\right).
			\end{aligned}
		\end{equation}
		In \eqref{err12}, the main part of the error which depends on $\eps$ and $\Dt$ is
		\beq
		\label{errd}
		\eps^{\ell}e^{-\mu(t-t_n)}\left(e^{\mu(t-t_n)}-\sum_{j=0}^\ell\frac{\big(\mu\,(t-t_n)\big)^j}{j!}\right).
		\eeq
		We can see the summation in the brackets is the first $(\ell+1)$-term Taylor expansion of $e^{\mu(t-t_n)}$ and
		\[
		0\le e^{-\mu(t-t_n)}\left(e^{\mu(t-t_n)}-\sum_{j=0}^\ell\frac{\big(\mu\,(t-t_n)\big)^j}{j!}\right) <1, \text{ for } t>t_n.
		\]
		The leading order term in \eqref{errd} is for $\ell=1$, that is
		\beq\label{errd1}
		\eps\,e^{-\mu(t-t_n)}\left(e^{\mu(t-t_n)}-\big(1+\mu\,(t-t_n)\big)\right)=\mathcal{O}\left( \frac{\eps\,\Dt^2}{2\eps^4+2\eps^2\Dt+\Dt^2}\right),
		\eeq
		where $e^{\mu(t-t_n)}$ has been expanded up to second order, and it results \eqref{eq315a}.
		Similar analysis holds for \eqref{eq315b}, so it is omitted.
	\end{proof}
	
	From above, we arrive at the following approximation model
	\begin{equation}
		\label{apprxmodel}
		\left\{
		\begin{array}{l}
			f_{t}+\frac{1}{\varepsilon} v f_{x} =\frac{1}{\varepsilon^{2}}(\rho-f), \\
			\rho_t\left(x,  t\right)+\frac{e^{-\mu \left(t-t_{n}\right)}}{\varepsilon} \langle v(f_x(\bX(t_n,x,t,v),  v,  t_{n})-\rho_x(\bX(t_n,x,t,v), t_{n}))\rangle
			-\langle v^2\rangle (1-e^{-\mu(t-t_n)})\rho_{xx}\left(x,  t\right) =0,\\
		\end{array}
		\right.
	\end{equation}
	for $t\in(t_n,t_{n+1}]$. Its approximation error \eqref{errd1} is of $\mathcal{O}(\Delta t^2)$ when $\eps=\mathcal{O}(1)$, and is of $\mathcal{O}(\eps)$ when $\eps$ is small, namely $\eps^2\ll\,C\,\Delta t$ for some positive number $C$. We can also easily see that the model has a correct formal asymptotic limit. As $\eps\ra0$, it converges to
	\begin{equation}
		\label{apprxmodellimit}
		\left\{
		\begin{array}{l}
			f\left(x,\,v,\,t\right)=\rho\left(x,  t\right), \\
			\rho_t\left(x,  t\right)-\langle v^2\rangle \rho_{xx}\left(x,  t\right) =0,\\
		\end{array}
		\right.
	\end{equation}
	since $\eps^{-k}\,e^{-\mu(t-t_n)}\ra0$ for any integer $k\ge0$ and $t\in(t_n,t_{n+1}]$.
	
	The approximation models for other collision operators in \eqref{Qf} can be similarly derived. To better illustrate our idea for developing numerical schemes in the next section, we put their approximation models in the numerical section in Section 5.
	
	\begin{rem}
		Our approximation model \eqref{apprxmodel} has a similar spirit to the AP implicit scheme developed in \cite{crestetto2019asymptotically} (see Eq. (3.17) in this paper), but the procedures are very different. For our approach, the approximation errors in \eqref{eq313} are clearly identified. {The model approximations may be extended to higher orders, however new challenges appear due to some emerging mixed derivatives, e.g., the last term in \eqref{eq312} will not be zero any more. We leave it for our future studies.}
	\end{rem}
	
	\section{Numerical scheme}
	\label{sec_scheme}
	\setcounter{equation}{0}
	\setcounter{figure}{0}
	\setcounter{table}{0}
	
	In this section, we will propose a first and a second order finite difference schemes for the approximation model \eqref{apprxmodel}. For time discretizations, we start with a first order backward Euler discretization, and extend it to second order with a BDF method. The advective terms follow the characteristics using a semi-Lagrangian scheme \cite{cho2021conservative1}. In space, an upwind discretization for advective terms, and central differences for diffusive terms are used. The resulting schemes can be shown to be uniformly unconditionally stable, in the sense that, the time step is not limited by any value of $\eps$ and spatial mesh sizes.
	
	For simplicity, we take a uniform mesh in space with grid points to be $x_j, ~~j=0,1,...,N_x-1$, and the mesh size is $\Delta x=x_{j+1}-x_j$. The time step is defined to be $\Delta t=t_{n+1}-t_n$ for the time interval $(t_n,t_{n+1}]$. For the second order BDF scheme, uniform time steps are used.
	
	We consider the discretization in a method-of-line transpose approach. Namely, time is discretized first, followed by spatial discretizations.
	
	\subsection{Temporal discretization}
	\label{sec_scheme_time}
	The first order backward Euler and the second order BDF temporal discretizations can be simply defined as follows
	\begin{itemize}
		\item \textbf{First order backward Euler scheme in time}
		\begin{subequations}
			\label{time_dis1}
			\begin{align}
				\label{time_dis11}
				&\frac{\rho^{n+1}-\rho^{n}}{\Delta t}+\frac{e^{-\mu \Delta t}}{\varepsilon} \langle v(f-\rho)_x\left(\bX(t_n,x,t_{n+1},v),  v,  t_{n}\right)\rangle
				-\langle v^2\rangle (1-e^{-\mu\Delta t})\,\rho_{xx}^{n+1} =0,\\
				\label{time_dis12}
				&\frac{f^{n+1}-f^{n}}{\Delta t}+\frac{v}{\varepsilon} \,f_x^{n+1}  =\frac{1}{\varepsilon^2} (\rho^{n+1}-f^{n+1});
			\end{align}
		\end{subequations}
		\item \textbf{Second order BDF scheme in time}
		\begin{subequations}
			\label{time_dis2}
			\begin{align}
				\label{time_dis21}
				&\frac{\rho^{n-1}-4\rho^{n}+3\rho^{n+1}}{2\Delta t}
				+\frac{e^{-\mu \Delta t}}{\varepsilon} \langle v(f-\rho)_x\left(\bX(t_n,x,t_{n+1},v),  v,  t_{n}\right)\rangle
				-\langle v^2\rangle (1-e^{-\mu\Delta t})\,\rho_{xx}^{n+1} =0,\\
				\label{time_dis22}
				&\frac{f^{n-1}-4f^{n}+3f^{n+1}}{2\Delta t}+\frac{v}{\varepsilon} \, f_x^{n+1}  =\frac{1}{\varepsilon^2} (\rho^{n+1}-f^{n+1}).
			\end{align}
		\end{subequations}
	\end{itemize}
	The first and second orders of accuracy in time can be easily verified.
	
	\subsection{Spatial discretization}
	\label{sec_scheme_space}
	
	Corresponding to the  first and second order temporal discretizations \eqref{time_dis1} and \eqref{time_dis2}, in the following we will introduce a first and a second order finite difference spatial discretizations respectively. There are mainly two parts in \eqref{time_dis11} and \eqref{time_dis21}: one is a transport part which can be traced back by following characteristics; the other is a diffusion part. For \eqref{time_dis12} and \eqref{time_dis22}, both of them only have one advection term which requires a careful treatment. For transport or advection  terms, we will use an upwind discretization, while a central difference will be applied to diffusion terms.
	
	\subsubsection{First order finite difference scheme}
	\label{subsec_1st}
	
\begin{figure}[!ht]
		\centering
		\includegraphics[width=1\textwidth]{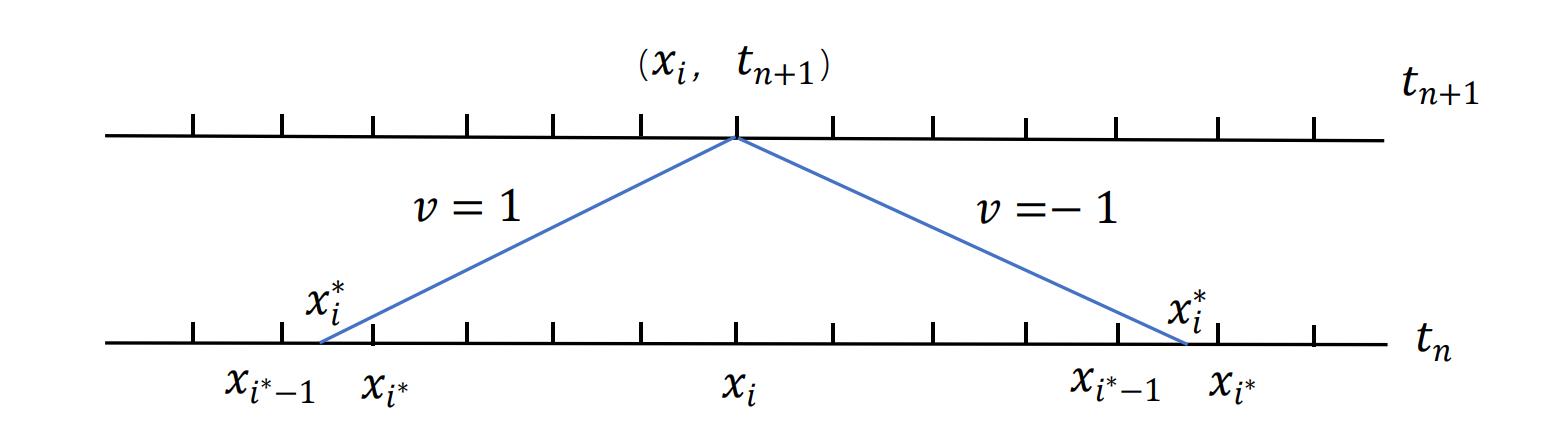}
		\caption{Locating the foot $x_i^{\ast}=\bX(t_n,v)$ by tracing back the characteristic line \eqref{eq43} when $v=1,-1$. }
		\label{shiyitu}
	\end{figure}

	We begin with the macroscopic equation of $\rho$. First,
	starting from the point $x_i$ at time level $t_{n+1}$, the characteristic line is defined as
	\begin{equation}
		\label{eq43}
		\left\{
		\begin{array}{l}
			\frac{d \bX(t)}{dt}=\frac{v}{\varepsilon}, \\
			\bX(t_{n+1})=x_i.
		\end{array}
		\right.
	\end{equation}
	The foot of the characteristics $\bX(t)$ at time level $t_n$, which is denoted as $x_i^{\ast}$, can be located exactly as
	\beq
	\label{eq44}
	\bX(t_n,v)=x_i^{\ast}= x_i-\frac{v}{\varepsilon} \Delta t
	\eeq
	when $v\in \{-1,1\}$, see Fig.\ref{shiyitu} for the illustration. {Here in \eqref{eq44}, according to the notation in \eqref{eq32}, we denote the characteristic curve as $\bX(t_n,v)$ and we omit the starting point $(x_i,t_{n+1})$ as variables for clarity.}
	
	Let us take the discrete-velocity model with $v\in\{-1,1\}$ as an example for the following description.
	The one-group velocity model with $v\in[-1,1]$ can be defined similarly, with a DOM at those Gaussian quadrature points.
	\begin{enumerate}
		\item {\bf Macroscopic equation of $\rho$ in \eqref{time_dis11}}. \\
		To be clear, we rewrite (\ref{time_dis11}) as
		\beq
		\label{eq45}
		\begin{aligned}
			\frac{\rho^{n+1}-\rho^{n}}{\Delta t}&+\frac{e^{-\mu \Delta t}}{2\varepsilon} \left((f-\rho)_x\left(x-\frac{1}{\varepsilon}\Delta t,  1,  t_{n}\right)-(f-\rho)_x\left(x+\frac{1}{\varepsilon}\Delta t,  -1,  t_{n}\right) \right) \\
			&- (1-e^{-\mu\Delta t})\,\rho_{xx}^{n+1} = 0.
		\end{aligned}
		\eeq
		\begin{itemize}
			\item
			\textbf{Case (I)}: $\bX(t_n,1)=x_i-\frac{1}{\varepsilon}\Delta t\in  [x_{i-1},x_{i})$ and $\bX(t_n,-1)=x_i+\frac{1}{\varepsilon}\Delta t\in  (x_{i},x_{i+1}]$. \\
			Namely, when the time step is small, the foots are located in the adjacent cells connecting to $x_i$. Upwind discretizations are defined according to the sign in front of these derivatives, {corresponding to the real advective directions of this equation, not the directions of characteristics from the microscopic equation}, which are
			\begin{subequations}
				\label{eq46}
				\begin{align}
					&f_x\left(x_i-\frac{1}{\varepsilon}\Delta t,  1,  t_{n}\right)=\frac{f_{i-1}^{n}-f_{i-2}^n}{\Delta x},~~~\rho_x\left(x_i-\frac{1}{\varepsilon}\Delta t, t_{n}\right)=\frac{\rho_{i+1}^{n}-\rho_{i}^n}{\Delta x}, \\
					&f_x\left(x_i+\frac{1}{\varepsilon}\Delta t,  -1,  t_{n}\right)=\frac{f_{i+2}^{n}-f_{i+1}^n}{\Delta x},~~~\rho_x\left(x_i+\frac{1}{\varepsilon}\Delta t,  t_{n}\right)=\frac{\rho_{i}^{n}-\rho_{i-1}^n}{\Delta x}.
				\end{align}
			\end{subequations}
			\item
			\textbf{Case (II)}: $\bX(t_n,1)=x_{i}^{\ast}=x_i-\frac{1}{\varepsilon}\Delta t< x_{i-1}$ or $\bX(t_n,-1)=x_{i}^{\ast}=x_i+\frac{1}{\varepsilon}\Delta t> x_{i+1}$. \\ In this case, a large time step is used.
			If we assume $x_{i}^{\ast}=x_i-\frac{1}{\varepsilon}\Delta t \in [x_{i^{\ast}-1},x_{i^{\ast}})$, then
			\beq
			\label{eq47}
			f_x\left(x_i-\frac{1}{\varepsilon}\Delta t,  1,  t_{n}\right)=\frac{f_{i^{\ast}-1}^{n}-f_{i^{\ast}-2}^n}{\Delta x},~~~\rho_x\left(x_i-\frac{1}{\varepsilon}\Delta t, t_{n}\right)=\frac{\rho_{i^{\ast}+1}^{n}-\rho_{i^{\ast}}^n}{\Delta x}.
			\eeq
			Similarly if $x_{i}^{\ast}=x_i+\frac{1}{\varepsilon}\Delta t \in (x_{i^{\ast}-1},x_{i^{\ast}}]$,
			we have
			\beq
			\label{eq48}
			f_x\left(x_i+\frac{1}{\varepsilon}\Delta t,  -1,  t_{n}\right)=\frac{f_{i^{\ast}+1}^{n}-f_{i^{\ast}}^n}{\Delta x},~~~\rho_x\left(x_i+\frac{1}{\varepsilon}\Delta t,  t_{n}\right)=\frac{\rho_{i^{\ast}-1}^{n}-\rho_{i^{\ast}-2}^n}{\Delta x}.
			\eeq
		\end{itemize}
		The diffusion term $\rho_{xx}$ at $x_i$ is simply discretized by a second order central difference, which is
		\beq
		\label{rhoxx}
		\rho_{xx}^{n+1}\Big|_{x_i}\approx \frac{\rho_{i-1}^{n+1}-2\rho_{i}^{n+1}+\rho_{i+1}^{n+1}}{\Delta x^2}.
		\eeq
		With \eqref{eq46} or \eqref{eq47}-\eqref{eq48}, and \eqref{rhoxx}, substituting them into \eqref{time_dis11}, we complete a fully discrete first order scheme for $\rho$.
		\item {\bf Microscopic equation of $f$ in \eqref{time_dis12}}. \\
		Here a pure upwind flux according to the sign of $v$ is used:
		\begin{subequations}
			\label{eq410}
			\begin{align}
				&\frac{f_i^{n+1}-f_i^{n}}{\Delta t}+ \frac{1}{\varepsilon} \, \frac{f_i^{n+1}-f_{i-1}^{n+1}}{\Delta x}=\frac{1}{\varepsilon^2} (\rho_i^{n+1}-f_i^{n+1}), \quad \text{ for } v=1; \label{eq410a}\\
				&\frac{f_i^{n+1}-f_i^{n}}{\Delta t}- \frac{1}{\varepsilon} \, \frac{f_{i+1}^{n+1}-f_{i}^{n+1}}{\Delta x}=\frac{1}{\varepsilon^2} (\rho_i^{n+1}-f_i^{n+1}), \quad \text{ for } v=-1. \label{eq410b}
			\end{align}
		\end{subequations}
		{With the obtained macroscopic variable $\rho$, the microscopic variable $f$ in \eqref{eq410a} and \eqref{eq410b} can be updated from solving a linear system for each separate discrete velocity $v=1$ and $v=-1$, respectively. This is crucial, and is another main advantage of our approach, which is different from an iterative approach involving $\rho$ and $f$ together \cite{li2017implicit}. It is very efficient especially for one-group velocity and high-dimensional problems.} 
		
		\item {{\bf A correction of the density $\rho$ for mass conservation}. \\
			After obtaining $f^{n+1}$, we further update $\rho^{n+1}=\vave{f^{n+1}}$ to ensure a
			right mass conservation. This is important especially when the solution $f$ is not smooth.}
		
	\end{enumerate}

	\subsubsection{Second order finite difference scheme}
	\label{subsec_2nd}
	Following the descriptions in the previous subsection, the second order scheme is defined as follows.
	
	\begin{enumerate}
		\item {\bf Macroscopic equation of $\rho$ in \eqref{time_dis21}}.
		\begin{itemize}
			
			\item \textbf{Case (I)}: $\bX(t_n,1)=x_i-\frac{1}{\varepsilon}\Delta t\in  [x_{i-1},x_{i})$ and $\bX(t_n,-1)=x_i+\frac{1}{\varepsilon}\Delta t\in  (x_{i},x_{i+1}]$. \\
			If $x_i^{\ast}=x_i-\frac{1}{\varepsilon}\Delta t \in  [x_{i-1},x_{i}) $, we define $\xi=\frac{x_{i}-x_i^{\ast}}{\Delta x}$, then
			\begin{subequations}
				\label{eq411}
				\begin{align}
					& f_x\left(x_i-\frac{1}{\varepsilon}\Delta t,  1,  t_{n}\right)=\frac{1-2 \xi}{2\Delta x}f_{i-2}^{n}-\frac{2-2 \xi}{\Delta x}f_{i-1}^{n}+\frac{3-2 \xi}{2\Delta x}f_{i}^{n}, \\
					& \rho_x\left(x_i-\frac{1}{\varepsilon}\Delta t,  t_{n}\right)=-\frac{1+2 \xi}{2\Delta x}\rho_{i-1}^{n}+\frac{2 \xi}{\Delta x}\rho_{i}^{n}+\frac{1-2 \xi}{2\Delta x}\rho_{i+1}^{n}.
				\end{align}
			\end{subequations}
			If $x_i^{\ast}=x_i+\frac{1}{\varepsilon}\Delta t  \in  (x_{i},x_{i+1}]$, we define  $\eta=\frac{x_{i+1}-x_i^{\ast}}{\Delta x}$, then
			\begin{subequations}
				\label{eq412}
				\begin{align}
					& f_x\left(x_i+\frac{1}{\varepsilon}\Delta t,  -1,  t_{n}\right)=-\frac{1+2 \eta}{2\Delta x}f_{i}^{n}+\frac{2 \eta}{\Delta x}f_{i+1}^{n}+\frac{1-2 \eta}{2\Delta x}f_{i+2}^{n}, \\
					& \rho_x\left(x_i+\frac{1}{\varepsilon}\Delta t,  t_{n}\right)=\frac{1-2 \eta}{2\Delta x}\rho_{i-1}^{n}-\frac{2-2 \eta}{\Delta x}\rho_{i}^{n}+\frac{3-2 \eta}{2\Delta x}\rho_{i+1}^{n}.
				\end{align}
			\end{subequations}
			
			\item \textbf{Case (II)}: $\bX(t_n,1)=x_{i}^{\ast}=x_i-\frac{1}{\varepsilon}\Delta t< x_{i-1}$ or $\bX(t_n,-1)=x_{i}^{\ast}=x_i+\frac{1}{\varepsilon}\Delta t> x_{i+1}$.\\
			Assume $x_i^{\ast}=x_i-\frac{1}{\varepsilon}\Delta t \in [x_{i^{\ast}-1},x_{i^{\ast}})$, we define $\xi=\frac{x_{i^{\ast}}-x_i^{\ast}}{\Delta x}$, then
			\begin{subequations}
				\label{eq413}
				\begin{align}
					\label{eq413a}
					&f_x\left(x_{i}-\frac{1}{\varepsilon}\Delta t,  1,  t_{n}\right)=\frac{1-2 \xi}{2\Delta x}f_{i^{\ast}-2}^{n}-\frac{2-2 \xi}{\Delta x}f_{i^{\ast}-1}^{n}+\frac{3-2 \xi}{2\Delta x}f_{i^{\ast}}^{n},\\
					\label{eq413b}
					&\rho_x\left(x_i-\frac{1}{\varepsilon}\Delta t,  t_{n}\right)=-\frac{1+2 \xi}{2\Delta x}\rho_{i^{\ast}-1}^{n}+\frac{2 \xi}{\Delta x}\rho_{i^{\ast}}^{n}+\frac{1-2 \xi}{2\Delta x}\rho_{i^{\ast}+1}^{n}.
				\end{align}
			\end{subequations}
			Assume $x_i^{\ast}=x_i+\frac{1}{\varepsilon}\Delta t \in (x_{i^{\ast}-1},x_{i^{\ast}}]$, we define $\eta=\frac{x_{i^{\ast}}-x_i^{\ast}}{\Delta x}$, then
			\begin{subequations}
				\label{eq414}
				\begin{align}
					&f_x\left(x_i+\frac{1}{\varepsilon}\Delta t,  -1,  t_{n}\right)=-\frac{1+2 \eta}{2\Delta x}f_{i^{\ast}-1}^{n}+\frac{2 \eta}{\Delta x}f_{i^{\ast}}^{n}+\frac{(1-2 \eta)}{2\Delta x}f_{i^{\ast}+1}^{n},\\
					&\rho_x\left(x_i+\frac{1}{\varepsilon}\Delta t,  t_{n}\right)=\frac{1-2 \eta}{2\Delta x}\rho_{i^{\ast}-2}^{n}-\frac{2-2 \eta}{\Delta x}\rho_{i^{\ast}-1}^{n}+\frac{3-2 \eta}{2\Delta x}\rho_{i^{\ast}}^{n}.
				\end{align}
			\end{subequations}
		\end{itemize}
		$\rho_{xx}$ is still discretized by \eqref{rhoxx}. With \eqref{eq411}-\eqref{eq412} or \eqref{eq413}-\eqref{eq414}, and \eqref{rhoxx}, substituting them into \eqref{time_dis21}, we get a second order scheme for $\rho$.
		\item {\bf Microscopic equation of $f$ in \eqref{time_dis22}}. \\
		With the obtained macroscopic variable $\rho$, {the microscopic variable $f$ can also be obtained from solving the following linear systems for each separate discrete velocity,}
		\begin{subequations}
			\label{eq415}
			\begin{align}
				\label{eq415a}
				&\frac{f^{n-1}-4f^{n}+3f^{n+1}}{2\Delta t}+\f{1}{\eps} \,\frac{f_{i-2}^{n+1}-4f_{i-1}^{n+1}+3f_i^{n+1}}{2\Delta x}=\frac{1}{\varepsilon^2} (\rho^{n+1}-f^{n+1}), \,\text{ for } v=1; \\
				\label{eq415b}
				&\frac{f^{n-1}-4f^{n}+3f^{n+1}}{2\Delta t}-\f{1}{\eps} \, \frac{-3f_{i}^{n+1}+4f_{i+1}^{n+1}-f_{i+2}^{n+1}}{2\Delta x}=\frac{1}{\varepsilon^2} (\rho^{n+1}-f^{n+1}),\,\text{ for } v=-1.
			\end{align}
		\end{subequations}
		
		\item {{\bf A correction of the density $\rho$ for mass conservation}. \\
			Similarly, after obtaining $f^{n+1}$, we further update $\rho^{n+1}=\vave{f^{n+1}}$ to ensure a
			right mass conservation. }
	\end{enumerate}
	
	\begin{rem} 
		\label{rem42}
		For non-smooth solutions, linear reconstructions may develop some numerical oscillations. Here we take a simple slope limiter to control such numerical oscillations. Taking $f_x\left(x_{i}-\frac{1}{\varepsilon}\Delta t,1,t_{n}\right)$ as an example, we rewrite \eqref{eq413a} in a conserved form
		\[
		f_x\left(x_{i}-\frac{1}{\varepsilon}\Delta t,  1,  t_{n}\right)= \frac{{\hat f}_{i+1/2}-{\hat f}_{i-1/2}}{\Delta x},
		\]
		where the limited numerical flux is
		\[
		{\hat f}_{i+1/2}= f_{i^{\ast}}^{n} + \frac{1-2\xi}{2}\,\phi(r^n_{i^{\ast}})\,(f_{i^{\ast}}^{n}-f_{i^{\ast}-1}^{n}), \quad
		r^n_{i^{\ast}}=\frac{f_{i^{\ast}-1}^{n}-f_{i^{\ast}-2}^{n}}{f_{i^{\ast}}^{n}-f_{i^{\ast}-1}^{n}},
		\]
		and $\phi(r)$ is a slope limiting function. In our numerical tests, the van albada slope limiter with $\phi(r)=  \frac{r^2+r}{r^2+1}$ from \cite{van1997comparative} is used.
		
		For the implicit time discretization in \eqref{eq415a}, to avoid the nonlinearity from the nonlinear limiting function, the following procedure is used:
		\[
		f_x\left(x_{i},  1,  t_{n+1}\right)= \frac{{\tilde f}_{i+1/2}-{\tilde f}_{i-1/2}}{\Delta x},
		\]
		where
		\[
		{\tilde f}_{i+1/2}= f_{i}^{n+1} + \frac{1}{2}\,\phi(r^n_i)\,(f_i^{n+1}-f_{i-1}^{n+1}), \quad r^n_{i}=\frac{f_{i-1}^{n}-f_{i-2}^{n}}{f_i^{n}-f_{i-1}^{n}}, \quad \phi(r)=  \frac{r^2+r}{r^2+1}.
		\]
		\eqref{eq415b} can be treated analogously. In this way, we still only need to solve a linear system.
		
	\end{rem}
	
	\begin{rem}
		\label{rem43}
		{For spatial discretizations, numerical boundary treatments are very important, especially with a long characteristic tracking. 
			For a periodic boundary condition, we directly take the value by periodicity no matter how far the characteristic is.}
		
		{For a Dirichlet boundary condition, a close-loop numerical boundary treatment is first applied \cite{peng2021asymptotic,peng2020stability,jang2015high}. For the distribution function $f$, on the left boundary, an inflow boundary condition is assigned when $v>0$, i.e.  $f(x_L,v,t)=f_L(v,t)$ at $x_L$. Otherwise, an outflow boundary condition with an extrapolation is applied, that is $f(x_L,v,t)=f(x^+_L,v,t)$ for $v<0$, where $f(x^+_L,v,t)$ is a second order extrapolated value from the interior solution as we have used. The right boundary $x_R$ can be treated similarly. After obtaining $f$ at the boundary, the macroscopic density $\rho$ is given as
			\begin{align*}
				\rho_L(t)= \frac{1}{2}\left(\int_{-1}^0 f(x_L^+,v,t)dv + \int_{0}^1 f_L(v,t)dv  \right),  \,\,
				\rho_R(t)= \frac{1}{2}\left(\int_{-1}^0 f_R(v,t)dv + \int_{0}^1 f(x_R^-,v,t)dv  \right),
			\end{align*}
			with $\rho_L(t)$ and $\rho_R(t)$ corresponding to the left and right boundaries, respectively. However, except these adjacent boundary values, we also have to deal with boundary values from a characteristic tracking, e.g., \eqref{eq46}-\eqref{eq48} for a first order scheme, and \eqref{eq411}-\eqref{eq414} for a second order scheme. It is troublesome to deal with boundary values with long characteristics, especially when $\eps$ is small. However, noticing that the transport term in \eqref{eq45} and \eqref{time_dis12} has a coefficient $e^{-\mu\Delta t}/\eps$, the contribution of this term is negligible when $\eps$ is small enough. Numerically we simply take the values obtained from the close-loop treatment when the characteristic reaches out of the computational domain. }
	\end{rem}
	
	\subsection{Numerical stability analysis}
	\label{sec:Fourier}
	
	We now analyze the numerical stability of our proposed first and second order schemes by a Fourier analysis, {for the two-discrete-velocity telegraph equation.} The one-group velocity case can be similarly analyzed, but the algebraic calculation is more complicated. We assume a uniform mesh with a periodic boundary condition along $x$. 
	
	
	{Let
		\begin{equation*}
			p(x,t)=f(x,1,t), \quad q(x,t)=f(x,-1,t),
		\end{equation*}
		and $m$ be the integer which satisfies
		\begin{equation}\label{eqn:m}
			m<\frac{\Delta t}{\varepsilon\Delta x}\leq m+1.
		\end{equation}
		Consequently,
		\begin{equation*}
			x_{i+m}<x_{i}+\frac{\Delta t}{\varepsilon}\leq x_{i+m+1}, \quad x_{i-m-1}\leq x_{i}-\frac{\Delta t}{\varepsilon}< x_{i-m}.
		\end{equation*}
		Since we take a correction step $\rho^{n+1}=\vave{f^{n+1}}$, 
		to distinguish the numerical approximations of $\rho$ before and after correction, we denote it as $\sigma_j^{n+1}$ for $j=i-1,i,i+1$ before correction, in the first and second order schemes. With such notations, our first order scheme \eqref{eq45}-\eqref{eq410} in Section \ref{subsec_1st} can be rewritten in the following form:
		\begin{subequations}
			\label{1st}
			\begin{align}
				&\label{1sta} \frac{\sigma_{i}^{n+1}-\rho_{i}^{n}}{\Delta t}+\frac{e^{-\mu\Delta t}}{2\varepsilon}\Big[\frac{p_{i-m-1}^{n}-p_{i-m-2}^{n}}{\Delta x}-\frac{\rho_{i-m+1}^{n}-\rho_{i-m}^{n}}{\Delta x}-\frac{q_{i+m+2}^{n}-q_{i+m+1}^{n}}{\Delta x}  \\
				& \hspace{1.64cm}+\frac{\rho_{i+m}^{n}-\rho_{i+m-1}^{n}}{\Delta x}\Big]-(1-e^{-\mu\Delta t})\frac{\sigma_{i-1}^{n+1}-2\sigma_{i}^{n+1}+\sigma_{i+1}^{n+1}}{\Delta x^2}=0, \notag\\
				&\label{1stb} \frac{p_{i}^{n+1}-p_{i}^{n}}{\Delta t}+\frac{p_{i}^{n+1}-p_{i-1}^{n+1}}{\varepsilon\Delta x}=\frac{\sigma_{i}^{n+1}-p_{i}^{n+1}}{\varepsilon^2},\\
				&\label{1stc} \frac{q_{i}^{n+1}-q_{i}^{n}}{\Delta t}-\frac{q_{i+1}^{n+1}-q_{i}^{n+1}}{\varepsilon\Delta x}=\frac{\sigma_{i}^{n+1}-q_{i}^{n+1}}{\varepsilon^2},\\
				&\label{1std} \rho_{i}^{n+1}=\frac{p_{i}^{n+1}+q_{i}^{n+1}}{2}.
			\end{align}
		\end{subequations}
		Denoting $\mathbf{u}=(\sigma,p,q,\rho)^{\top}$ and taking the ansatz $\mathbf{u}_{i}^{n}=\mathbf{\hat{u}}^{n}e^{I\,\kappa x_{i}}$ with $I$ being the imaginary unit, substituting them into \eqref{1st}, we obtain
		\begin{equation}
			\mathbf{\hat{u}}^{n+1}=G^{(1)}(\Delta t, \Delta x, \varepsilon, \omega)\mathbf{\hat{u}}^{n},
		\end{equation}
		where $G^{(1)}=L^{(1)^{-1}}R^{(1)}$ is the amplification matrix with
		\begin{equation}\label{FAmatrix1}
			L^{(1)}=\left[\begin{array}{cccc}
				L_{11}^{(1)} & 0 & 0 & 0 \\
				-\frac{\Delta t}{\varepsilon^2} & L_{22}^{(1)} & 0 & 0 \\
				-\frac{\Delta t}{\varepsilon^2} & 0 & L_{33}^{(1)} & 0 \\
				0 & -\frac{1}{2} & -\frac{1}{2} & 1
			\end{array}\right],\quad
			R^{(1)}=\left[\begin{array}{cccc}
				0 & R_{12}^{(1)} & R_{13}^{(1)} & R_{14}^{(1)}\\
				0 & 1 & 0 & 0 \\
				0 & 0 & 1 & 0 \\
				0 & 0 & 0 & 0
			\end{array}\right],
		\end{equation}
		and
		\begin{eqnarray*}
			L_{11}^{(1)} &=& 1+\frac{\Delta t}{\Delta x^2}(1-e^{-\mu\Delta t})(2-2\cos(\omega)), \\
			L_{22}^{(1)} &=& 1+\frac{\Delta t}{\varepsilon\Delta x}(1-e^{-\omega I})+\frac{\Delta t}{\varepsilon^2}, \\
			L_{33}^{(1)} &=& 1+\frac{\Delta t}{\varepsilon\Delta x}(1-e^{\omega I})+\frac{\Delta t}{\varepsilon^2}, \\
			R_{12}^{(1)} &=& -\frac{e^{-\mu\Delta t}\Delta t}{2\varepsilon\Delta x}(e^{-(m+1)\omega I}-e^{-(m+2)\omega I}), \\
			R_{13}^{(1)} &=& -\frac{e^{-\mu\Delta t}\Delta t}{2\varepsilon\Delta x}(e^{(m+1)\omega I}-e^{(m+2)\omega I}), \\
			R_{14}^{(1)} &=& 1+\frac{e^{-\mu\Delta t}\Delta t}{\varepsilon\Delta x}(\cos((m-1)\omega)-\cos(m\omega)).
		\end{eqnarray*}
		Here and below, $\omega=\kappa\Delta x$ is used to denote the discrete wave number.}
	
	{Similarly, our second order scheme \eqref{eq411}-\eqref{eq415} in Section \ref{subsec_2nd} can be rewritten as
		\begin{subequations}
			\label{2nd}
			\begin{align}
				&\label{2nda} \frac{\rho_{i}^{n-1}-4\rho_{i}^{n}+3\sigma_{i}^{n+1}}{2\Delta t}+\frac{e^{-\mu\Delta t}}{2\varepsilon}
				\Big[\left(\frac{1-2\xi}{2\Delta x}p_{i-m-2}^{n}-\frac{2-2\xi}{\Delta x}p_{i-m-1}^{n}+\frac{3-2\xi}{2\Delta x}p_{i-m}^{n}\right)  \\
				&\hspace{1.6cm}-\left(-\frac{1+2\xi}{2\Delta x}\rho_{i-m-1}^{n}+\frac{2\xi}{\Delta x}\rho_{i-m}^{n}+\frac{1-2\xi}{2\Delta x}\rho_{i-m+1}^{n}\right) \notag \\
				&\hspace{1.6cm}-\left(-\frac{3-2\xi}{2\Delta x}q_{i+m}^{n}+\frac{2-2\xi}{\Delta x}q_{i+m+1}^{n}-\frac{1-2\xi}{2\Delta x}q_{i+m+2}^{n}\right) \notag \\
				&\hspace{1.6cm}+\left(-\frac{1-2\xi}{2\Delta x}\rho_{i+m-1}^{n}-\frac{2\xi}{\Delta x}\rho_{i+m}^{n}+\frac{1+2\xi}{2\Delta x}\rho_{i+m+1}^{n}\right)\Big] \notag \\
				&\hspace{1.6cm}-(1-e^{-\mu\Delta t})\frac{\sigma_{i-1}^{n+1}-2\sigma_{i}^{n+1}+\sigma_{i+1}^{n+1}}{\Delta x^2}=0, \notag \\
				&\label{2ndb} \frac{p_{i}^{n-1}-4p_{i}^{n}+3p_{i}^{n+1}}{2\Delta t}+\frac{p_{i-2}^{n+1}-4p_{i-1}^{n+1}+3p_{i}^{n+1}}{2\eps\Delta x}=\frac{\sigma_{i}^{n+1}-p_{i}^{n+1}}{\varepsilon^2},\\
				&\label{2ndc} \frac{q_{i}^{n-1}-4q_{i}^{n}+3q_{i}^{n+1}}{2\Delta t}-\frac{-3q_{i}^{n+1}+4q_{i+1}^{n+1}-q_{i+2}^{n+1}}{2\eps\Delta x}=\frac{\sigma_{i}^{n+1}-q_{i}^{n+1}}{\varepsilon^2},\\
				&\label{2ndd} \rho_{i}^{n+1}=\frac{p_{i}^{n+1}+q_{i}^{n+1}}{2}.
			\end{align}
		\end{subequations}
		Here $\xi=\frac{\Delta t}{\varepsilon \Delta x}-m$, which is consistent with the definition in Section \ref{subsec_2nd} and (\ref{eqn:m}). Notice that $\eta=1-\xi$ in (\ref{eq412}) and (\ref{eq414}) in this scenario.}
	
	{The second order scheme is a two-step scheme. To determine the amplification matrix, we introduce three auxiliary variables $a_{i}^{n+1}=\rho_{i}^{n}, b_{i}^{n+1}=p_{i}^{n}, c_{i}^{n+1}=q_{i}^{n}$ and let $\mathbf{u}=(\sigma,p,q,a,b,c,\rho)^{\top}$. Using the same ansatz as above, similar to the first order scheme, we get
		\begin{equation}
			\mathbf{\hat{u}}^{n+1}=G^{(2)}(\Delta t, \Delta x, \varepsilon, \omega)\mathbf{\hat{u}}^{n},
		\end{equation}
		where the amplification matrix $G^{(2)}=L^{(2)^{-1}}R^{(2)}$ and
		\begin{equation}\label{FAmatrix2L}
			L^{(2)}=\left[\begin{array}{ccccccc}
				L_{11}^{(2)} & 0 & 0 & 0 & 0 & 0 & 0 \\
				-\frac{2\Delta t}{\varepsilon^2} & L_{22}^{(2)} & 0 & 0 & 0 & 0 & 0 \\
				-\frac{2\Delta t}{\varepsilon^2} & 0 & L_{33}^{(2)} & 0 & 0 & 0 & 0 \\
				0 & 0 & 0 & 1 & 0 & 0 & 0 \\
				0 & 0 & 0 & 0 & 1 & 0 & 0 \\
				0 & 0 & 0 & 0 & 0 & 1 & 0 \\
				0 & -\frac{1}{2} & -\frac{1}{2} & 0 & 0 & 0 & 1
			\end{array}\right], \quad
			R^{(2)}=\left[\begin{array}{ccccccc}
				0 & R_{12}^{(2)} & R_{13}^{(2)} & -1 & 0 & 0 & R_{17}^{(2)}\\
				0 & 4 & 0 & 0 & -1 & 0 & 0 \\
				0 & 0 & 4 & 0 & 0 & -1 & 0 \\
				0 & 0 & 0 & 0 & 0 & 0 & 1 \\
				0 & 1 & 0 & 0 & 0 & 0 & 0 \\
				0 & 0 & 1 & 0 & 0 & 0 & 0 \\
				0 & 0 & 0 & 0 & 0 & 0 & 0
			\end{array}\right],
		\end{equation}
		\begin{eqnarray*}
			L_{11}^{(2)} &=& 3+\frac{2\Delta t}{\Delta x^2}(1-e^{-\mu\Delta t})(2-2\cos(\omega)), \\
			L_{22}^{(2)} &=& 3+\frac{\Delta t}{\varepsilon\Delta x}(3-4e^{-\omega I}+e^{-2\omega I})+\frac{2\Delta t}{\varepsilon^2}, \\
			L_{33}^{(2)} &=& 3+\frac{\Delta t}{\varepsilon\Delta x}(3-4e^{\omega I}+e^{2\omega I})+\frac{2\Delta t}{\varepsilon^2}, \\
			R_{12}^{(2)} &=& -\frac{e^{-\mu\Delta t}\Delta t}{2\varepsilon\Delta x}\left[(3-2\xi)e^{-m\omega I}-(4-4\xi)e^{-(m+1)\omega I}+(1-2\xi)e^{-(m+2)\omega I}\right], \\
			R_{13}^{(2)} &=& -\frac{e^{-\mu\Delta t}\Delta t}{2\varepsilon\Delta x}\left[(3-2\xi)e^{m\omega I}-(4-4\xi)e^{(m+1)\omega I}+(1-2\xi)e^{(m+2)\omega I}\right], \\
			R_{17}^{(2)} &=& 4+\frac{e^{-\mu\Delta t}\Delta t}{\varepsilon\Delta x}[(1-2\xi)\cos((m-1)\omega)+4\xi\cos(m\omega)-(1+2\xi)\cos((m+1)\omega)].
		\end{eqnarray*}
		The stability is essentially determined by the eigenvalues of the amplification matrix.} We use the following principle to study the numerical stability \cite{peng2020stability}:
	
	\textbf{Principle for Numerical Stability:} Given the values of $\Delta x, \Delta t$ and  $\varepsilon$, let $\lambda_j(\omega), j=1,2,\dots,Q$ ($Q=4$ for the first order scheme and $Q=7$ for the second order scheme) be the eigenvalues of the amplification matrix $G$. Our scheme is ``stable'' if for all $\omega\in[-\pi,\pi]$, it satisfies either
	\begin{enumerate}
		\item $\max_{1\leq j\leq Q} |\lambda_j(\omega)|<1$, or
		\item $\max_{1\leq j\leq Q} |\lambda_j(\omega)|=1$, and $G$ is diagonalizable.
	\end{enumerate}
	In our numerical stability analysis, we test several different settings covering $\Delta x=10^{-j}\,(1\leq j\leq 4)$, $\Delta t=10^{k}\Delta x\, (-3\leq k\leq 3)$ and $\varepsilon=10^{l}\,(-10\leq l\leq 5)$. The discrete wave number $\omega$ is taken from $[-\pi,\pi]$ with $500$ samples uniformly distributed. Based on such a stability principle, we find that our first and second order schemes are uniformly unconditionally stable, which is consistent with what we will observe numerically.
	
	{\subsection{Formal asymptotic analysis}
		\label{sec:ap}
		We now formally prove that our first and second order fully discrete schemes are AP. We take the two-discrete-velocity telegraph equation as an example, and assume all solutions are of $\mathcal{O}(1)$. It can be straightforwardly extended to the one-group velocity and other collision operator models.}
	
	{We follow the notations in the previous subsection. For the first order scheme \eqref{1st}, as $\eps\rightarrow 0$, due to $e^{-\mu\Delta t}\rightarrow 0$ and $e^{-\mu\Delta t}/\eps\rightarrow 0$, from \eqref{1st}, at $x_i$ we have 
		\begin{equation}
			\label{heatns}
			\frac{\sigma_i^{n+1}-\rho_i^n}{\Delta t}-\frac{\sigma_{i-1}^{n+1}-2\sigma_{i}^{n+1}+\sigma_{i+1}^{n+1}}{\Delta x^2}=0.
		\end{equation}
		From \eqref{1stb} and \eqref{1stc}, it is easy to find that $p^{n+1}_i=\sigma^{n+1}_i$ and $q^{n+1}_i=\sigma^{n+1}_i$. The correction step \eqref{1std} yields $\rho^{n+1}_i=\frac{p^{n+1}_i+q^{n+1}_i}{2}=\sigma^{n+1}_i$, so that we have $\rho^{n+1}_i=f^{n+1}_i=\sigma^{n+1}_i$, where $\sigma^{n+1}_i$ is solved from \eqref{heatns}, which is a consistent discretization for the limiting equation \eqref{apprxmodellimit}.}
	
	{For the second order scheme \eqref{2nd}, it can be similarly analyzed. From \eqref{2nda}, as $\eps\rightarrow 0$, we have
		\begin{equation}
			\label{heatns2}
			\frac{\rho_{i}^{n-1}-4\rho_{i}^{n}+3\sigma_{i}^{n+1}}{2\Delta t}-\frac{\sigma_{i-1}^{n+1}-2\sigma_{i}^{n+1}+\sigma_{i+1}^{n+1}}{\Delta x^2}=0,
		\end{equation}
		while \eqref{2ndb} and \eqref{2ndc} still give $p^{n+1}_i=\sigma^{n+1}_i$ and $q^{n+1}_i=\sigma^{n+1}_i$, so that \eqref{2ndd} yields $\rho^{n+1}_i=\frac{p^{n+1}_i+q^{n+1}_i}{2}=\sigma^{n+1}_i$, with $\sigma^{n+1}_i$ from \eqref{heatns2} as a second order consistent discretization to \eqref{apprxmodellimit}. So both our first and second order schemes are AP.}
	
	\section{Numerical tests}
	\label{sec_numer}
	\setcounter{equation}{0}
	\setcounter{figure}{0}
	\setcounter{table}{0}
	
	In this section, we will verify the convergence orders of our proposed first and second order schemes, and demonstrate their good performance with large time step conditions for both smooth and discontinuous solutions. We take $N$ uniform grid points along each direction in space. In all figures, the numerical solutions of the macroscopic density $\rho$ are presented. A reference time step is taken to be
	\[
	\Delta t= \text{CFL}\, \Delta x,
	\]
	where CFL is a referred CFL number, which can be large enough if only stability is concerned.
	For the  second order scheme, if discontinuous solutions appear, the limiter described in Remark 4.2 will be applied. We consider 1D in space discrete-velocity models in Section \ref{sec_numer_dv}, 1D in space one-group velocity models in Section \ref{sec_numer_oneg}, followed by two-dimensional (2D) space problems in Section \ref{sec_numer_2d}.
	
	\subsection{Discrete-velocity transport equations}
	\label{sec_numer_dv}
	\subsubsection{Telegraph equation}
	\label{tele}
	First, we consider the telegraph equation with an exact solution
	\beq
	\label{eq51}
	f(x,v,t)= \frac{1 }{r}e^{rt}\sin(x)+v\varepsilon e^{rt}\cos(x), \,\, \rho(x,t)=\frac{1 }{r}e^{rt}\sin(x),\,\,r=\frac{-2}{1+\sqrt{(1-4\varepsilon^2)}}.
	\eeq
	We consider a periodic boundary condition on the domain $[-\pi,\pi]$, and compute the solution up to a final time $T=1$ for several different $\varepsilon$'s, $\varepsilon =0.5, 0.1, 10^{-2}, 10^{-6}$. In Table \ref{table1} and Table \ref{table2}, we show the $L^1$ and $L^\infty$ errors and corresponding orders with a large time step $\Delta t=3\Delta x$. In Table \ref{table_TG_time1} and Table \ref{table_TG_time2}, we show the temporal accuracy and orders on a   fixed fine mesh with $N=5000$. For large (rarefied regime) and small (diffusive regime) $\eps$'s, we can observe the corresponding first and second orders of accuracy both in space and in time as expected. For intermediate $\eps$'s, such as $\eps=0.1$, due to the modeling error as shown in Theorem \ref{thm1}, order reductions can be observed when $\eps^2\approx\Delta t$, especially for the second order scheme (Table \ref{table2} and Table \ref{table_TG_time2}). We have also tried larger time steps with $\text{CFL} = 10, 20$. Similar observations can be found, which verify that our schemes work well for very large time steps.
	
	\begin{table}[!ht]
		\scriptsize
		\centering
		\caption{$L^{\infty}$ and $L^{1}$ errors and convergence orders in space for  $\rho$ and $f(x,1,T)$. First order scheme for telegraph equation.  \label{table1}}
		\begin{tabular}{|c|c|c|c|c|c|c|c|c|c|}
			\hline	
			$\eps$ &N  & $L^{\infty}$ error of $\rho$ & Order & $L^{\infty}$ error of $f$ & Order  & $L^{1}$ error of $\rho$ & Order & $L^{1}$ error of $f$ & Order        \\ \hline
			\multirow{6}{*}{$0.5$}         	
			& 40 & 7.35E-2 & -- &  7.33E-2& --  & 4.68E-2 & -- &  4.72E-2& --                        \\ \cline{2-10}
			& 80 & 2.89E-2 & 1.35 &  3.13E-2& 1.23  & 1.84E-2 & 1.35 &  2.01E-2& 1.23              \\ \cline{2-10}
			& 160 & 1.02E-2 & 1.50 &  1.24E-2& 1.33 & 6.50E-3 & 1.50 &  7.92E-3& 1.34       \\ \cline{2-10}
			& 320 & 3.55E-3 & 1.52 &  4.92E-3& 1.33   & 2.26E-3 & 1.52 &  3.14E-3& 1.34     \\ \cline{2-10}
			& 640 & 1.33E-3 & 1.41   &  2.09E-3& 1.24   & 8.49E-4 & 1.41 &  1.33E-3& 1.23       \\  \hline
			\multirow{6}{*}{$0.1$}         	
			& 40 & 7.24E-2 & -- &  6.92E-2& -- & 4.60E-2 & -- &  4.44E-2& --                       \\ \cline{2-10}
			& 80 & 4.12E-2 & 0.81 &  3.98E-2& 0.80 & 2.62E-2 & 0.81 &  2.54E-2& 0.81              \\ \cline{2-10}
			& 160 & 2.36E-2 & 0.81 &  2.29E-2& 0.80 & 1.50E-2 & 0.81 &  1.46E-2& 0.80       \\ \cline{2-10}
			& 320 & 1.41E-2 & 0.74 &  1.38E-2& 0.73 & 8.99E-3 & 0.74 &  8.77E-3& 0.73      \\ \cline{2-10}
			& 640 & 7.27E-3 & 0.96   &  7.12E-3& 0.95  & 4.63E-3 & 0.96   &  4.54E-3& 0.95       \\  \hline
			\multirow{5}{*}{$10^{-2}$}
			& 40 & 7.30E-2 & -- &  7.26E-2& --  & 4.64E-2 & -- &  4.62E-2& --                       \\ \cline{2-10}
			& 80 & 3.95E-2 & 0.88 &  3.94E-2& 0.88 & 2.52E-2 & 0.88 &  2.51E-2& 0.88              \\ \cline{2-10}
			& 160 & 2.07E-2 & 0.94 &  2.06E-2& 0.94 & 1.32E-2 & 0.94 &  1.31E-2& 0.94       \\ \cline{2-10}
			& 320 & 1.06E-2 & 0.96 &  1.06E-2& 0.96  & 6.74E-3 & 0.96 &  6.72E-3& 0.96     \\ \cline{2-10}
			& 640 & 5.38E-3 & 0.98   &  5.37E-3& 0.98 & 3.43E-3 & 0.98   &  3.42E-3& 0.98     \\  \hline
			\multirow{5}{*}{$10^{-6}$}
			& 40 & 7.29E-2 & -- &  7.26E-2& --  & 4.64E-2 & -- &  4.62E-2& --                       \\ \cline{2-10}
			& 80 & 3.95E-2 & 0.88 &  3.94E-2& 0.88 & 2.52E-2 & 0.88 &  2.51E-2& 0.88              \\ \cline{2-10}
			& 160 & 2.07E-2 & 0.94 &  2.06E-2& 0.94 & 1.32E-2 & 0.94 &  1.31E-2& 0.94       \\ \cline{2-10}
			& 320 & 1.06E-2 & 0.97 &  1.06E-2& 0.96  & 6.72E-3 & 0.97 &  6.72E-3& 0.97     \\ \cline{2-10}
			& 640 & 5.35E-3 & 0.98   &  5.35E-3& 0.98 & 3.41E-3 & 0.98   &  3.41E-3& 0.98     \\  \hline
			
		\end{tabular}
	\end{table}
	
	\begin{table}[!ht]
		\scriptsize
		\centering
		\caption{$L^{\infty}$ and $L^{1}$  errors and convergence orders in space for $\rho$ and $f(x,1,T)$. Second order scheme for telegraph equation. \label{table2}}
		\begin{tabular}{|c|c|c|c|c|c|c|c|c|c|}
			\hline	
			$\eps$ &N  & $L^{\infty}$ error of $\rho$ & Order & $L^{\infty}$ error of $f$ & Order  & $L^{1}$ error of $\rho$ & Order & $L^{1}$ error of $f$ & Order        \\ \hline
			
			\multirow{6}{*}{$ 0.5$}         	
			
			& 40 & 1.27E-1 & -- &  1.17E-1& --  & 8.06E-2 & -- &  7.63E-2& --                        \\ \cline{2-10}
			& 80 & 5.37E-2 & 1.24 &  4.85E-2& 1.27  & 3.42E-2 & 1.23 &  3.12E-2& 1.29              \\ \cline{2-10}
			& 160 & 1.64E-2 & 1.71 &  1.49E-2& 1.70 & 1.04E-2 & 1.71 &  9.53E-3& 1.71      \\ \cline{2-10}
			& 320 & 4.45E-3 & 1.88 &  4.02E-3& 1.89   & 2.84E-3 & 1.88 &  2.56E-3& 1.89     \\ \cline{2-10}
			& 640 & 1.14E-3 & 1.97   &  1.03E-3& 1.97   & 7.24E-4 & 1.97 &  6.54E-4& 1.97       \\  \hline
			\multirow{8}{*}{$ 0.1$}         	
			& 40 & 5.02E-2 & -- &  5.03E-2& --  & 3.19E-2 & -- &  3.22E-2& --                        \\ \cline{2-10}
			& 80 & 1.56E-2 & 1.69 &  1.56E-2& 1.69  & 9.91E-3 & 1.69 &  1.00E-2& 1.69              \\ \cline{2-10}
			& 160 & 6.33E-3 & 1.30 & 6.32E-3& 1.30 & 4.03E-3 &1.30 &  4.03E-3& 1.30      \\ \cline{2-10}
			& 320 & 4.42E-3 & 0.52 &  4.40E-3& 0.52   & 2.82E-3 & 0.52 & 2.81E-3& 0.52     \\ \cline{2-10}
			& 640 & 3.93E-3 & 0.17   &  3.91E-3& 0.17   & 2.50E-3 & 0.17 &  2.49E-3& 0.17  \\\cline{2-10}
			& 1280 & 2.57E-3 & 0.61   &  2.56E-3& 0.61   & 1.64E-3 & 0.61 &  1.63E-3& 0.61  \\\cline{2-10}
			& 2560 & 1.12E-3 & 1.20   &  1.11E-3& 1.20   & 7.13E-4 & 1.20 &  7.09E-4& 1.20  \\\cline{2-10}
			& 5120 & 3.75E-4 & 1.58   &  3.73E-3& 1.58   & 2.38E-4 & 1.58 &  2.37E-4& 1.58  \\\hline
			\multirow{5}{*}{$ 10^{-2}$}
			& 40 & 4.70E-2 & -- &  4.70E-2& --  & 2.99E-2 & -- &  2.99E-2& --                       \\ \cline{2-10}
			& 80 & 1.20E-2 & 1.97 &  1.20E-2& 1.97 & 7.66E-3 & 1.97 & 7.66E-3 & 1.97            \\ \cline{2-10}
			& 160 & 2.70E-3 & 2.15 &  2.70E-3& 2.15 & 1.72E-3 & 2.15 & 1.72E-3 & 2.15       \\ \cline{2-10}
			& 320 & 6.64E-4 & 2.03 &  6.64E-4& 2.03  & 4.23E-4 & 2.03 & 4.23E-4 & 2.03     \\ \cline{2-10}
			& 640 & 1.82E-4 & 1.87   &  1.82E-4& 1.87 & 1.16E-4 & 1.87 & 1.16E-4 & 1.87     \\  \hline
			\multirow{5}{*}{$ 10^{-6}$}
			& 40 & 4.70E-2 & -- &  4.70E-2& --  & 2.98E-2 & -- &  2.98E-2& --                       \\ \cline{2-10}
			& 80 & 1.20E-2 & 1.97 &  1.20E-2& 1.97 & 7.63E-3 & 1.97 & 7.63E-3 & 1.97            \\ \cline{2-10}
			& 160 & 2.67E-3 & 2.17 &  2.67E-3& 2.15 & 1.70E-3 & 2.17 &1.70E-3 & 2.17       \\ \cline{2-10}
			& 320 & 6.27E-4 & 2.09 &  6.27E-4& 2.03  & 4.00E-4 & 2.09 & 4.00E-4 & 2.09    \\ \cline{2-10}
			& 640 & 1.46E-4 & 2.11   &  1.46E-4& 2.11 & 9.27E-5 & 2.11 & 9.27E-5 & 2.11     \\  \hline
			
		\end{tabular}
	\end{table}

	\begin{table}[htbp]
		\scriptsize
		\centering
		\caption{ $L^{\infty}$ and $L^{1}$  errors and convergence orders in time for $\rho$ and $f(x,1,T)$. First order scheme  for telegraph equation with a fixed mesh $N=5000$.  \label{table_TG_time1} }
		\begin{tabular}{|c|c|c|c|c|c|c|c|c|c|}
			\hline	
			$\eps$  & $\Delta t$  & $L^{\infty}$ error of $\rho$ & Order & $L^{\infty}$ error of $f$ & Order & $L^1$ error of $\rho$ & Order & $L^1$ error of $f$ & Order  \\ \hline
			\multirow{6}{*}{$ 0.5$}         	
			
			& T/8 & 3.61E-2 & -- &  3.92E-2& --    & 2.30E-2 & -- &  2.50E-2& --     \\ \cline{2-10}
			& T/16 & 1.40E-2 & 1.37 & 1.64E-2& 1.26  &  8.93E-3& 1.36 &  1.05E-2& 1.25   \\ \cline{2-10}
			& T/32 & 5.29E-3 & 1.40 &  6.79E-3 & 1.27  &  3.37E-3& 1.41 &  4.33E-3& 1.28   \\ \cline{2-10}
			& T/64 & 1.89E-3 & 1.48 &  2.75E-3 & 1.30   &  1.20E-3& 1.49 &  1.75E-3& 1.31   \\ \cline{2-10}
			& T/128 & 4.88E-4 & 1.95 &  9.84E-4 & 1.48   &  3.11E-4& 1.95 &  6.27E-4& 1.48   \\ \hline
			\multirow{6}{*}{$ 0.1$}
			& T/8 & 2.53E-2 & -- &  2.54E-2& --    & 1.61E-2 & -- &  1.62E-2& --     \\ \cline{2-10}
			& T/16 & 1.49E-2 & 0.76 & 1.48E-2& 0.78  &  9.46E-3& 0.77 &  9.47E-3& 0.77   \\ \cline{2-10}
			& T/32 & 9.09E-3 & 0.71 &  9.08E-3 & 0.70  &  5.79E-3& 0.71 &  5.78E-3& 0.71   \\ \cline{2-10}
			& T/64 & 3.68E-3 & 1.30 &  3.67E-3 & 1.31   &  2.34E-3& 1.31 &  2.33E-3& 1.31   \\ \cline{2-10}
			& T/128 & 1.52E-3 & 1.28 &  1.52E-3 & 1.27   &  9.66E-4& 1.28 &  9.67E-4& 1.27    \\ \hline
			\multirow{6}{*}{$ 10^{-2}$}
			& T/8 & 2.19E-2 & -- &  2.19E-2  & --    & 1.69E-2 & -- &  1.69E-2& --     \\ \cline{2-10}
			& T/16 & 1.12E-2 & 0.96 &  1.12E-2 & 0.96  & 7.16E-3& 0.96 & 7.16E-3& 0.96   \\ \cline{2-10}
			& T/32 & 5.71E-3 & 0.98 &  5.71E-3 & 0.98  & 3.64E-3& 0.98 &  3.64E-3& 0.98   \\ \cline{2-10}
			& T/64 & 2.89E-3 & 0.98 &  2.89E-3 & 0.98   & 1.84E-3& 0.98 &  1.84E-3& 0.98 \\ \cline{2-10}
			& T/128 & 1.47E-3 & 0.98 & 1.47E-3 & 0.98   &  9.35E-4& 0.98 & 9.34E-4& 0.98   \\ \hline
			\multirow{6}{*}{$ 10^{-6}$}
			& T/8 & 2.19E-2 & -- &  2.19E-2  & --    & 1.39E-2 & -- &  1.39E-2& --     \\ \cline{2-10}
			& T/16 & 1.12E-2 & 0.96 &  1.12E-2 & 0.96  & 7.13E-3& 0.96 & 7.13E-3& 0.96   \\ \cline{2-10}
			& T/32 & 5.67E-3 & 0.98 &  5.67E-3 & 0.98  & 3.61E-3& 0.98 &  3.61E-3& 0.98   \\ \cline{2-10}
			& T/64 & 2.86E-3 & 0.99 &  2.86E-3 & 0.99   & 1.82E-3& 0.99 &  1.82E-3& 0.99 \\ \cline{2-10}
			& T/128 & 1.43E-3 & 1.00 & 1.43E-3 & 1.00   &  9.12E-4& 1.00 & 9.12E-4& 1.00   \\ \hline
		\end{tabular}
	\end{table}
	
	\begin{table}[htbp]
		\scriptsize
		\centering
		\caption{$L^{\infty}$ and $L^{1}$  errors and convergence orders in time for $\rho$ and $f(x,1,T)$. Second order scheme   for telegraph equation with a fixed mesh $N=5000$.   \label{table_TG_time2} }
		\begin{tabular}{|c|c|c|c|c|c|c|c|c|c|}
			\hline	
			$\eps$  & $\Delta t$  & $L^{\infty}$ error of $\rho$ & Order & $L^{\infty}$ error of $f$ & Order & $L^1$ error of $\rho$ & Order & $L^1$ error of $f$ & Order  \\ \hline
			\multirow{6}{*}{$0.5$}
			& T/8 & 1.77E-2 & -- &  1.64E-2& --    & 1.13E-2 & -- &  1.04E-2& --     \\ \cline{2-10}
			& T/16 & 4.86E-3 & 1.86 & 4.47E-3& 1.87  &  3.09E-3& 1.87 &  2.84E-3& 1.87   \\ \cline{2-10}
			& T/32 & 1.26E-3 & 1.95 &  1.15E-3 & 1.96  &  8.05E-4& 1.94 & 7.35E-4& 1.95   \\ \cline{2-10}
			& T/64 & 3.21E-4 & 1.97 &  2.92E-4 & 1.98   &  2.05E-4& 1.97 &  1.86E-4& 1.98   \\ \cline{2-10}
			& T/128 & 8.11E-5 & 1.99 &  7.34E-5 & 1.99   &  5.16E-5& 1.99 &  4.67E-5& 1.99    \\ \hline
			\multirow{8}{*}{$0.1$}         	
			
			& T/2 &  5.20E-2 & -- &  5.22E-2& --    &  3.31E-2 & -- &  3.32E-2& --   \\ \cline{2-10}
			& T/4 &  1.59E-2 & 1.71 &  1.59E-2& 1.72    &  1.01E-2 & 1.71 &  1.01E-2& 1.72   \\ \cline{2-10}
			& T/8 & 6.36E-3 & 1.32   &  6.35E-3& 1.32     & 4.05E-3 &1.32    &  4.04E-3& 1.32      \\ \cline{2-10}
			& T/16 & 4.41E-3 & 0.52 & 4.39E-3& 0.53  &  2.81E-3& 0.53 &  2.80E-3& 0.53   \\ \cline{2-10}
			
			& T/32 & 3.99E-3 & 0.15 &  3.97E-3 & 0.15  &  2.54E-3& 0.15 &  2.53E-3& 0.15   \\ \cline{2-10}
			
			& T/64 & 2.71E-3 & 0.56 &  2.70E-3 & 0.56   &  1.73E-3& 0.55 &  1.72E-3& 0.56   \\ \cline{2-10}
			& T/128 & 1.22E-3 & 1.15 &  1.21E-3 & 1.16   &  7.75E-4& 1.16 &  7.74E-4& 1.16   \\ \cline{2-10}
			& T/256 &  4.14E-4 & 1.56  &  4.12E-4 & 1.55    &  2.64E-4& 1.55  &  2.63E-4& 1.56  \\ \cline{2-10}
			& T/512 &  1.21E-4 & 1.77  &  1.21E-4 & 1.77    &  7.71E-5& 1.78  &  7.68E-5& 1.78 \\ \hline
			\multirow{7}{*}{$10^{-2}$}         	
			& T/2 &  4.88E-2 & -- &  4.88E-2& --    &  3.11E-2 & -- &  3.11E-2& --   \\ \cline{2-10}
			& T/4 &  1.23E-2 & 1.98 & 1.23E-2 & 1.98    &7.85E-3 & 1.98 & 7.85E-3 & 1.98    \\ \cline{2-10}
			& T/8 & 2.72E-3 & 2.18   &  2.72E-3 & 2.18  &1.73E-3 &2.18    &  1.73E-3 &2.18     \\ \cline{2-10}
			& T/16 & 6.64E-4 & 2.04 & 6.64E-4 & 2.04  &  4.23E-4& 2.04 & 4.23E-4& 2.04  \\ \cline{2-10}
			& T/32 & 1.90E-4 & 1.81 &  1.90E-4 & 1.81  &  1.21E-4& 1.81 &  1.21E-4& 1.81   \\ \cline{2-10}
			& T/64 & 7.46E-5 & 1.35 &  7.46E-5 & 1.35   &  4.75E-5& 1.35 & 4.75E-5& 1.35  \\ \cline{2-10}
			& T/128 & 4.62E-5 & 0.69 &  4.62E-5 & 0.69   &  2.94E-5& 0.69 &  2.94E-5& 0.69 \\  \hline
			\multirow{5}{*}{$10^{-6}$}
			& T/8 & 2.69E-3 & -- &   2.69E-3& --    & 1.71E-3 & -- &  1.71E-3& --     \\ \cline{2-10}
			& T/16 & 6.28E-4 & 2.09 & 6.28E-4 & 2.09  & 4.00E-4& 2.09 &  4.00E-4& 2.09   \\ \cline{2-10}
			& T/32 & 1.53E-4 & 2.04 &  1.53E-4 & 2.04 &  9.74E-5& 2.04 &  9.74E-5& 2.04   \\ \cline{2-10}
			& T/64 & 3.78E-5 & 2.01 &  3.78E-5 & 2.01   &  2.41E-5& 2.01 &   2.41E-5& 2.01   \\ \cline{2-10}
			& T/128 & 9.45E-6 & 2.00 &  9.45E-6 & 2.00  &  6.02E-6& 2.00 & 6.02E-6& 2.00   \\ \hline
		\end{tabular}
	\end{table}
	
	Then we consider a Riemann problem with the following initial data
	\begin{equation}
		\label{tele_dis}
		\left\{\begin{array}{llll}
			\rho(x,0)=2.0, & f(x,v,0)=2.0,   & -1<x<0, \\
			\rho(x,0)=1.0, & f(x,v,0)=1.0,   & 0<x<1.
		\end{array}\right.
	\end{equation}
	Inflow and outflow boundary conditions are taken.
	
\begin{figure}[!ht]
		\centering
		\subfigure[]
		{
			\begin{minipage}{7cm}
				\centering
				\includegraphics[scale=0.36]{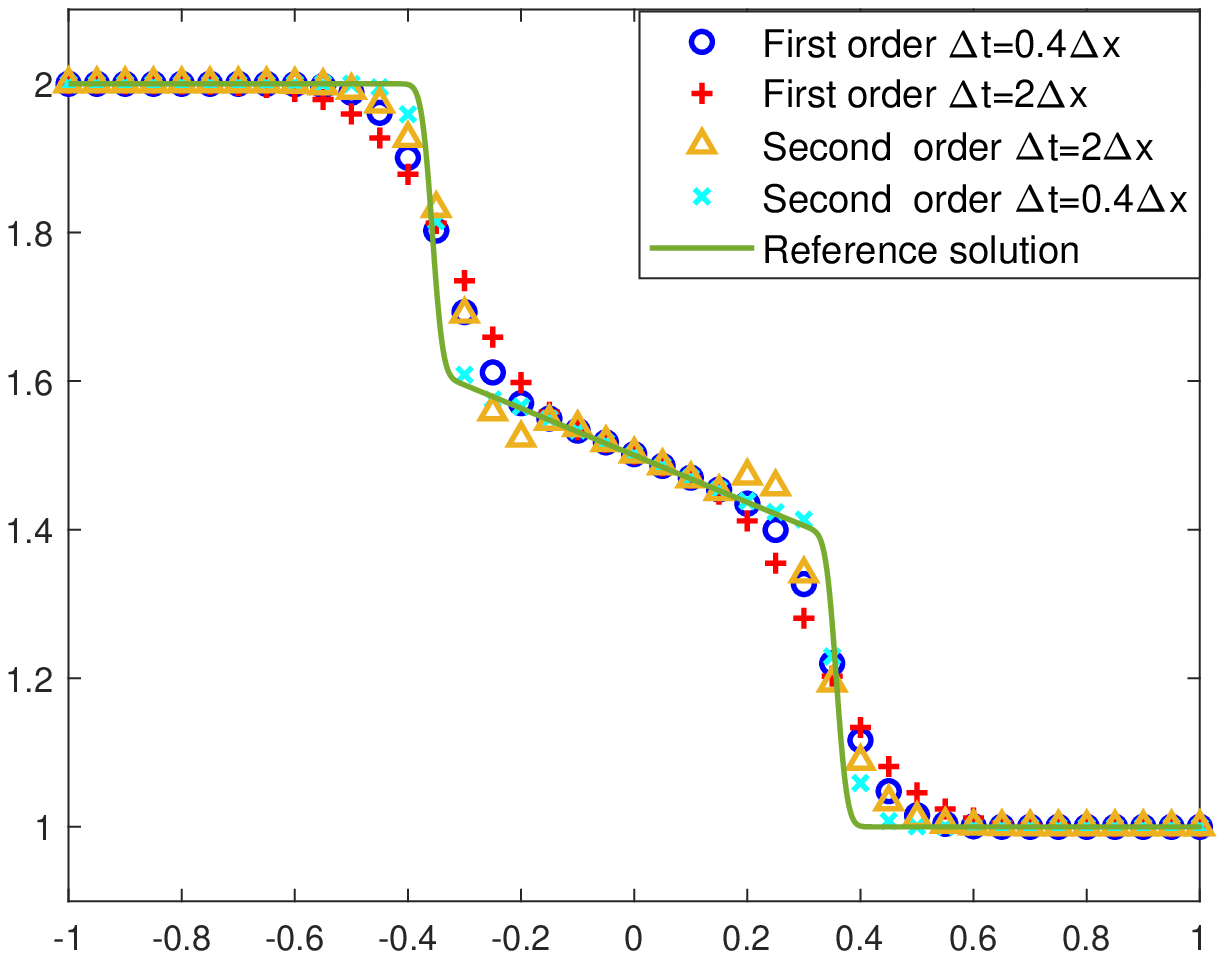}
			\end{minipage}
		}
		\subfigure[]
		{
			\begin{minipage}{7cm}
				\centering
				\includegraphics[scale=0.36]{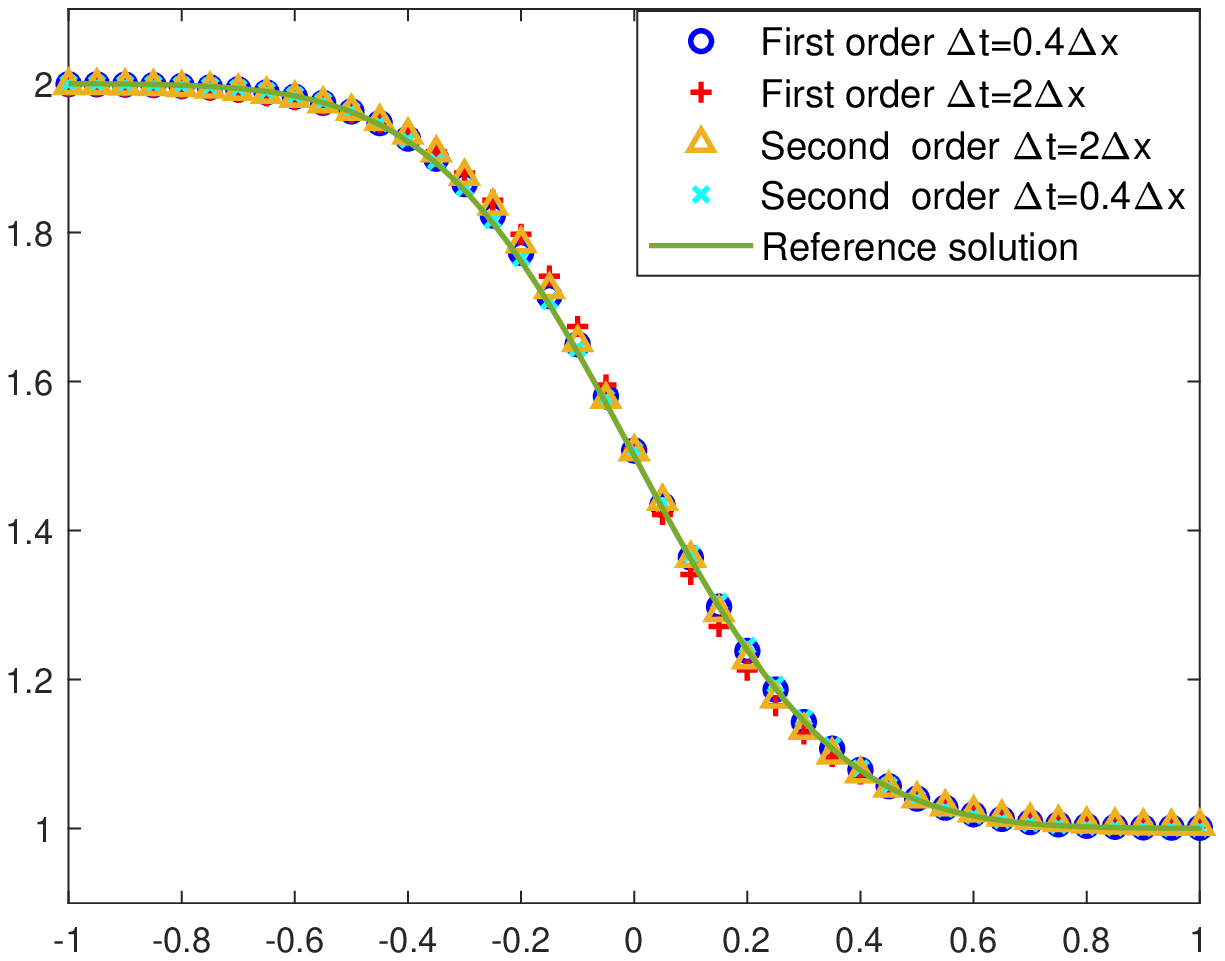}
			\end{minipage}
		}
		\caption{Numerical solution of $\rho$ for the Riemann problem \eqref{tele_dis} of the telegraph equation. (a):  a kinetic regime with $\varepsilon=0.7$ at $T=0.25$; (b): a  parabolic regime with $\varepsilon=10^{-6}$ at $T=0.04$.  \label{TG_1order} }
	\end{figure}

	Here the reference solution is computed by the first order scheme using a much refined mesh with $N=5000$, $\Delta t=0.2\Delta x$. In Fig. \ref{TG_1order}, we show the solutions from both the first and second order schemes on a uniform mesh with $N=200$. Fig. \ref{TG_1order}(a) is in the kinetic regime with $\varepsilon=0.7$ at $T=0.25$. Fig. \ref{TG_1order}(b) is in the diffusive regime with $\varepsilon=10^{-6}$ at $T=0.04$. Two time steps are used: one is $\Delta t=0.4\Delta x$ and the other is $\Delta t=2\Delta x$. We can see both schemes with both time steps agree well with the reference solution in these two different regimes.
	
\begin{figure}[!ht]
		\centering
		\subfigure[]
		{
			\begin{minipage}{7cm}
				\centering
				\includegraphics[scale=0.36]{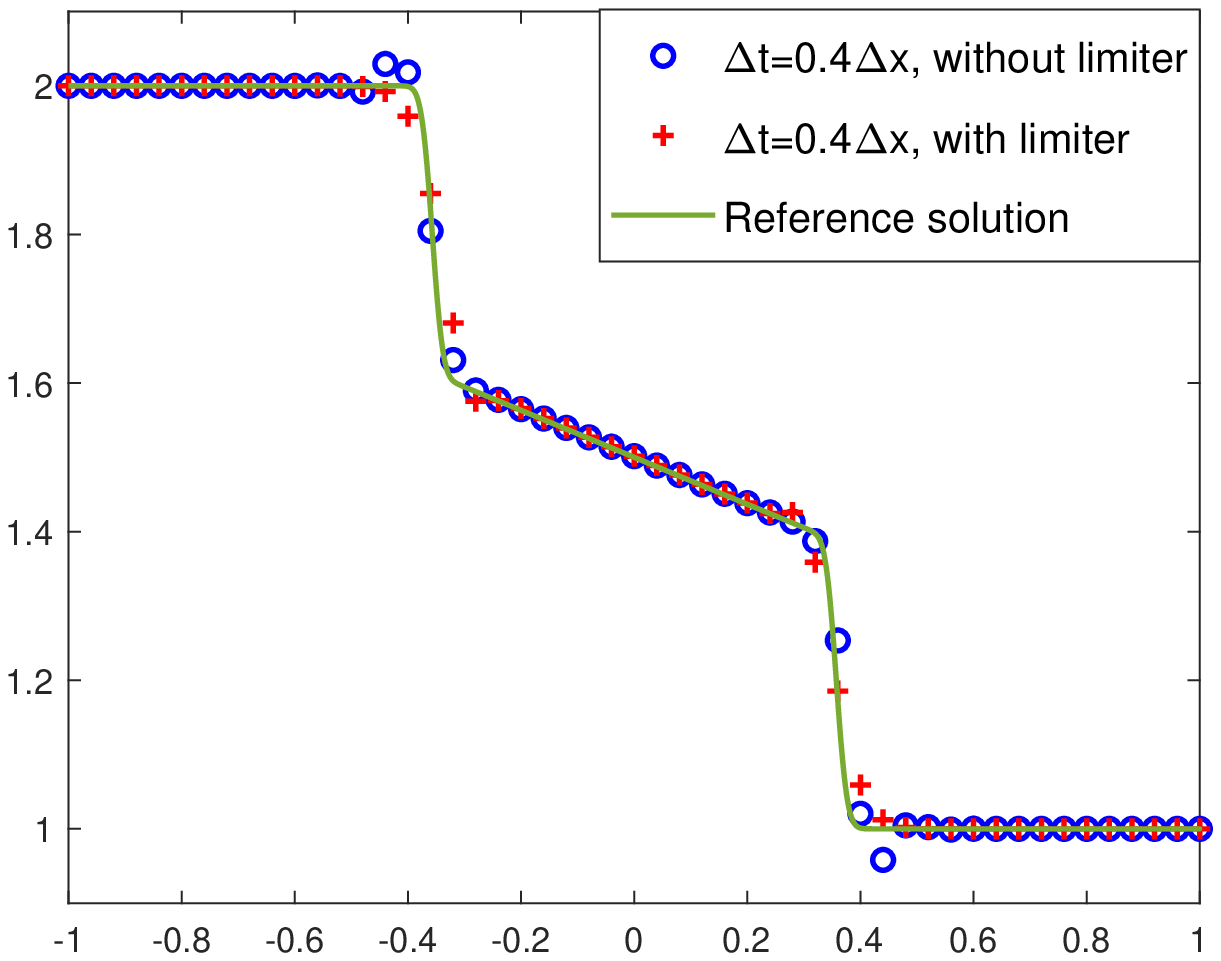}
			\end{minipage}
		}
		\subfigure[]
		{
			\begin{minipage}{7cm}
				\centering
				\includegraphics[scale=0.36]{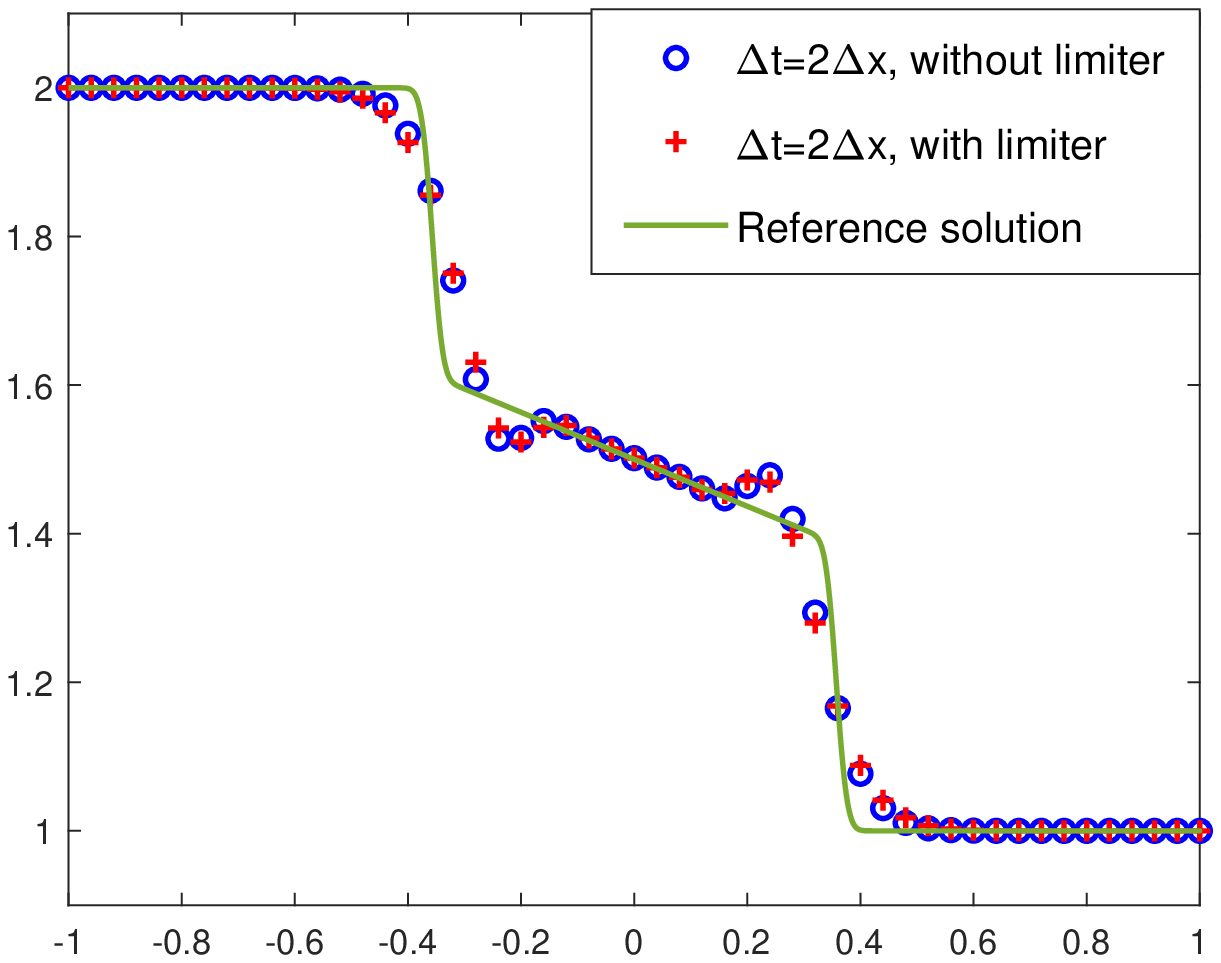}
			\end{minipage}
		}
		\caption{Second order numerical results of $\rho$ for the Riemann problem \eqref{tele_dis} of the telegraph equation in the kinetic regime with $\varepsilon=0.7$ at $T=0.25$, with and without limiters. (a): $\Delta t=0.4 \Delta x$; (b): $\Delta t=2\Delta x$.    \label{TG_ALg2} }
	\end{figure}

	For the second order scheme, due to discontinuous solutions in the kinetic regime, limiters are needed to control numerical oscillations. In Fig. \ref{TG_ALg2}, we show the results with and without limiters, using two different time steps $\Delta t=0.4 \Delta x$ and $\Delta t=2 \Delta x$. We can see that the oscillations can be well controlled by the limiters.
	
\begin{figure}[!ht]
		\centering
		\subfigure[]
		{
			\begin{minipage}{7cm}
				\centering
				\includegraphics[scale=0.36]{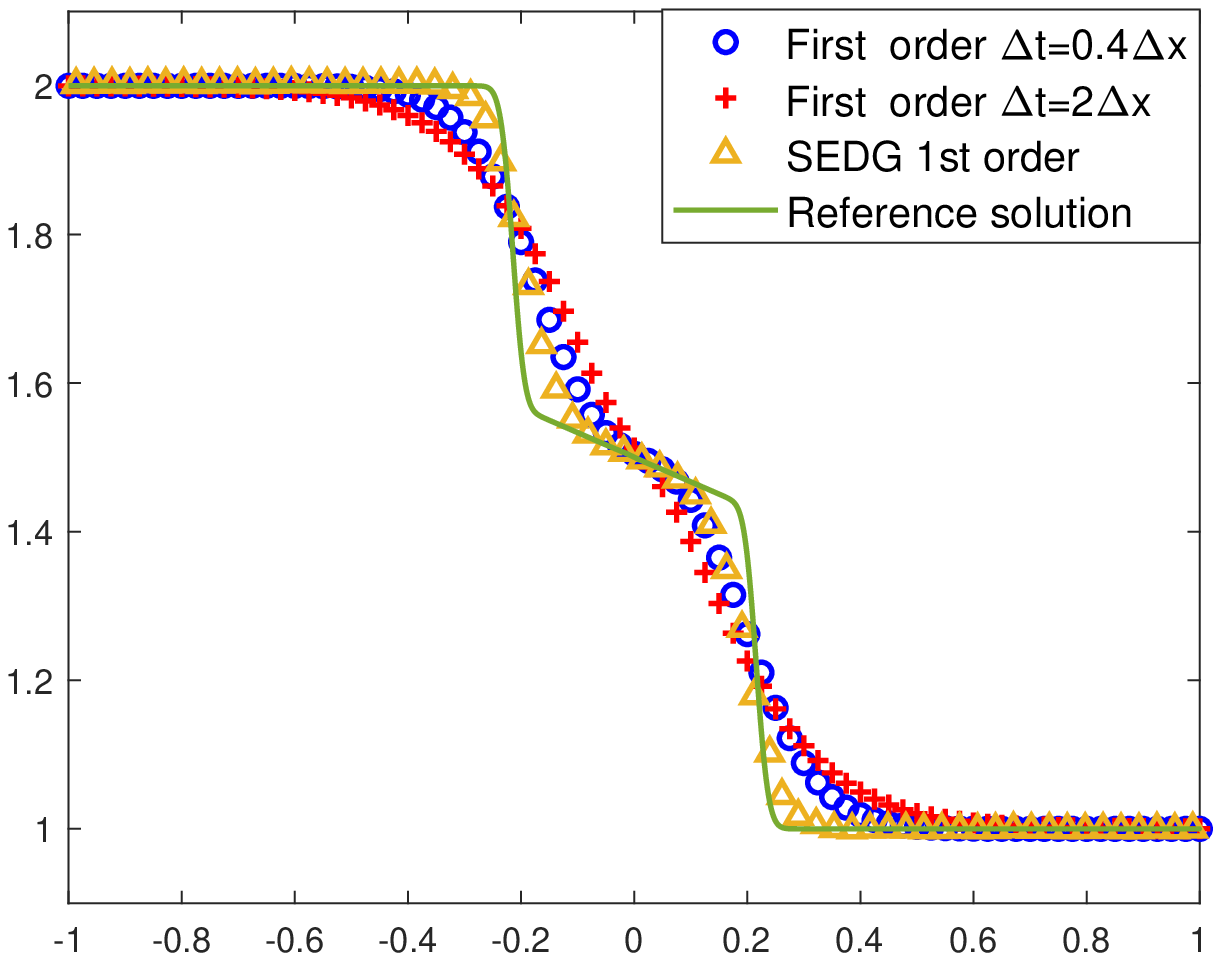}
			\end{minipage}
		}
		\subfigure[]
		{
			\begin{minipage}{7cm}
				\centering
				\includegraphics[scale=0.36]{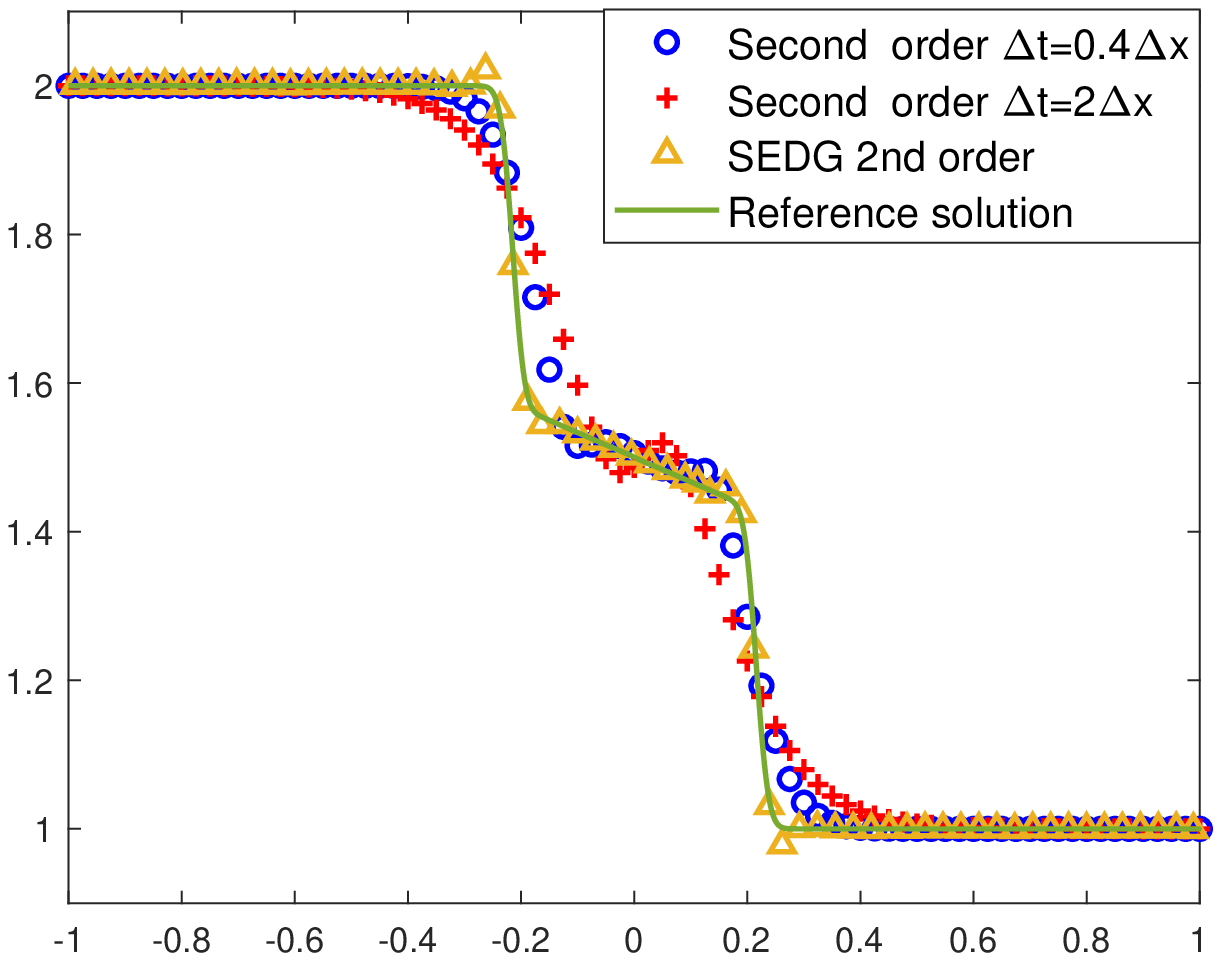}
			\end{minipage}
		}
		\caption{The comparison of our first and second order schemes with the SEDG method at the corresponding orders for the telegraph equation in the kinetic regime with $\eps=0.7$ at $T=0.15$, mesh size $\Delta x = 0.025$. (a): first order schemes; (b): second order schemes. The symbol ``$\triangle$" is the SEDG method, and the time step is $\Delta t = 0.0063 \Delta x$ in (a), $\Delta t = 0.0078\Delta x$ in (b).   \label{CompairPeng}}
	\end{figure}

	We also compare our schemes with the stability-enhanced discontinuous Galerkin (SEDG) method in \cite{peng2020stability} under a micro-macro decomposition framework, in the kinetic regime with $\varepsilon=0.7$ at time $T=0.15$. The two schemes are very similar in the diffusive regime. In \cite{peng2020stability} (Example 6.1.2), they took a mesh size $\Delta x = 0.025$ ($N=80$) $\Delta t= 0.0063\Delta x$ for the first order scheme and  $\Delta t= 0.0078\Delta x$ for the second order scheme. We take the same mesh size, but use much larger time steps $\Delta t=0.4\Delta x$ and $\Delta t=2\Delta x$. The  results are shown in Fig. \ref{CompairPeng}. As we can see, both schemes can  capture the reference solution well. Our schemes have larger deviation errors due to larger time steps, but use much less computational cost.
	
	\subsubsection{Advection-diffusion equation}
	We now consider the collision operator (\ref{Q3}), which converges to the advection-diffusion equation \eqref{adv-diff} when $\eps \rightarrow 0$.
	The approximation model, corresponding to \eqref{apprxmodel}, can be given as
	\begin{equation}
		\left\{\begin{aligned}
			&f_{t}+\frac{1}{\varepsilon} v f_{x} =\frac{1}{\varepsilon^{2}}(
			\rho-f+A \varepsilon v\rho),    \\
			&\rho_t+\frac{e^{-\mu \left(t-t_{n}\right)}}{\varepsilon} \langle v\,(f-\rho)_x\left(\bX(t_n,v),  v,  t_{n}\right)\rangle -\langle v^2\rangle  (1-e^{-\mu(t-t_n)})\rho_{xx}\left(x,t\right)\\
			&\hspace{0.3cm} +\langle v^2\rangle A  (1-e^{-\mu(t-t_n)})\rho_{x}\left(x,t\right)=0.
		\end{aligned}\right.
	\end{equation}
	If we take $A=1$, the limiting equation \eqref{adv-diff} admits the following exact solution
	\begin{equation}
		\label{241}
		\rho(x, t)=e^{-t} \sin (x-t), \quad j(x, t)=e^{-t}(\sin (x-t)-\cos (x-t)),
	\end{equation}
	on the domain $[-\pi, \pi]$ with periodic boundary conditions, where $j(x, t):=\frac{1}{2\varepsilon} (f(x,v=1,t)-f(x,v=-1,t))$.
	On the other hand for a Riemann initial data $\rho_L$ and $\rho_R$ with $j_L=j_R=0$, the limiting equation has an exact solution
	\begin{equation} \label{24}
		\rho(x, t)=\frac{1}{2}\left(\rho_{L}+\rho_{R}\right)+\frac{1}{2}\left(\rho_{L}-\rho_{R}\right) \text{erf}\left(\frac{(t-x)}{2 \sqrt{t}}\right),
	\end{equation}
	where $\text{erf}$ denotes the error function.
	
	First, we use the smooth solution \eqref{241} to test the orders of accuracy for our schemes. In order to compare with the exact solution \eqref{241} in the limit, we take a very small $\eps=10^{-6}$, so that the error between two models can be ignored. In Table \ref{Table3}, we show the $L^\infty$ and $L^1$ errors for first and second order schemes with their corresponding orders at $T=1$, with $\Delta t=3\Delta x$. We can observe the expected first and second convergence orders. We also test much larger time steps, such as $\text{CFL}=10,\,20$. The designed orders are still observed (we omit them to save space), which shows the ability of our approach with large time steps.
	
	\begin{table}[htbp]
		\scriptsize
		\centering
		\caption{ $L^{\infty}$ and $L^{1}$  errors and convergence orders in space for $\rho$   and $f(x,1,T)$. $\eps=10^{-6}$. $ \Delta t=3\Delta x$. \label{Table3}}
		\begin{tabular}{|c|c|c|c|c|c|c|c|c|c|}
			\hline	
			\multirow{6}{1.2cm}{1st order scheme}         	
			&N  & $L^{\infty}$ error of f & Order & $L^{\infty}$ error of $\rho$ & Order & $L^1$ error of f & Order & $L^1$ error of $\rho$ & Order  \\ \hline \hline
			
			& 40 & 3.47E-1 & -- &  3.47E-1& --    & 2.23E-1 & -- &  2.23E-1& --     \\ \cline{2-10}
			& 80 & 1.57E-1 & 1.15 &  1.57E-1& 1.15 &  1.01E-1& 1.15   &  1.01E-1& 1.15      \\ \cline{2-10}
			& 160 & 6.52E-2 & 1.27 & 6.52E-2& 1.27  &  4.17E-2& 1.27 &  4.17E-2& 1.27   \\ \cline{2-10}
			& 320 & 2.15E-2 & 1.60 &  2.15E-2 & 1.60   &  1.37E-2& 1.60 &  1.37E-2& 1.60   \\ \cline{2-10}
			& 640 & 1.12E-2 & 0.95   & 1.12E-2 & 0.95   &  7.14E-3& 0.95   &  7.14E-3& 0.95    \\ \hline \hline
			\multirow{6}{1.2cm}{2nd order scheme}         	
			
			& 40 & 8.91E-2 & -- &  8.91E-2& --    & 5.78E-2 & -- & 5.78E-2 & --    \\ \cline{2-10}
			& 80 & 2.56E-2 & 1.82 & 2.56E-2 & 1.82 &  1.64E-2& 1.82  &  1.64E-2& 1.82    \\ \cline{2-10}
			& 160 & 6.33E-3 & 2.02 &6.33E-3 & 2.02   &  4.03E-3& 2.02 &  4.03E-3& 2.02   \\ \cline{2-10}
			& 320 & 1.54E-3 & 2.04 & 1.54E-3 & 2.04   &  9.79E-4& 2.04 &  9.79E-4& 2.04  \\ \cline{2-10}
			& 640 & 3.65E-4 & 2.07  & 3.65E-4 & 2.07    &  2.33E-4& 2.07 &  2.33E-4& 2.07    \\ \hline
		\end{tabular}
	\end{table}
	
\begin{figure}[htbp]
		\centering
		\subfigure[]
		{
			\begin{minipage}{7cm}
				\centering
				\includegraphics[scale=0.36]{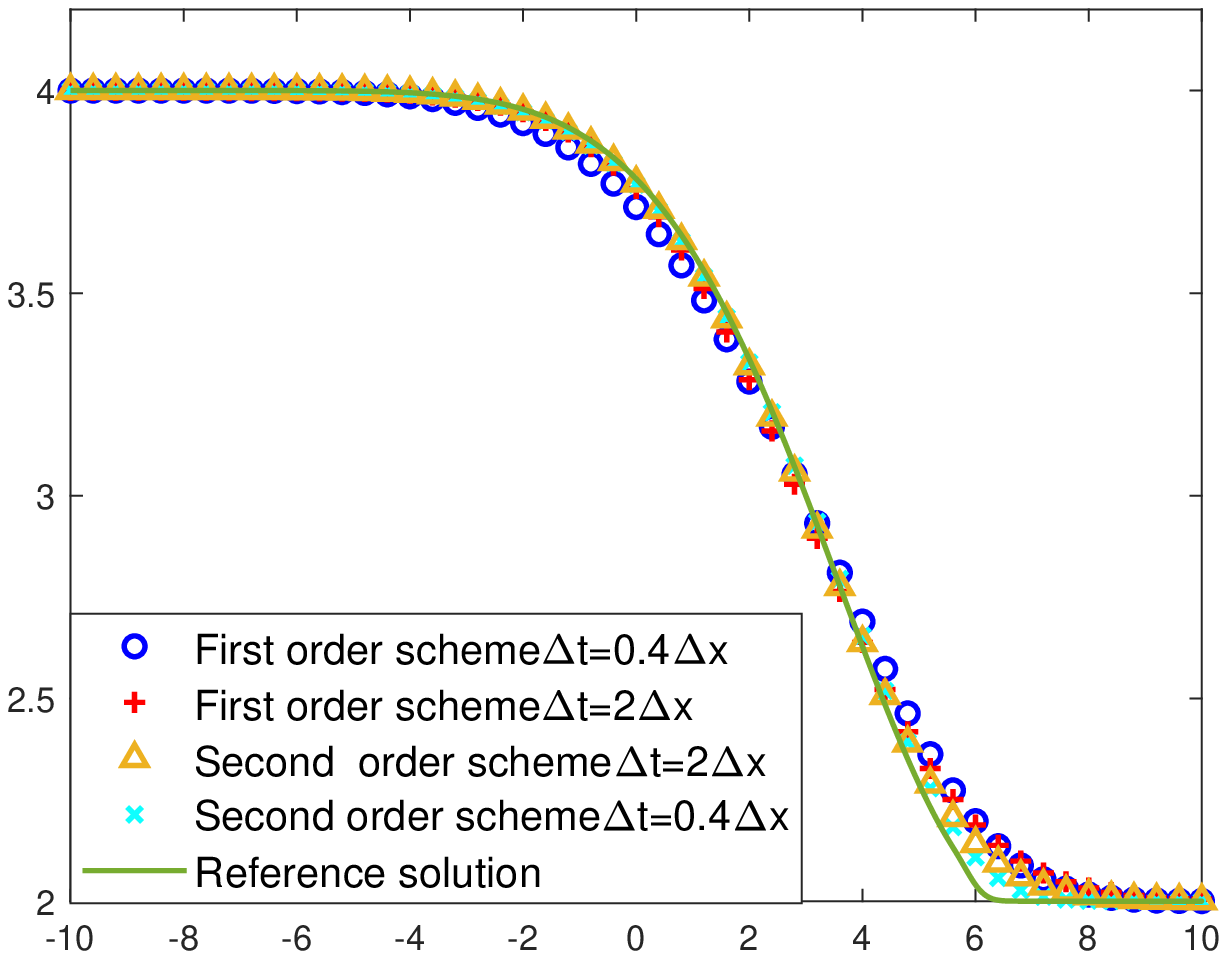}
			\end{minipage}
		}
		\subfigure[]
		{
			\begin{minipage}{7cm}
				\centering
				\includegraphics[scale=0.36]{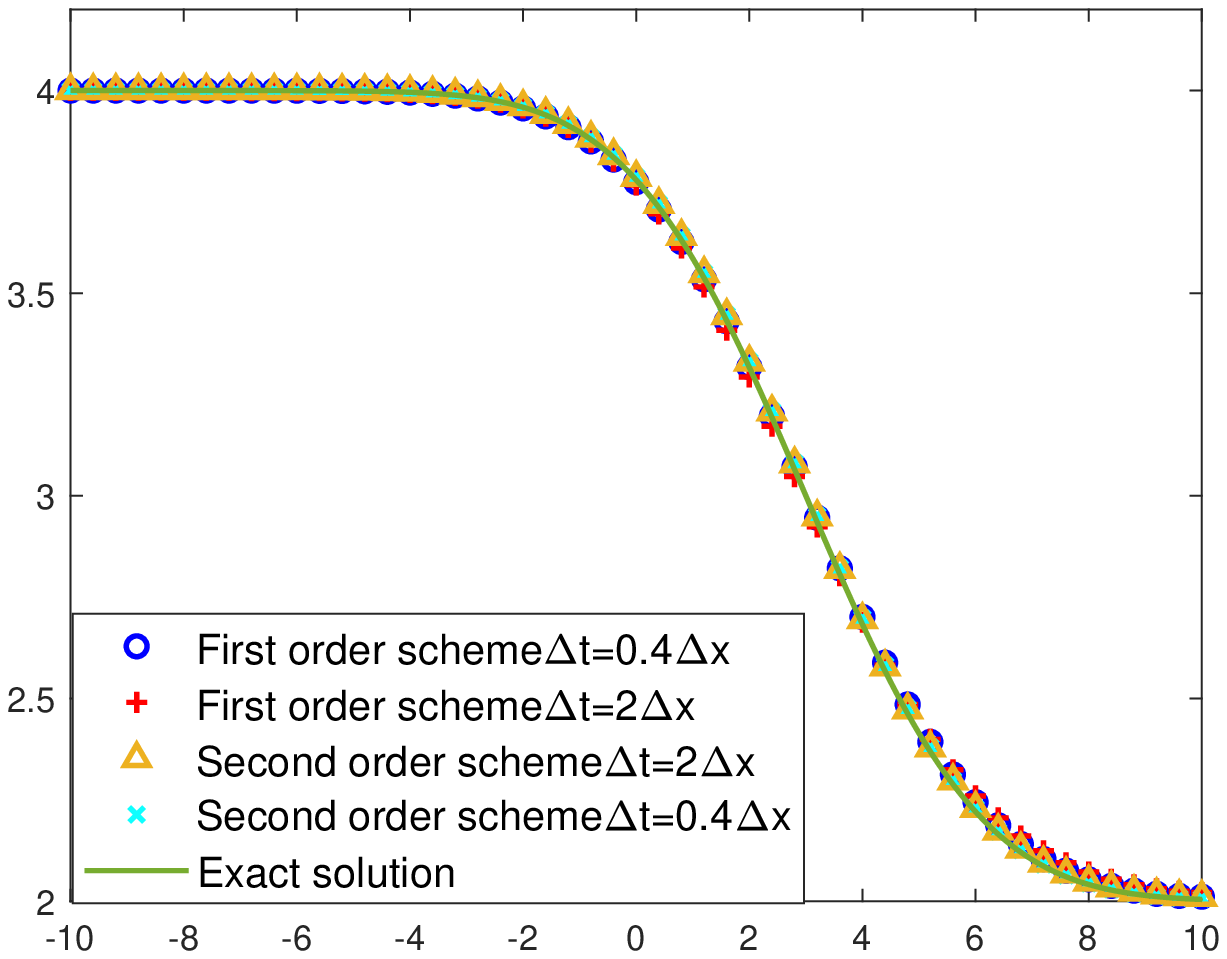}
			\end{minipage}
		}
		\caption{Numerical results of $\rho$ for the Riemann problem of the advection-diffusion equation with different $\varepsilon$'s. (a):  a rarefied regime with $\varepsilon=0.5$ at $T=3$; (b):  a diffusive regime with $\varepsilon=10^{-6}$ at $T=3$.  \label{AD_two}}
	\end{figure}

	Next, we test the Riemann problem with   the initial condition
	\begin{equation}
		\left\{\begin{array}{llll}
			\rho(x,0)=4.0, & f(x,v,0)=4.0,   & -10<x<0, \\
			\rho(x,0)=2.0, & f(x,v,0)=2.0,  & 0<x<10,
		\end{array}\right.
	\end{equation}
	with inflow and  outflow boundary conditions. We take $\varepsilon =0.5$ and $\varepsilon =10^{-6}$, which correspond to a rarefied regime and a diffusive regime, respectively. The numerical solutions computed with $N=200$ are displayed in Fig. \ref{AD_two}.  Fig. \ref{AD_two}(a) gives the numerical solutions in the rarefied regime with $\varepsilon=0.5$ at $T=3$. The solid line is a reference solution obtained by the first order scheme using a much refined mesh $N=5000$ and $\Delta t=0.2\Delta x$. We observe that both the first and second order schemes, taking either a small time step $\Delta t= 0.4 \Delta x$ or a large time step $\Delta t=2\Delta x$, can capture the reference solution well. The larger time step solutions show a little larger errors as expected. In addition,  the second order scheme with a large time step has comparable results to the first order scheme with a small time step, which shows the advantage of high order schemes even  with large time steps.
	Fig. \ref{AD_two}(b) presents the numerical solutions in the diffusive regime with $\varepsilon =10^{-6}$ at $T=3$. The solid line is the exact solution (\ref{24}). In this regime, two different time steps lead to very close solutions. Namely in this regime, large time steps would be more efficient.
	
	\subsubsection{Viscous nonlinear Burgers' equation}
	Here we consider the collision operator (\ref{Q4}), which converges to the viscous Burgers' equation \eqref{visburg} as $\eps\ra0$. The corresponding approximation model is
	\begin{equation}
		\left\{\begin{aligned} \label{burgers_mode2}
			f_{t}&+\frac{1}{\varepsilon} v f_{x} =\frac{1}{\varepsilon^{2}}(
			\rho-f+C\varepsilon v(\rho^2-(\rho-f)^2)),  \\
			\rho_t&+\frac{e^{-\mu \left(t-t_{n}\right)}}{\varepsilon} \langle v\,(f-\rho)_x\left(\bX(t_n,v),  v,  t_{n}\right)\rangle -\langle v^2\rangle (1-e^{-\mu(t-t_n)})\rho_{xx}\left(x,t\right)\\
			&+ C(1-e^{-\mu(t-t_n)})\langle v^2 \rangle (\rho^2(x,t))_x -C(1-e^{-\mu(t-t_n)}) \left\langle v^2  \left(\rho-f  \right)^2(\bX(t_n,v),v,t_n) \right\rangle_x =0.
		\end{aligned}\right.
	\end{equation}
	Due to the nonlinear advection term, we use a Picard iterative method to solve this simple nonlinear system. Taking the first order time discretization \eqref{time_dis1} as an example, the Picard iterative procedure is as follows:
	\begin{equation}\label{burgers_model}
		\left\{\begin{array}{l}
			\frac{f^{n+1,(\ell)}-f^n}{\Delta t}+\frac{1}{\varepsilon} v f^{n+1,(\ell)}_{x} =\frac{1}{\varepsilon^{2}}\Big(
			\rho^{n+1}-f^{n+1,(\ell)}+C\varepsilon v\big((\rho^{n+1})^2-(\rho^{n+1}-f^{n+1,(\ell-1)})^2\big)\Big),  \\
			\frac{\rho^{n+1,(\ell)}-\rho^n}{\Delta t}+\frac{e^{-\mu \Delta t}}{\varepsilon} \langle v\,(f-\rho)_x\left(\bX(t_n,v),  v,  t_{n}\right)\rangle -\langle v^2\rangle (1-e^{-\mu\Delta t})\rho^{n+1,(\ell)}_{xx}\left(x,t\right)\\
			\hspace{0.5cm}+ C(1-e^{-\mu\Delta t})\langle v^2 \rangle ((\rho^{n+1,(\ell-1)})^2(x,t))_x -C(1-e^{-\mu\Delta t}) \left\langle v^2  \left(\rho-f  \right)^2(\bX(t_n,v),v,t_n) \right\rangle_x =0,
		\end{array}\right.
	\end{equation}
	where $\ell$ is the iterative variable. We first solve the second equation in \eqref{burgers_model} separately to get $\rho^{n+1}$. With $\rho^{n+1}$, we then solve the first equation to update $f^{n+1}$.
	
	Taking $C=1/2$, we test the problem with a Riemann initial data
	\begin{equation}
		\left\{\begin{array}{llll}
			\rho(x,0)=2.0, & f(x,1,0)=\rho+\varepsilon j, & f(x,-1,0)=\rho-\varepsilon j,  & -10<x<0, \\
			\rho(x,0)=1.0, & f(x,1,0)=\rho+\varepsilon j, & f(x,-1,0)=\rho-\varepsilon j, & 0<x<10,
		\end{array}\right.
	\end{equation}
	where $j=\rho^{2} /\left(1+\sqrt{1+\rho^{2} \varepsilon^{2}}\right)$. Inflow and outflow boundary conditions are used.
	
\begin{figure}[htbp]
		\centering
		\subfigure[]
		{
			\begin{minipage}{7cm}
				\centering
				\includegraphics[scale=0.36]{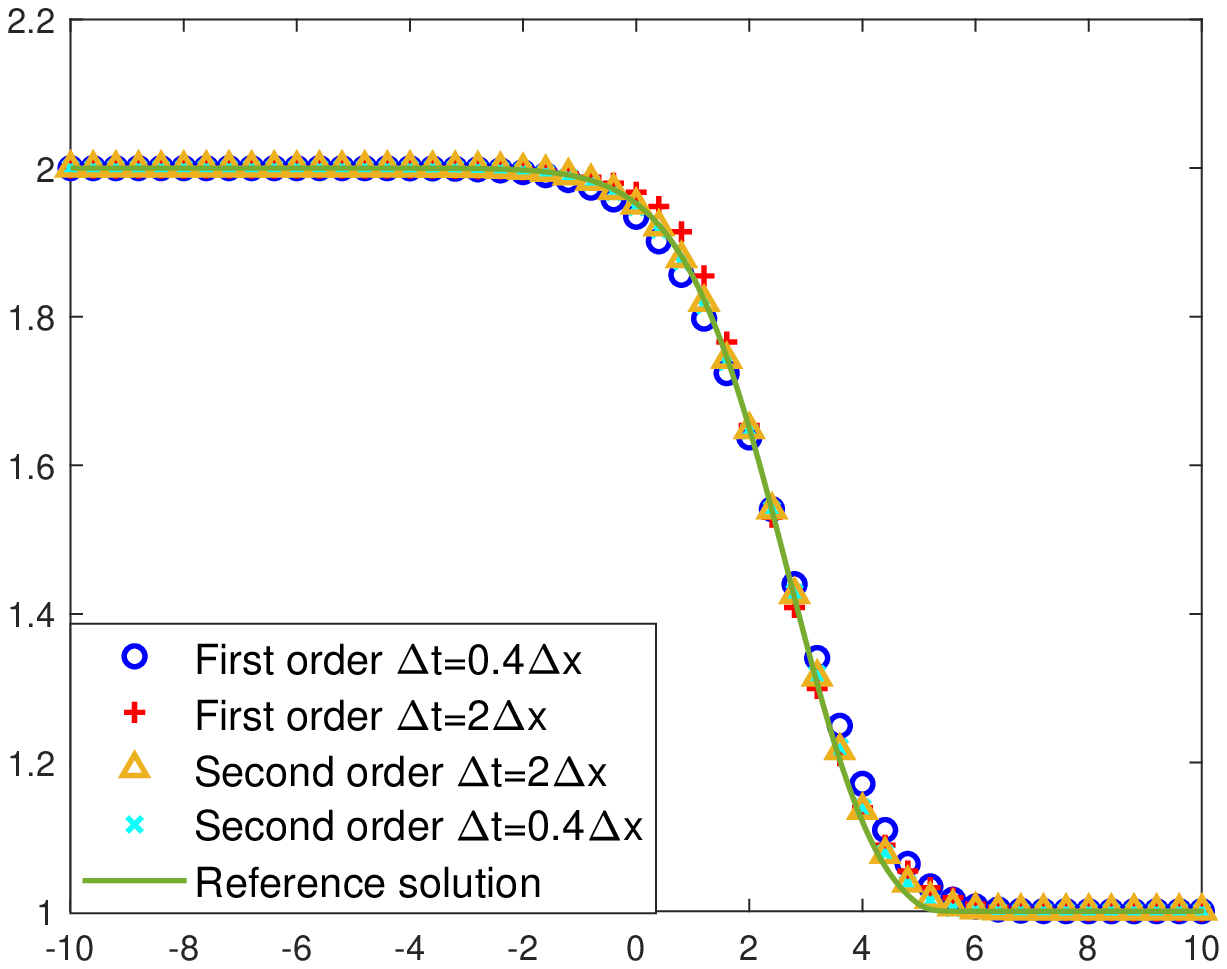}
			\end{minipage}
		}
		\subfigure[]
		{
			\begin{minipage}{7cm}
				\centering
				\includegraphics[scale=0.36]{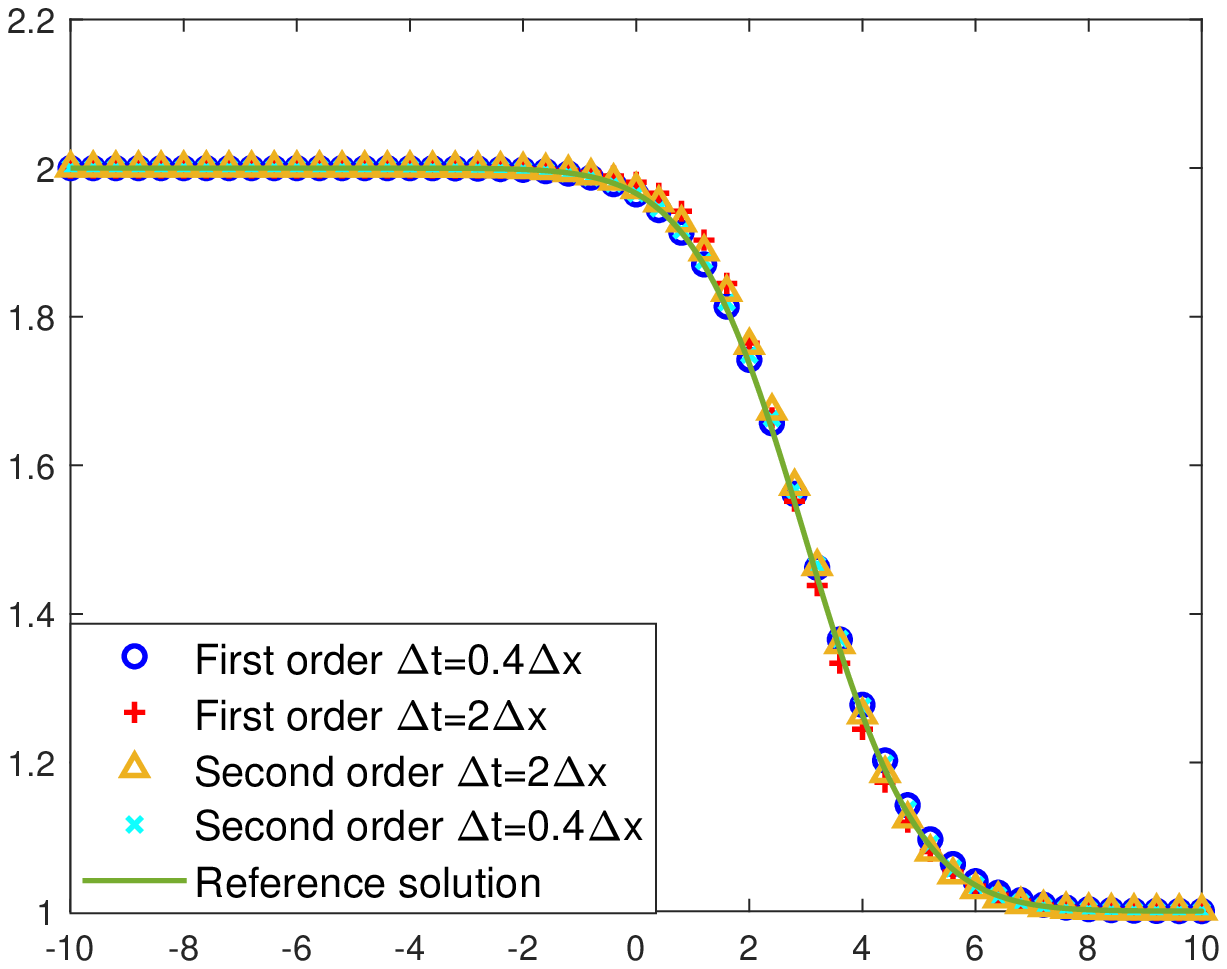}
			\end{minipage}
		}
		\caption{Numerical results of $\rho$ for the Riemann problem of the viscous Burgers' equation with different $\varepsilon$'s. (a):  a rarefied regime with $\varepsilon=0.4$ at $T=2$, (b):  a diffusive regime with $\varepsilon=10^{-6}$ at $T=2$. \label{Burgers_kinetic} }
	\end{figure}

	The numerical solutions at $T=2$  are displayed in Fig. \ref{Burgers_kinetic}, using $N=200$.  Fig. \ref{Burgers_kinetic}(a) is for a rarefied regime with $\varepsilon=0.4$.
	Fig. \ref{Burgers_kinetic}(b) is for a diffusive regime with $\varepsilon =10^{-6}$. The solid lines are reference solutions obtained by the first order scheme using a much refined mesh $N=5000$ and $\Delta t=0.2\Delta x$. We can see that our schemes, with a small time step $\Delta t= 0.4 \Delta x$ or a large time step $\Delta t=2\Delta x$, agree well with the reference solution. All results do not show obvious differences, so large time steps would be more efficient. {For this problem with a nonlinear diffusive limit, using a Picard iteration, by taking a termination criteria $ \text{max} (|\rho^{n+1,l}-\rho^{n+1,l-1}|)<10^{-8}$, the number of iterations are around $15$ to $19$.}
	
	\subsection{ One-group transport equation in a slab geometry}
	\label{sec_numer_oneg}
	In this section, we focus on a one-group kinetic transport equation
	\begin{align}
		\label{one-group}
		\varepsilon \partial_{t} f+v \partial_{x} f=\frac{\sigma_S}{\varepsilon}(\langle f\rangle-f)-\varepsilon \sigma_{A} f+\varepsilon G,
	\end{align}
	where $v$ is a continuous velocity within $\Omega_v = [-1,1]$. $\sigma_S$ and $\sigma_{A}$ are the scattering and absorbing coefficients, respectively. $G=G(x)$ is a source term. If we consider $\sigma_S$, $\sigma_A$ as constants and define $\mu= \frac{\sigma_S}{\varepsilon^2}+\sigma_A$,
	following the approach in Section \ref{sec_model_apprx}, we can derive an approximation model for \eqref{one-group}:
	
	\begin{equation}
		\label{apprxmodeburgers}
		\left\{
		\begin{array}{l}
			f_{t}+\frac{1}{\varepsilon} v f_{x} =\frac{\sigma_S}{\varepsilon^2}(\rho-f)- \sigma_{A} f+G(x), \\ \, \\
			\rho_t\left(x,  t\right)+\frac{e^{-\mu \left(t-t_{n}\right)}}{\varepsilon} \langle v(f-\rho)_x\left(\bX(t_n,v),  v,  t_{n}\right)\rangle
			-\langle v^2\rangle \frac{1}{\sigma_S+\eps^2\,\sigma_A} (1-e^{-\mu(t-t_n)})\rho_{xx}\left(x,  t\right) =0.\\
		\end{array}
		\right.
	\end{equation}
	For the one-group velocity model, a DOM with $16$ Gaussian quadrature points is used for $v\in[-1,1]$. In this case, our first and second order schemes developed in Section \ref{sec_scheme} can be directly applied.
	
	\subsubsection{Accuracy and large time step tests for a smooth solution}
	We consider the following initial condition
	\beq
	\label{eq521}
	\left\{\begin{array}{l}
		\rho(x, 0)=2+\sin (x), \\
		f(x, v, 0)=2+\sin (x)-\varepsilon v \cos (x)
	\end{array}\right.
	\eeq
	on the domain $\Omega_x=[-\pi, \pi]$ with a periodic boundary condition. We compute the solution up to a final time $T=1$. We take $\sigma_S= 1$ and $\sigma_{A} = G = 0$. Several different $\eps$'s are tested and we show the results for $\varepsilon=0.5, 0.1, 10^{-2}, 10^{-6}$. Since the exact solutions are not available, we compute a reference solution using the second order scheme with $N=5120$ and $\Delta t=0.0005$.
	
	In Table \ref{table_one_space1} and Table \ref{table_one_space2}, we show the $L^1$ and $L^{\infty}$ errors and orders for $\rho$, and  $f(x,v_1,T)$ which corresponds to the value of $f$ at the first Gaussian point in velocity. 
	We take a number of grid points $N=40\cdot 2^k \quad (k=0,1,...,5)$, and a corresponding time step $\Delta t=T/(3\cdot2^k)$. The CFL number is approximately $2.12$. The expected first and second orders of accuracy are observed in the kinetic and diffusive regimes.
	
	To test the accuracy in time, we take a fixed mesh with $N=5120$ but adopt different time steps. The $L^1$ and $L^{\infty}$ errors and orders are tabulated in Table \ref{table_one_time1} and Table \ref{table_one_time2}.
	Compared with the spatial mesh size, the time steps here are extremely large. We can clearly observe first and second orders of accuracy in both the kinetic and diffusive regimes. Here for the second order scheme in Table \ref{table_one_space2} and Table \ref{table_one_time2}, in the intermediate regime ($\eps=0.1$), an order reduction occurs when $\eps$ and $\Delta t$ are comparable, which is similar to that in Table \ref{table2} and Table \ref{table_TG_time2}.
	
	\begin{table}[!ht]
		\scriptsize
		\centering
		\caption{$L^{\infty}$ and $L^{1}$  errors and convergence orders in space for $\rho$   and $f(x,v_1,T)$.  First order scheme for one-group transport equation.}
		\label{table_one_space1}
		\begin{tabular}{|c|c|c|c|c|c|c|c|c|c|} \hline	
			$\eps$ &N  & $L^{\infty}$ error of $\rho$ & Order & $L^{\infty}$ error of $f$ & Order  & $L^{1}$ error of $\rho$ & Order & $L^{1}$ error of $f$ & Order        \\ \hline
			\multirow{6}{*}{ 0.5}         	
			
			& 40 & 2.96E-2 & -- &  1.93E-2& --  &1.89E-2 & -- &  1.23E-2& --                        \\ \cline{2-10}
			& 80 & 1.73E-2 & 0.78 &  1.06E-2& 0.87  & 1.10E-2 & 0.77 & 6.74E-3& 0.86             \\ \cline{2-10}
			& 160 & 8.70E-3 & 0.99 &  4.89E-3& 1.11 & 5.54E-3 & 1.00 & 3.12E-3& 1.11     \\ \cline{2-10}
			& 320 & 4.07E-3 & 1.09 &  2.15E-3& 1.18   &2.59E-3 & 1.09 &  1.37E-3& 1.19     \\ \cline{2-10}
			& 640 & 1.90E-3 & 1.10   & 9.28E-4& 1.21   & 1.21E-3 & 1.10 & 5.91E-4& 1.21       \\  \hline
			\multirow{5}{*}{ 0.1}
			& 40 & 1.69E-2 & -- &  1.14E-2& --  &1.08E-2 & -- &  7.23E-3& --                        \\ \cline{2-10}
			& 80 & 8.86E-3 & 0.93 & 6.07E-3& 0.90 & 5.63E-3 & 0.93& 3.87E-3& 0.90             \\ \cline{2-10}
			& 160 & 4.64E-3 & 0.93 & 3.25E-3& 0.90 & 2.95E-3 & 0.93 & 2.07E-3& 0.90     \\ \cline{2-10}
			& 320 & 2.52E-3 & 0.88 & 1.82E-3& 0.83   &1.60E-3 & 0.88 & 1.16E-3& 0.83     \\ \cline{2-10}
			& 640 & 1.34E-3 & 0.92   & 9.88E-4& 0.85   & 8.50E-4 & 0.92 & 6.29E-4& 0.89       \\  \hline
			\multirow{5}{*}{ $10^{-2}$}
			& 40 &1.65E-2 & --    &1.59E-2 & --     & 1.05E-2 & -- &  1.01E-2 & --                        \\ \cline{2-10}
			& 80 & 8.41E-3 & 0.97 &   8.13E-3 & 0.97   & 5.35E-3 & 0.97 & 5.18E-3 & 0.97          \\ \cline{2-10}
			& 160 &4.19E-3 & 1.01 & 4.06E-3 & 1.00   & 2.67E-3   &1.00    & 2.58E-3   &1.00    \\ \cline{2-10}
			& 320 &2.03E-3 & 1.04 & 1.97E-3 & 1.04   & 1.29E-3 & 1.04 & 1.25E-3 & 1.04   \\ \cline{2-10}
			& 640 & 8.77E-4 & 1.22  & 8.47E-4 & 1.22   & 5.56E-4 & 1.22 & 5.40E-4 & 1.22      \\  \hline
			\multirow{5}{*}{$ 10^{-6}$}
			& 40 &1.65E-2 & --    &1.65E-2 & --     & 1.05E-2 & -- &  1.05E-2 & --                        \\ \cline{2-10}
			& 80 & 8.41E-3 & 0.97 &   8.41E-3 & 0.97   & 5.35E-3 & 0.97 & 5.35E-3 & 0.97          \\ \cline{2-10}
			& 160 &4.19E-3 & 1.01 & 4.19E-3 & 1.00   & 2.67E-3   &1.00    & 2.67E-3   &1.00    \\ \cline{2-10}
			& 320 &2.03E-3 & 1.04 & 2.03E-3 & 1.04   & 1.29E-3 & 1.04 & 1.29E-3 & 1.04   \\ \cline{2-10}
			& 640 & 8.74E-4 & 1.22  & 8.74E-4 & 1.22   & 5.56E-4 & 1.22 & 5.56E-4 & 1.22      \\  \hline
			
		\end{tabular}
	\end{table}

	\begin{table}[!ht]
		\scriptsize
		\centering
		\caption{$L^{\infty}$ and $L^{1}$  errors and convergence orders in space for   $\rho$ and $f(x,v_1,T)$.  Second order scheme for one-group transport equation. \label{table_one_space2}}
		\begin{tabular}{|c|c|c|c|c|c|c|c|c|c|}
			\hline	
			$\eps$ &N  & $L^{\infty}$ error of $\rho$ & Order & $L^{\infty}$ error of $f$ & Order  & $L^{1}$ error of $\rho$ & Order & $L^{1}$ error of $f$ & Order        \\ \hline
			\multirow{6}{*}{ 0.5}         	
			
			& 40 & 1.63E-2 & -- &  1.20E-2& --  &1.04E-2 & -- & 7.89E-3& --                        \\ \cline{2-10}
			& 80 & 6.29E-3 & 1.37 & 3.07E-3& 1.96  & 4.01E-3 & 1.37 & 1.98E-3& 1.99              \\ \cline{2-10}
			& 160 &2.02E-3 & 1.64 &  1.04E-3& 1.56 & 1.29E-3 & 1.64 &  6.66E-4& 1.58     \\ \cline{2-10}
			& 320 & 5.69E-4 & 1.83 &  3.25E-4& 1.68   & 3.62E-4 & 1.83 & 2.07E-4& 1.68   \\ \cline{2-10}
			& 640 & 1.38E-4 & 2.04   & 8.20E-5& 1.99   & 8.78E-5 & 2.05 &5.22E-5& 1.99       \\  \hline
			\multirow{8}{*}{ 0.1}         	
			& 40 & 6.24E-3 & - & 6.12E-3 & - & 3.97E-3 & - & 3.93E-3 & - \\           \cline{2-10}
			& 80 & 1.54E-3 & 2.02 & 1.52E-3 & 2.01 & 9.81E-4 & 2.02 & 9.73E-4 & 2.01     \\ \cline{2-10}
			& 160 & 3.50E-4 & 2.14 & 3.47E-4 & 2.13 & 2.23E-4 & 2.14 & 2.21E-4 & 2.14    \\ \cline{2-10}
			& 320 & 2.20E-4 & 0.67 & 2.18E-4 & 0.67 & 1.40E-4 & 0.67 & 1.39E-4 & 0.67     \\ \cline{2-10}
			& 640 & 2.92E-4 & -0.40 & 2.88E-4 & -0.40 & 1.86E-4 & -0.40 & 1.83E-4 & -0.40    \\ \cline{2-10}
			& 1280 &1.08E-4 & 1.43 & 1.07E-4 & 1.43 & 6.89E-5 & 1.43 & 6.81E-5 & 1.43    \\   \hline
			\multirow{5}{*}{$ 10^{-2}$}
			& 40   &  6.29E-3 & -- &   6.29E-3& --  &4.00E-3 & -- &  4.00E-3 & --                      \\ \cline{2-10}
			& 80     & 1.59E-3 & 1.98  & 1.59E-3 & 1.98  & 1.01E-3 & 1.98& 1.01E-3 & 1.98         \\ \cline{2-10}
			& 160   & 3.88E-4   &2.03   & 3.88E-4   &2.03   & 2.47E-4 & 2.03& 2.47E-4 & 2.03   \\ \cline{2-10}
			& 320    & 9.39E-5 & 2.05 & 9.39E-5 & 2.05  & 5.98E-5 & 2.05 & 5.98E-5 & 2.05 \\ \cline{2-10}
			& 640    & 2.12E-5 & 2.15  & 2.12E-5 & 2.15  & 1.35E-5 & 2.15  & 1.35E-5 & 2.15      \\  \hline
			\multirow{5}{*}{$ 10^{-6}$}
			& 40   &  6.29E-3 & -- &   6.29E-3& --  &4.00E-3 & -- &  4.00E-3 & --                      \\ \cline{2-10}
			& 80     & 1.59E-3 & 1.98  & 1.59E-3 & 1.98  & 1.01E-3 & 1.98& 1.01E-3 & 1.98         \\ \cline{2-10}
			& 160   & 3.88E-4   &2.03   & 3.88E-4   &2.03   & 2.47E-4 & 2.03& 2.47E-4 & 2.03   \\ \cline{2-10}
			& 320    & 9.39E-5 & 2.05 & 9.39E-5 & 2.05  & 5.98E-5 & 2.05 & 5.98E-5 & 2.05 \\ \cline{2-10}
			& 640    & 2.12E-5 & 2.15  & 2.12E-5 & 2.15  & 1.35E-5 & 2.15  & 1.35E-5 & 2.15      \\  \hline
			
		\end{tabular}
	\end{table}

	\begin{table}[htbp]
		\scriptsize
		\centering
		\caption{ $L^{\infty}$ and $L^{1}$  errors and convergence orders in time for $\rho$  and $f(x,v_1,T)$. First order scheme   for one-group equation with a fixed mesh $N=5120$.  \label{table_one_time1} }
		\begin{tabular}{|c|c|c|c|c|c|c|c|c|c|}
			\hline	
			$\eps$ &$\Delta t$  & $L^{\infty}$ error of $\rho$ & Order & $L^{\infty}$ error of $f$ & Order  & $L^{1}$ error of $\rho$ & Order & $L^{1}$ error of $f$ & Order        \\ \hline
			\multirow{6}{*}{$ 0.5$}         	
			& T/2 & 3.10E-2 & -- &  4.20E-2& --    & 1.97E-2 & -- &  2.67E-2& --     \\ \cline{2-10}
			& T/4 & 1.79E-2 & 0.79 & 2.22E-2& 0.92  &  1.14E-2& 0.79 &  1.41E-2& 0.92  \\ \cline{2-10}
			& T/8 & 8.26E-3 & 1.11 &  1.02E-2 & 1.13  &  5.26E-3& 1.11 & 6.48E-3& 1.13   \\ \cline{2-10}
			& T/16 & 3.68E-3 & 1.17 &  4.49E-3 & 1.18  &  2.34E-3& 1.17 & 2.86E-3& 1.18   \\ \cline{2-10}
			& T/32 & 1.40E-3 & 1.39 &  1.71E-3 & 1.39   &  8.93E-4& 1.39 &  1.09E-3& 1.39    \\ \hline
			\multirow{6}{*}{$ 0.1$}
			& T/2 & 1.82E-2 & -- &  1.82E-2& --    & 1.16E-2 & -- &  1.16E-2& --     \\ \cline{2-10}
			& T/4 & 9.47E-3 & 0.94 & 9.42E-3& 0.95  &  6.03E-3& 0.94 &  5.99E-3& 0.95  \\ \cline{2-10}
			& T/8 & 4.68E-3 & 1.02 &  4.63E-3 & 1.02  & 2.98E-3& 1.02 & 2.95E-3& 1.02   \\ \cline{2-10}
			& T/16 & 2.03E-3 & 1.20 &  1.99E-3 & 1.22  & 1.29E-3& 1.20 & 1.26E-3& 1.22   \\ \cline{2-10}
			& T/32 & 8.65E-4 & 1.23 &  8.47E-4 & 1.23  &  5.46E-4& 1.24 &  5.35E-4& 1.24   \\ \hline
			\multirow{6}{*}{$ 10^{-2}$}
			& T/2 & 1.79E-2 & -- &  1.79E-2& --    & 1.14E-2 & -- &  1.14E-2& --     \\ \cline{2-10}
			& T/4 & 9.18E-3 & 0.96 & 9.19E-3& 0.96  &  5.85E-3& 0.96 &  5.85E-3& 0.96  \\ \cline{2-10}
			& T/8 & 4.53E-3 & 1.02 &  4.53E-3 & 1.02  & 2.88E-3& 1.02 & 2.89E-3& 1.02   \\ \cline{2-10}
			& T/16 & 2.10E-3 & 1.11 &  2.10E-3 & 1.11  & 1.33E-3& 1.11 & 1.34E-3& 1.11   \\ \cline{2-10}
			& T/32 & 8.03E-4 & 1.38 &  8.08E-4 & 1.38  &  5.11E-4& 1.38 &  5.14E-4& 1.38   \\ \hline
			\multirow{6}{*}{$ 10^{-6}$}
			& T/2 & 1.82E-2 & -- &  1.82E-2& --    & 1.16E-2 & -- &  1.16E-2& --     \\ \cline{2-10}
			& T/4 & 9.49E-3 & 0.94 & 9.49E-3 & 0.94  & 6.04E-3& 0.94  &6.04E-3& 0.94 \\ \cline{2-10}
			& T/8 & 4.86E-3 & 0.97 &  4.86E-3 & 0.97  & 3.09E-3& 0.97 & 3.09E-3& 0.97    \\ \cline{2-10}
			& T/16 & 2.46E-3 & 0.98 &  2.46E-3 & 0.98  & 1.56E-3& 0.98 & 1.56E-3& 0.98    \\ \cline{2-10}
			& T/32 & 1.24E-3 & 0.99 & 1.24E-3 & 0.99  &  7.87E-4& 0.99 &7.87E-4& 0.99   \\ \hline
		\end{tabular}
	\end{table}
	
	\begin{table}[htbp]
		\scriptsize
		\centering
		\caption{$L^{\infty}$ and $L^{1}$  errors and convergence orders in time for $\rho$ and $f(x,v_1,T)$. Second order scheme   for one-group equation with a fixed mesh $N=5120$.   \label{table_one_time2} }
		\begin{tabular}{|c|c|c|c|c|c|c|c|c|c|}
			\hline	
			$\eps$ &$\Delta t$    & $L^{\infty}$ error of $\rho$ & Order & $L^{\infty}$ error of $f$ & Order  & $L^{1}$ error of $\rho$ & Order & $L^{1}$ error of $f$ & Order        \\ \hline
			\multirow{6}{*}{$ 0.5$}
			& T/2 & 2.46E-2 & -- &  2.60E-2& --    & 1.57E-2 & -- &  1.66E-2& --     \\ \cline{2-10}
			& T/4 & 1.13E-3 & 1.13 & 7.02E-3& 1.89  &  7.17E-3& 1.13 &  4.47E-3& 1.89   \\ \cline{2-10}
			& T/8 & 4.08E-3 & 1.47 &  2.13E-3 & 1.72  & 2.60E-3& 1.47 & 1.36E-3& 1.72   \\ \cline{2-10}
			& T/16 & 1.27E-3 & 1.69 &  7.47E-4 & 1.51   &  8.05E-4& 1.69 &4.75E-4& 1.51    \\ \cline{2-10}
			& T/32 & 3.56E-4 & 1.83 &  2.37E-4 & 1.74   &  2.26E-4& 1.83 &  1.42E-4& 1.74    \\ \hline
			\multirow{8}{*}{$ 0.1$}                  	
			& T/2 &  1.22E-2 & -- & 1.21E-2& --    &  7.76E-3 & -- &  7.72E-3& --   \\ \cline{2-10}
			& T/4 &  3.49E-3 & 1.81 &  3.48E-3& 1.80    & 2.22E-3 & 1.81 &  2.22E-3& 1.80   \\ \cline{2-10}
			& T/8 & 9.77E-4 & 1.84   &  9.74E-4& 1.84     & 6.22E-4 &1.84    &  6.20E-4& 1.84      \\ \cline{2-10}
			& T/16 & 4.09E-4 & 1.26 & 4.06E-4& 1.26     &  2.60E-4& 1.26 &  2.59E-4& 1.26  \\ \cline{2-10}
			& T/32 &  4.82E-4 & -0.23 &  4.77E-4& -0.23    &  3.07E-4 & -0.23 &  3.04E-4& -0.23   \\ \cline{2-10}
			& T/64 & 4.47E-4 & 0.11   &  4.42E-4& 0.11     & 2.84E-4 & 0.11    & 2.81E-4& 0.11      \\ \cline{2-10}
			& T/128 & 2.26E-4 & 0.98 & 2.24E-4& 0.98  &  1.44E-4& 0.98   &  1.42E-4& 0.99   \\ \cline{2-10}
			& T/256 & 7.73E-5 & 1.55 &  7.64E-5 & 1.55  &  4.92E-5& 1.55 &  4.86E-5& 1.55   \\ \hline
			\multirow{6}{*}{$ 10^{-2}$}
			& T/2 & 1.20E-2 & --  &1.20E-2 & --   & 7.66E-3 & -- &  7.66E-3& --     \\ \cline{2-10}
			& T/4 & 3.33E-3 & 1.85 &3.33E-3 & 1.85  &  2.12E-3& 1.85 & 2.12E-3& 1.85   \\ \cline{2-10}
			& T/8 & 8.17E-4 & 2.03 & 8.17E-4 & 2.03  & 5.20E-4& 2.03  &5.20E-4& 2.03   \\ \cline{2-10}
			& T/16 & 1.99E-4 & 2.04 &  1.99E-4 & 2.04    & 1.27E-4& 2.04 &1.27E-4& 2.04   \\ \cline{2-10}
			& T/32 & 4.73E-5 & 2.07 &4.73E-5 & 2.07    &  3.01E-5& 2.07  & 3.01E-5& 2.07    \\ \hline
			\multirow{6}{*}{$ 10^{-6}$}
			& T/2 & 1.20E-2 & --  &1.20E-2 & --   & 7.66E-3 & -- &  7.66E-3& --     \\ \cline{2-10}
			& T/4 & 3.33E-3 & 1.85 &3.33E-3 & 1.85  &  2.12E-3& 1.85 & 2.12E-3& 1.85   \\ \cline{2-10}
			& T/8 & 8.17E-4 & 2.03 & 8.17E-4 & 2.03  & 5.20E-4& 2.03  &5.20E-4& 2.03   \\ \cline{2-10}
			& T/16 & 1.99E-4 & 2.04 &  1.99E-4 & 2.04    & 1.27E-4& 2.04 &1.27E-4& 2.04   \\ \cline{2-10}
			& T/32 & 4.73E-5 & 2.07 &4.73E-5 & 2.07    &  3.01E-5& 2.07  & 3.01E-5& 2.07    \\ \hline
			
		\end{tabular}
	\end{table}	
	
	\subsubsection{ Diffusive and kinetic regimes with an isotropic boundary condition}
	We then consider \eqref{one-group} with an isotropic boundary condition. The boundary and initial values are given as
	\beq
	\label{eq522}
	\begin{gathered}
		f_{L}(v, t)=1, \quad f_{R}(v, t)=0 ; \quad f(x, v, 0)=0, \quad x \in \Omega_{x}=[0, 1].
	\end{gathered}
	\eeq
	The parameters $\sigma_S=1, \, \sigma_{A}=0, \, G=0$. We compute the numerical solutions in the kinetic regime with $\varepsilon=1$ and in the diffusive regime with $\varepsilon=10^{-4}$.
	
\begin{figure}[!ht]
		\centering
		\subfigure[]
		{
			\begin{minipage}{7cm}
				\centering
				\includegraphics[scale=0.36]{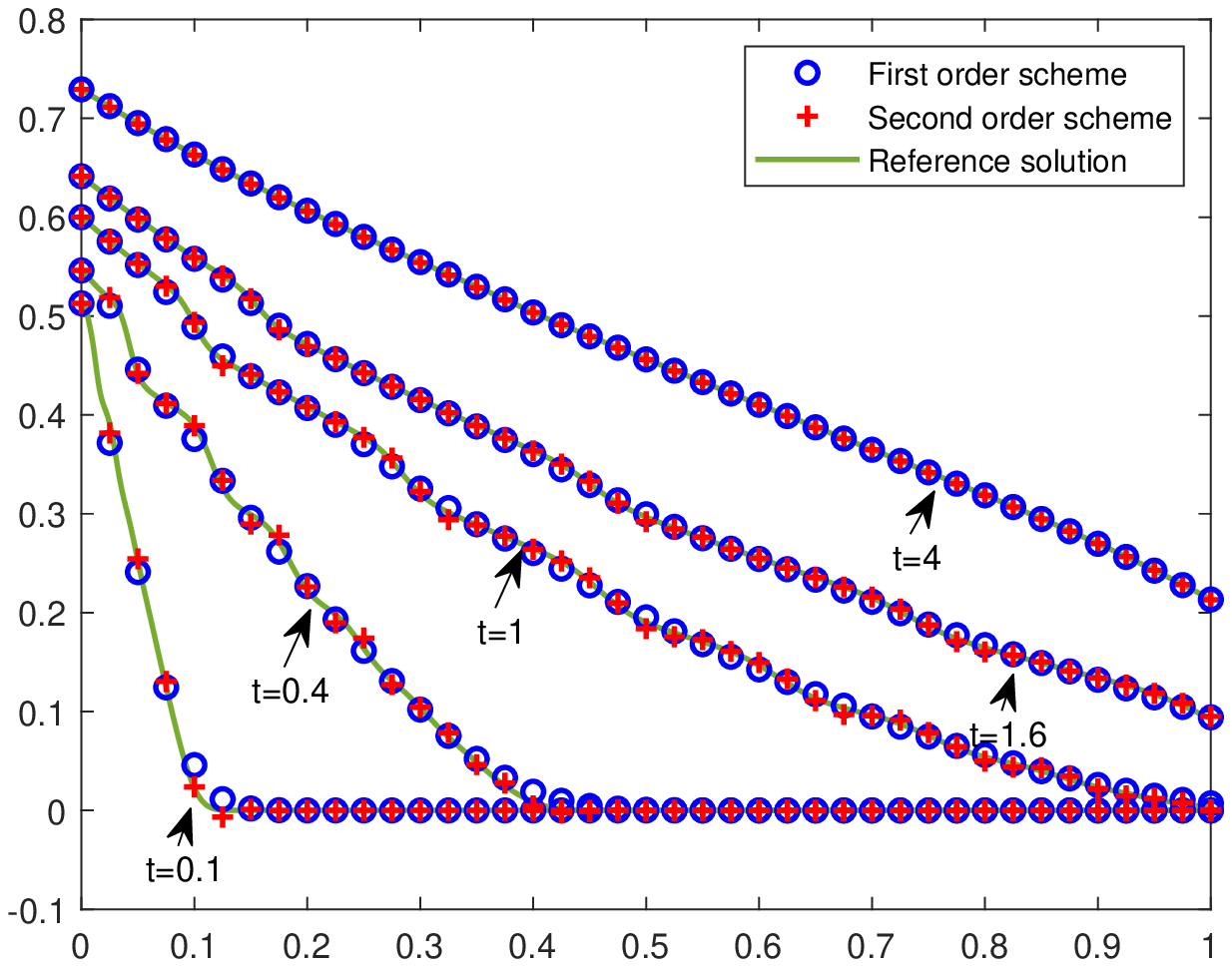}
			\end{minipage}
		}
		\subfigure[]
		{
			\begin{minipage}{7cm}
				\centering
				\includegraphics[scale=0.36]{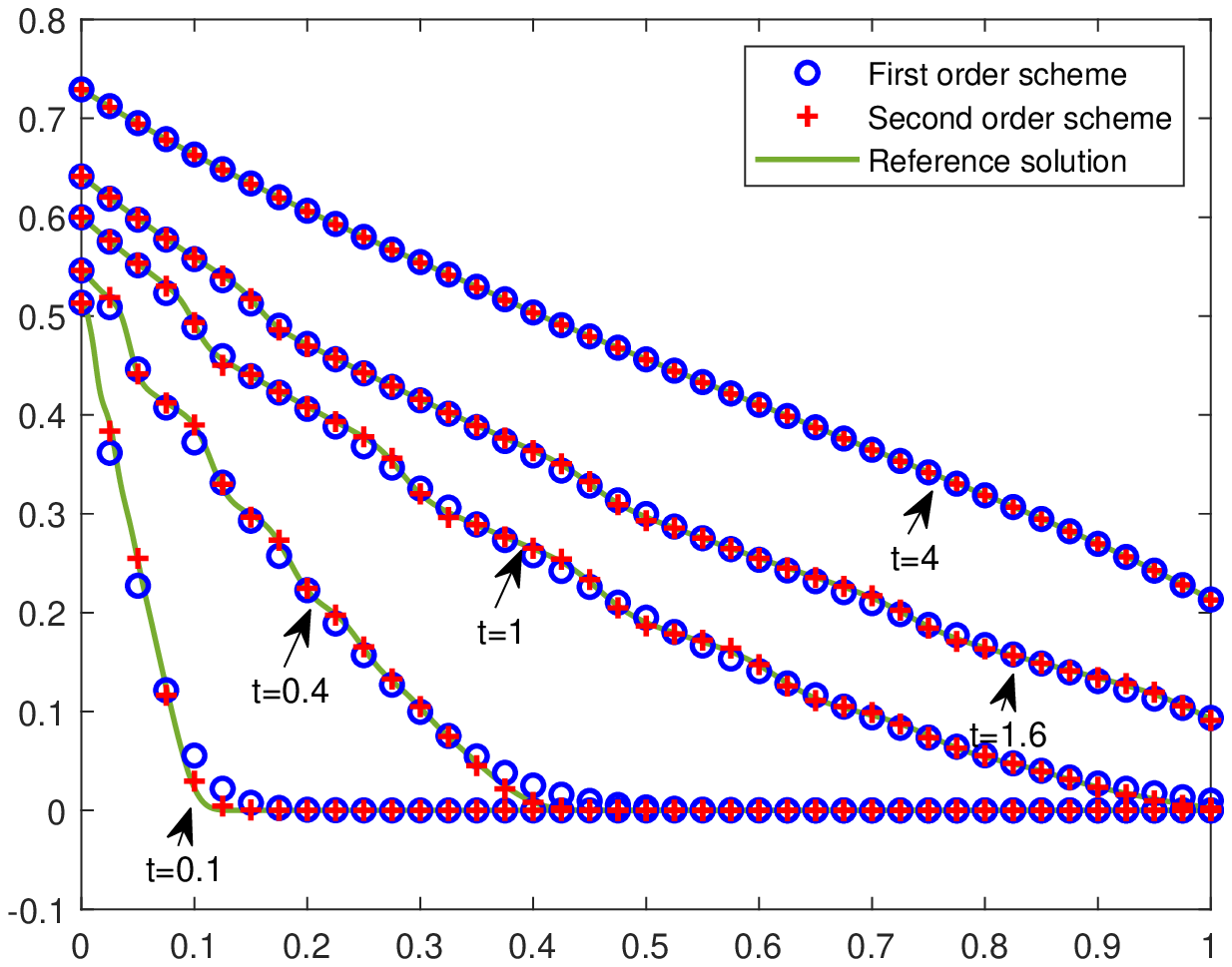}
			\end{minipage}
		}
		\caption{ Numerical solutions of $\rho$ for first and second order schemes in the kinetic regime $\varepsilon=1$ for the one-group transport equation with initial data \eqref{eq522}. (a): $\Delta t=0.4 \Delta x$, (b): $\Delta t=2\Delta x$.   }
		\label{One_kinetic}
	\end{figure}

	Fig. \ref{One_kinetic} shows the numerical solutions at $t = 0.1, 0.4, 1.0, 1.6$ and $4.0$, by taking $\varepsilon=1$ in the kinetic regime and $N=200$.
	The solid lines are the reference solutions obtained by the first order scheme using   $N=1000$ and $\Delta t=0.2\Delta x$. In Fig. \ref{One_kinetic}(a) and (b), we show the numerical solutions with $\Delta t=0.4\Delta x$ and $\Delta t=2\Delta x$ respectively. Both the results from first and second order schemes match the reference solution well. The second order results are slightly better than the first order ones. Besides, the differences between small and large time steps are almost negligible, which shows the higher efficiency of large time steps.
	
\begin{figure}[!ht]
		\centering
		\subfigure[]
		{
			\begin{minipage}{7cm}
				\centering
				\includegraphics[scale=0.36]{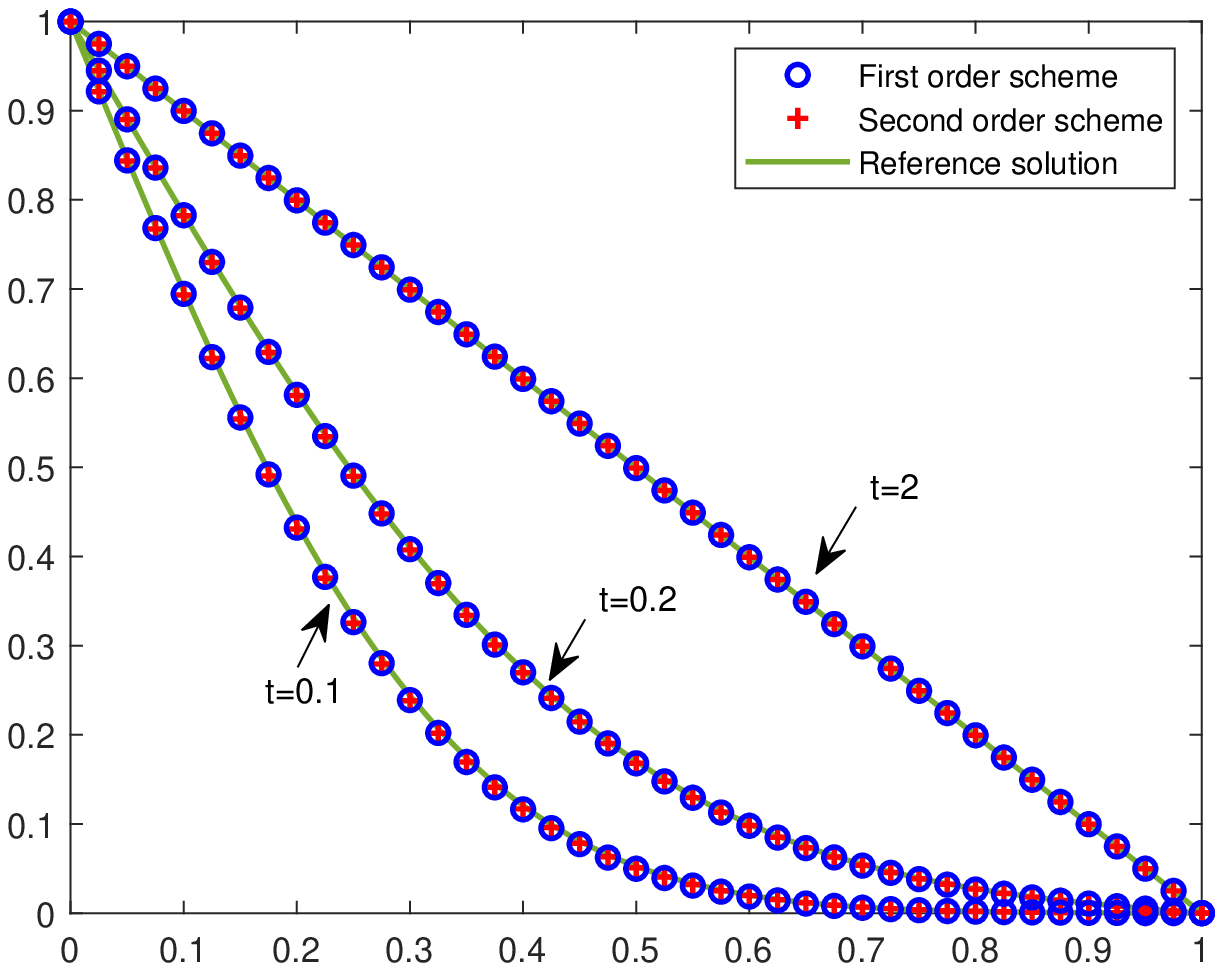}
			\end{minipage}
		}
		\subfigure[]
		{
			\begin{minipage}{7cm}
				\centering
				\includegraphics[scale=0.36]{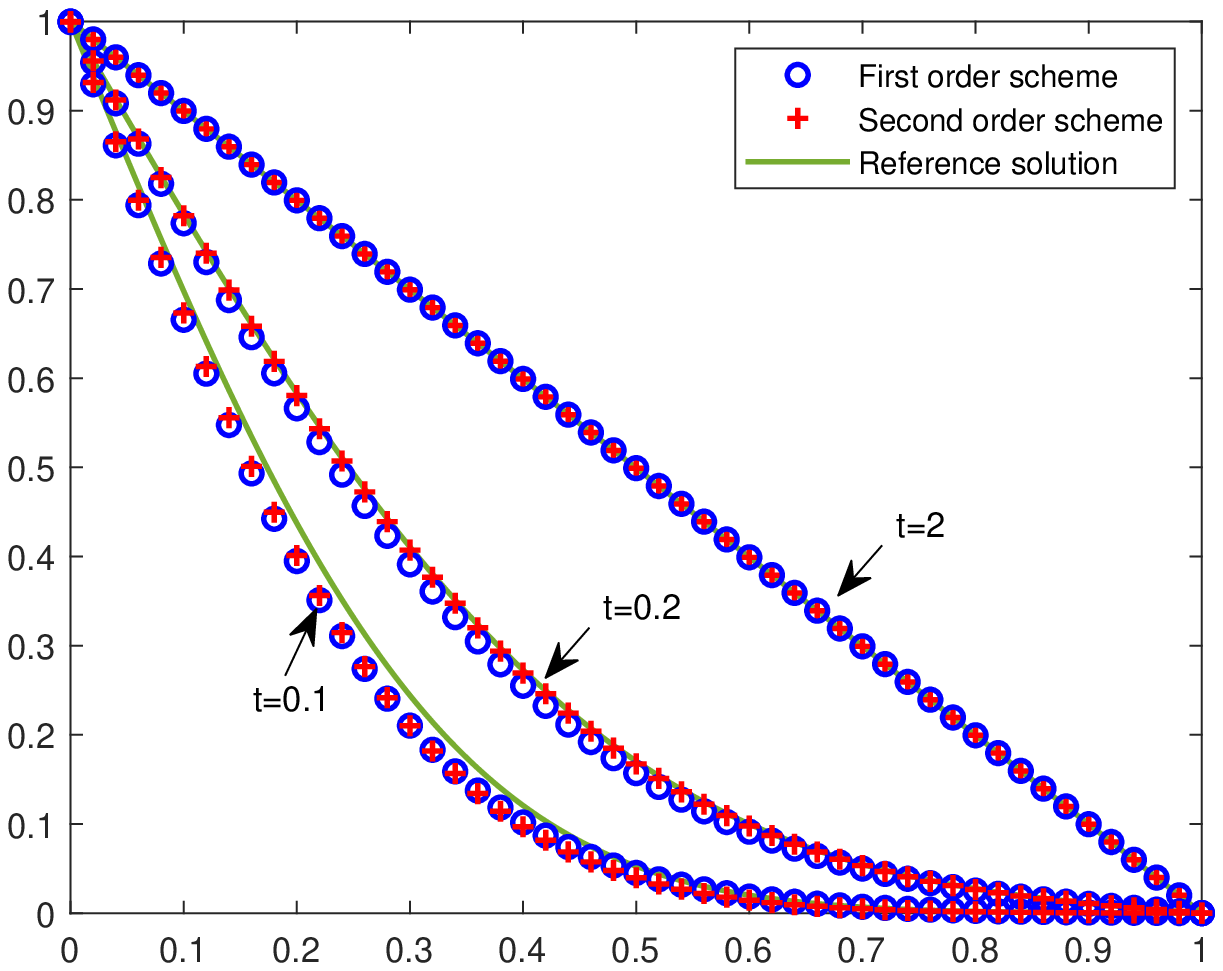}
			\end{minipage}
		}
		\caption{ Numerical solutions of $\rho$ for first and second order schemes in the diffusive regime $\varepsilon=10^{-4}$ for the one-group transport equation with initial data \eqref{eq522}. (a): $\Delta t=0.4 \Delta x$, (b): $\Delta t=2\Delta x$.   \label{One_parabolic} }
		
	\end{figure}

	Fig. \ref{One_parabolic} depicts the numerical solution at $t = 0.1, 0.2$ and $2.0$, by taking $\varepsilon =10^{-4}$ in the diffusive regime and   $N=200$. The solid lines are the  reference solutions obtained in the same way as above. We can see that the results from both first and second order schemes match the reference solution well. The second order scheme has slightly better results than the first order one when using a large time step. Both the first and second order results with a large time step are also very close to the reference solutions, which again shows the advantage of our schemes with large time steps. Additionally, all results are comparable to those in \cite{peng2021asymptotic,peng2020stability}.
	
	{\subsection{Two-dimension in space problems}
		\label{sec_numer_2d}
		In this section, we consider a 2D in space problem
		\begin{align}
			\label{two-1}
			\varepsilon \partial_{t} f+\ \mathbf{v}  \cdot  \nabla_{\mathbf{x}} f=\frac{\sigma_S}{\varepsilon}(\langle f\rangle-f)-\varepsilon \sigma_{A} f+\varepsilon G,
		\end{align}
		where $\langle f\rangle=\f{1}{4\pi}\int_{\mathbb{S}^2} f d\bv$, and the velocity $\bv:=(\xi,\eta,\gamma)=(\sin \theta \cos  \varphi, \sin \theta \sin  \varphi,\cos \theta )$ is confined on a unit sphere 
		$\mathbb{S}^2$, with $\theta$ and $\varphi$ being the polar and azimuthal angles respectively. Here the space is in two dimensions with $\bx=(x,y)\in \mathbb{R}^2$. $\sigma_S$ and $\sigma_{A}$ are the scattering and absorbing coefficients, and $G=G(\bx)$ is a source term, the same as in \eqref{one-group}. For simplicity, we take $\sigma_S$, $\sigma_A$ as constants and define $\mu= \frac{\sigma_S}{\varepsilon^2}+\sigma_A$,
		following the approach in Section \ref{sec_model_apprx}, we can similarly derive a new approximation model for \eqref{two-1}:
		\begin{equation}
			\label{model2D}
			\left\{
			\begin{array}{l}
				\varepsilon \partial_{t} f+\ \xi \partial_x f+\eta \partial_y f=\frac{\sigma_S}{\varepsilon}(\langle f\rangle-f)-\varepsilon \sigma_{A} f+\varepsilon G, \\ \, \\
				\rho_t(x,y,t)+\frac{e^{-\mu \left(t-t_{n}\right)}}{\varepsilon} \langle \xi(f-\rho)_x\left(\bX_n,\bY_n, \xi, \eta, \gamma, t_{n}\right)\rangle+\frac{e^{-\mu \left(t-t_{n}\right)}}{\varepsilon} \langle \eta(f-\rho)_y\left(\bX_n,\bY_n, \xi,\eta,\gamma, t_{n}\right)\rangle  \\ \, \\
				\hspace{1.5cm}-\langle \xi^2\rangle \frac{1}{\sigma_S+\eps^2\,\sigma_A} (1-e^{-\mu(t-t_n)})\rho_{xx}\left(x,y,t\right)-\langle \eta^2\rangle \frac{1}{\sigma_S+\eps^2\,\sigma_A} (1-e^{-\mu(t-t_n)})\rho_{yy}\left(x,y,t\right) =0,
			\end{array}
			\right.
		\end{equation}
		where we denote $\bX_n=x-(t-t_n)\xi/\eps $ and $\bY_n=y-(t-t_n)\eta/\eps$.
		For this 2D problem, a DOM with $86$ Lebedev quadrature points is used for the angular variables $\theta$ and $\varphi$ on the unit sphere $\mathbb{S}^2$ \cite{sphere}. The first and second order schemes proposed in Section \ref{sec_scheme} can be easily extended to the 2D in space problems here.}
	
	{\subsubsection{Accuracy test}
		We will first verify the convergence orders of our schemes with an exact solution
		\begin{align}
			&f(t, x, y, \xi, \eta, \gamma)=\exp (-t) \sin ^2(2 \pi x) \sin ^2(2 \pi y)\left(1+\varepsilon\left(\frac{\eta+\eta^3}{3}\right)\right), \quad(x, y)\in \Omega_{\bx} =[0,1]^2,\nonumber\\
			&\rho(t, x, y)=\exp (-t) \sin ^2(2 \pi x) \sin ^2(2 \pi y), \nonumber
		\end{align}
		where we take the coefficients to be $\sigma_S=1, \sigma_A=0$, with a source term $G$ given by
		$$
		G(t, x, y, \xi, \eta,\gamma)=\partial_t f+\frac{1}{\varepsilon}\left( \xi  \nabla_x f+\eta  \nabla_y f\right)+\frac{1}{\varepsilon} \left( \exp (-t) \sin ^2(2 \pi x) \sin ^2(2 \pi y)\left(\frac{\eta+\eta^3}{3}\right)  \right).
		$$
		We consider periodic boundary conditions along the $x$ and $y$ directions,  and compute the solution up to a final time $T=1$. In Table \ref{two_t1} and Table \ref{two_t2}, we show the $L^1$ and $L^\infty$ errors and corresponding orders with a large time step $\Delta t=3\Delta x$. $f(x,y,\xi_1,\eta_1,\gamma_1,T)$ is the solution at the first quadrature point. One easily finds that our scheme can reach the desired orders and capture the correct asymptotic limit. Besides, our schemes are fine for a uniformly large time step.}
	\begin{table}[!ht]
		\scriptsize
		\centering
		\caption{$L^{\infty}$ and $L^{1}$  errors and convergence orders for $\rho$   and $f(x,y,\xi_1,\eta_1,\gamma_1,T)$.  First order scheme for 2D in space accuracy test.}
		\label{two_t1}
		\begin{tabular}{|c|c|c|c|c|c|c|c|c|c|} \hline	
			$\eps$ &N  & $L^{\infty}$ error of $\rho$ & Order & $L^{\infty}$ error of $f$ & Order  & $L^{1}$ error of $\rho$ & Order & $L^{1}$ error of $f$ & Order        \\ \hline
			\multirow{4}{*}{ 1}         	
			
			& 8 & 7.19E-2 & -- &  2.21E-1& --  &5.01E-2 & -- &  6.61E-2& --                        \\ \cline{2-10}
			& 16 & 4.71E-2 & 0.61 &  9.86E-2& 1.17  & 2.12E-2 & 1.24 &3.13E-2& 1.08            \\ \cline{2-10}
			& 32 & 2.64E-2 & 0.83 &  5.12E-2& 0.95 & 1.03E-2 & 1.05 &1.43E-2& 1.13    \\ \cline{2-10}
			& 64 & 1.47E-2 & 0.85 &  2.51E-2& 1.03   &5.50E-2 & 0.90 &6.79E-2& 1.07     \\   \hline
			\multirow{4}{*}{ $10^{-2}$}
			& 8 & 8.65E-2 & -- &  8.66E-2& --  &2.47E-2 & -- &  2.75E-2& --                        \\ \cline{2-10}
			& 16 & 2.03E-2 & 2.09 & 2.15E-2 & 2.01  & 1.52E-2 & 0.70 & 1.52E-2 & 0.86          \\ \cline{2-10}
			& 32 & 8.43E-3 & 1.23 &  1.07E-2 & 1.01 & 7.46E-3 & 1.03 & 7.46E-3 & 1.03     \\ \cline{2-10}
			& 64 & 4.65E-3 & 0.86 &  5.81E-3 & 0.89   &3.77E-3 & 0.99 & 3.77E-3 & 0.99     \\  \hline
			\multirow{4}{*}{$ 10^{-6}$}
			& 8 & 1.09E-1 & -- &  1.09E-1& --     &2.95E-2 & -- &  2.95E-2& --                        \\ \cline{2-10}
			& 16 & 3.02E-2 & 1.85 & 3.02E-2 & 1.85  &1.49E-3 & 0.99 & 1.49E-3 & 0.99          \\ \cline{2-10}
			& 32 & 1.11E-2 & 1.44 & 1.11E-2 & 1.44 & 7.38E-3 & 1.01 &7.38E-3 & 1.01    \\ \cline{2-10}
			& 64 & 4.67E-3 & 1.25 & 4.67E-3 & 1.25   &3.75E-3 & 0.98 &3.75E-3 & 0.98    \\    \hline
		\end{tabular}
	\end{table}

	\begin{table}[!ht]
		\scriptsize
		\centering
		\caption{$L^{\infty}$ and $L^{1}$  errors and convergence orders for $\rho$   and $f(x,y,\xi_1,\eta_1,\gamma_1,T)$.  Second order scheme for 2D in space accuracy test.}
		\label{two_t2}
		\begin{tabular}{|c|c|c|c|c|c|c|c|c|c|} \hline	
			$\eps$ &N  & $L^{\infty}$ error of $\rho$ & Order & $L^{\infty}$ error of $f$ & Order  & $L^{1}$ error of $\rho$ & Order & $L^{1}$ error of $f$ & Order        \\ \hline
			\multirow{4}{*}{ 1}         	
			& 8 & 5.04E-2 & -- &  1.22E-1& --  &4.13E-2 & -- &  6.27E-2& --                        \\ \cline{2-10}
			& 16 & 1.59E-2 & 1.66 & 3.65E-2 & 1.74  & 1.26E-2 & 1.71 & 1.84E-2 & 1.76         \\ \cline{2-10}
			& 32 & 4.53E-3 & 1.81 &  9.62E-3 & 1.92 &  3.21E-3 & 1.97  & 4.65E-3 & 1.98    \\ \cline{2-10}
			& 64 & 1.11E-3 & 2.03 &  2.38E-3 & 2.02  & 8.04E-4 & 2.00 & 1.14E-3 & 2.03    \\   \hline
			\multirow{4}{*}{ $10^{-2}$}
			& 8 & 7.41E-2 & -- &  7.44E-2& --  &2.62E-2 & -- &  2.68E-2& --                        \\ \cline{2-10}
			& 16 & 1.64E-2 & 2.17 & 1.64E-2 & 2.18  & 5.12E-3 & 2.36 & 5.12E-3 & 2.39          \\ \cline{2-10}
			& 32 & 4.01E-3 & 2.03 &  4.03E-3 & 2.02 & 1.22E-3 & 2.07 & 1.22E-3 & 2.07     \\ \cline{2-10}
			& 64 & 1.03E-3 & 1.96 &  1.04E-3 & 1.95  &3.13E-4 & 1.97 & 3.12E-4 & 1.97   \\   \hline
			\multirow{4}{*}{$ 10^{-6}$}
			& 8 & 1.04E-1 & -- &  1.04E-1& --  &3.12E-2 & -- &  3.12E-2& --                        \\ \cline{2-10}
			& 16 & 2.03E-2 & 2.35 & 2.03E-2 & 2.35  & 5.02E-3 & 2.64 & 5.02E-3 & 2.64          \\ \cline{2-10}
			& 32 & 5.02E-3 & 2.02 &  5.02E-3 & 2.02 & 1.22E-3 & 2.04 & 1.22E-3 & 2.04     \\ \cline{2-10}
			& 64 & 1.21E-3 & 2.03 &  1.21E-3 & 2.03  &2.97E-4 & 2.04 & 2.97E-4 & 2.04    \\   \hline
		\end{tabular}
	\end{table}

	{\subsubsection{Gaussian initial value}
		We then consider a Gaussian initial value \cite{einkemmer2021asymptotic}
		\begin{equation}\label{Gau}
			f(t=0, x, y, \xi, \eta, \gamma)=\frac{1}{4 \pi \varsigma^2} \exp \left(-\frac{x^2+y^2}{4 \varsigma^2}\right), \quad \varsigma^2=10^{-2}, \quad(x, y) \in[-1,1]^2,
		\end{equation}
		and the parameters $\sigma_A=G=0$. We take two different scattering coefficients as follows.}
	
\begin{figure}[!ht]
		\centering
		\subfigure[]
		{
			\begin{minipage}{5cm}
				\centering
				\includegraphics[scale=0.3]{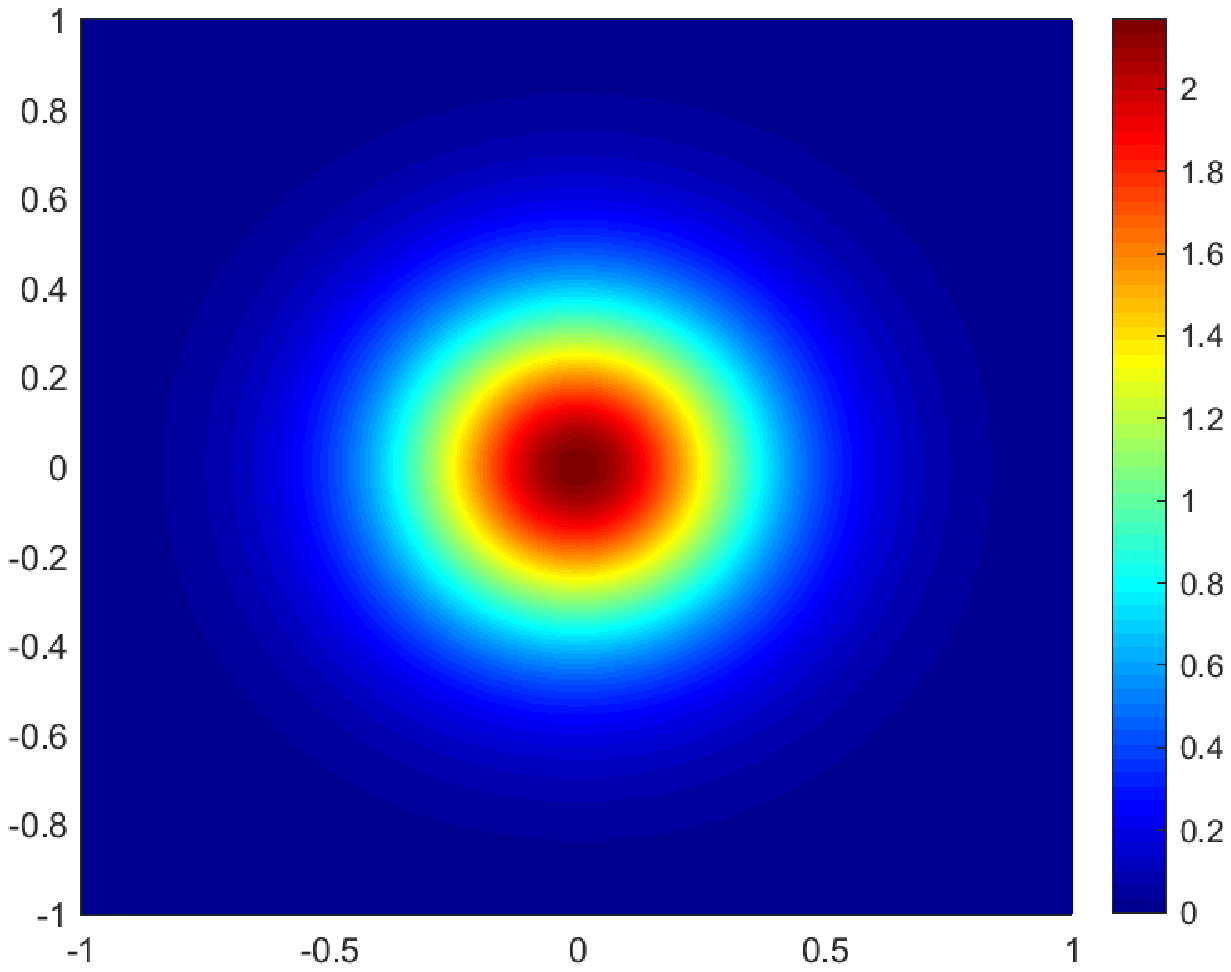}
			\end{minipage}
		}
		\subfigure[]
		{
			\begin{minipage}{5cm}
				\centering
				\includegraphics[scale=0.3]{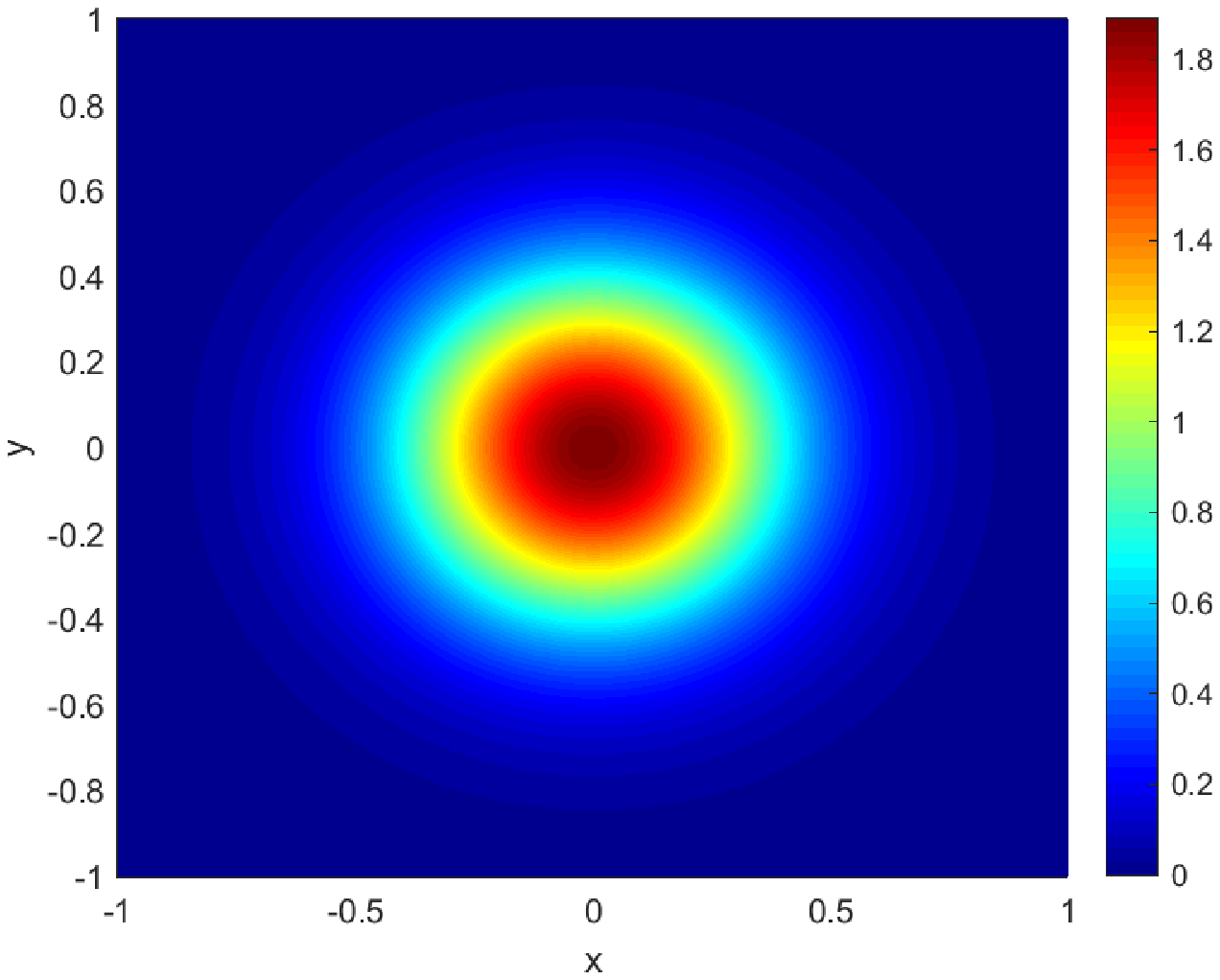}
			\end{minipage}
		}
		\subfigure[]
		{
			\begin{minipage}{5cm}
				\centering
				\includegraphics[scale=0.3]{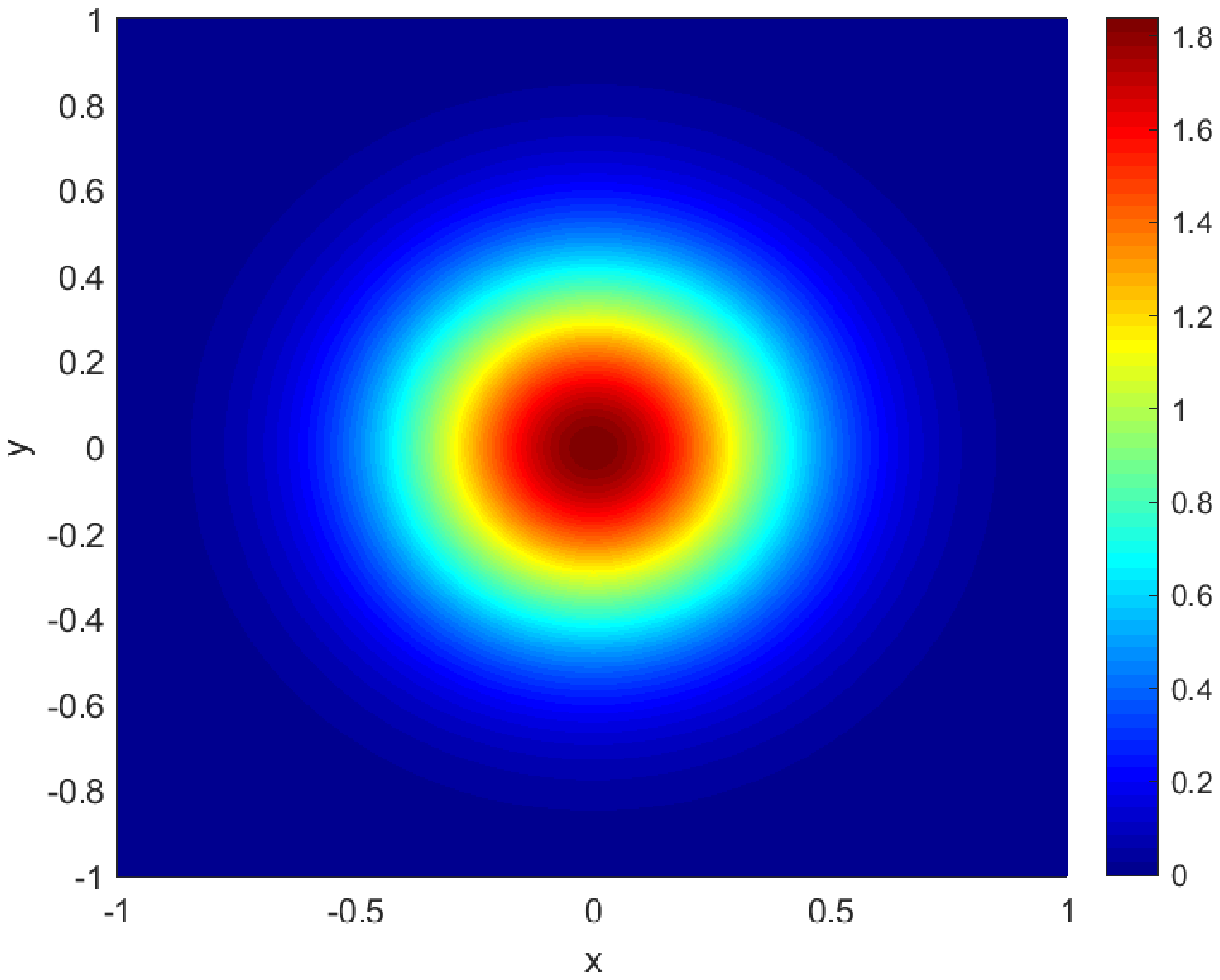}
			\end{minipage}
		}
		\caption{ Numerical solutions of $\rho$  in the diffusive regime $\varepsilon=10^{-6}$ with $\Delta t=2\Delta x$. (a): first order scheme, (b): second order scheme,  (c): reference solution.  }
		\label{two_G_1}
	\end{figure}

\begin{figure}[!ht]
		\centering
		\subfigure[]
		{
			\begin{minipage}{7cm}
				\centering
				\includegraphics[scale=0.36]{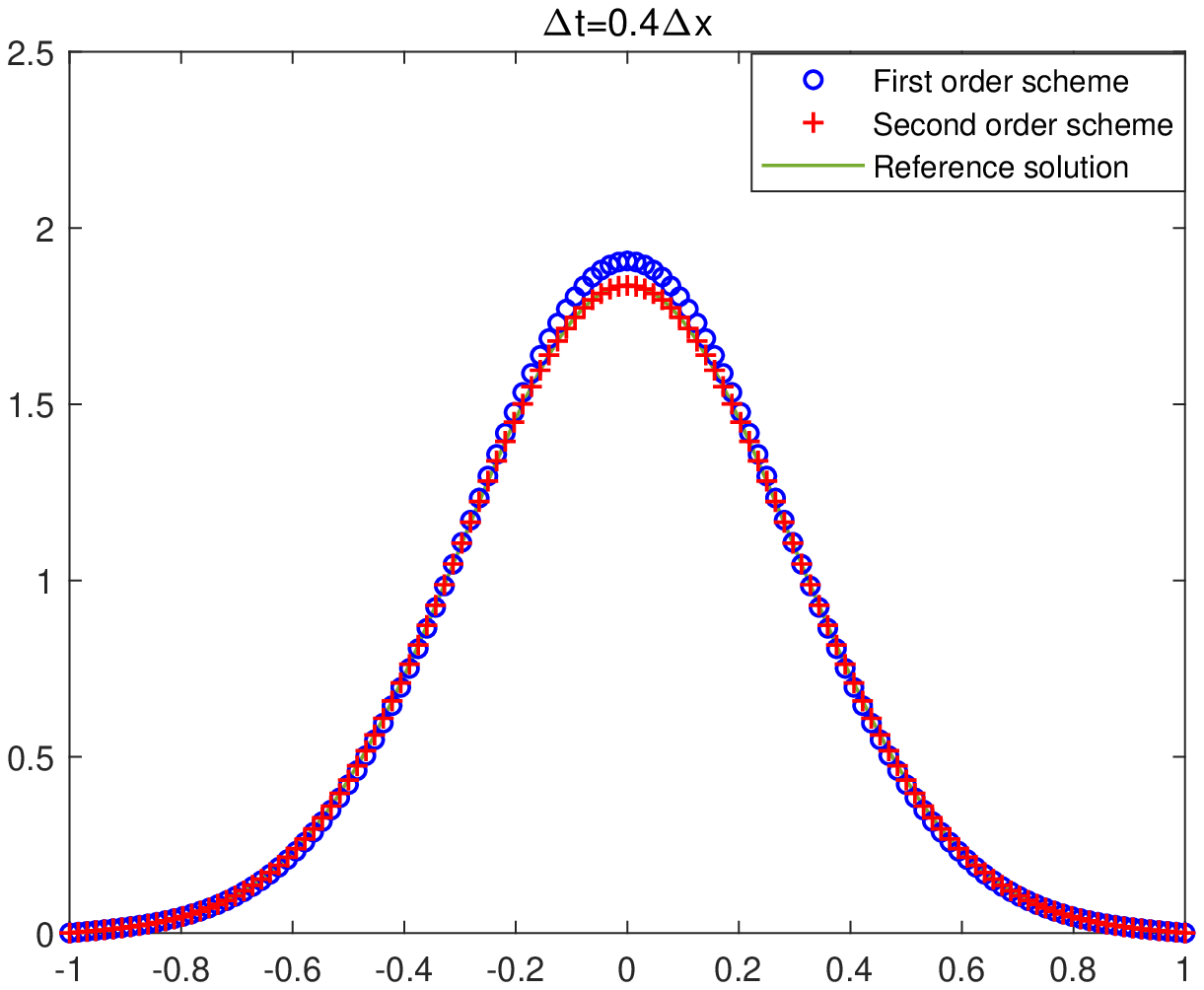}
			\end{minipage}
		}
		\subfigure[]
		{
			\begin{minipage}{7cm}
				\centering
				\includegraphics[scale=0.36]{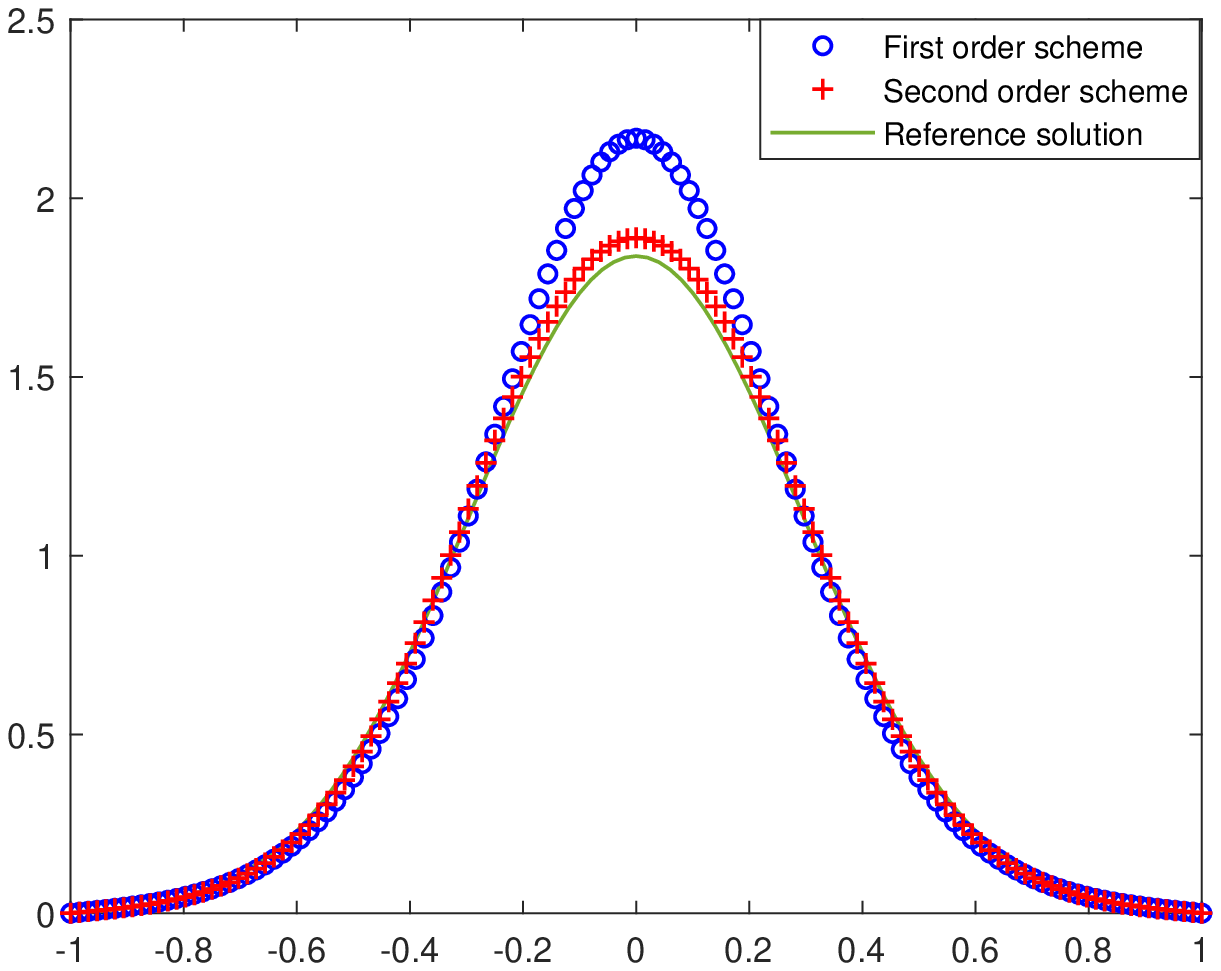}
			\end{minipage}
		}
		\caption{  Comparison of first and second order results with a reference solution along $y = 0$. (a): $\Delta t=0.4 \Delta x$, (b): $\Delta t=2\Delta x$.   }
		\label{two_G_2}
	\end{figure}

	{{\bf Case (a): Constant scattering coefficient $\sigma_S$.}
		Firstly, we consider a constant coefficient $\sigma_S=1$. We take a large time step $\Delta t= 2\Delta x$, with $N=128$ along each spatial direction. In Fig. \ref{two_G_1}, we show the numerical density $\rho$ in the diffusive regime $\varepsilon =10^{-6}$ at $T = 0.1$. (a) and (b) are the solutions obtained by first and second order schemes, respectively. (c) is a reference solution computed for its diffusive limiting equation by a backward Euler method with a second order central difference discretization, using a small time step $\Delta t= 0.001\Delta x$ and $N=256$ along each spatial direction. We can see both first and second order results with a large time step are close to the reference solution. In Fig. \ref{two_G_2}, we show the cutting plots along $y = 0$ for these results, a smaller time step $\Delta t= 0.4\Delta x$ is also considered for comparison. We can find that with a small time step $\Delta t= 0.4\Delta x$, both first and second order results are almost the same as the reference solution. For a large time step $\Delta t= 2\Delta x$, the first order scheme has a relatively large error, while the second order scheme still matches the reference solution. From this example, we can see that our schemes also work well in the high-dimensional case with large time steps, only large errors may occur for a first order scheme.}
	
\begin{figure}[!ht]
		\centering
		\includegraphics[scale=0.3]{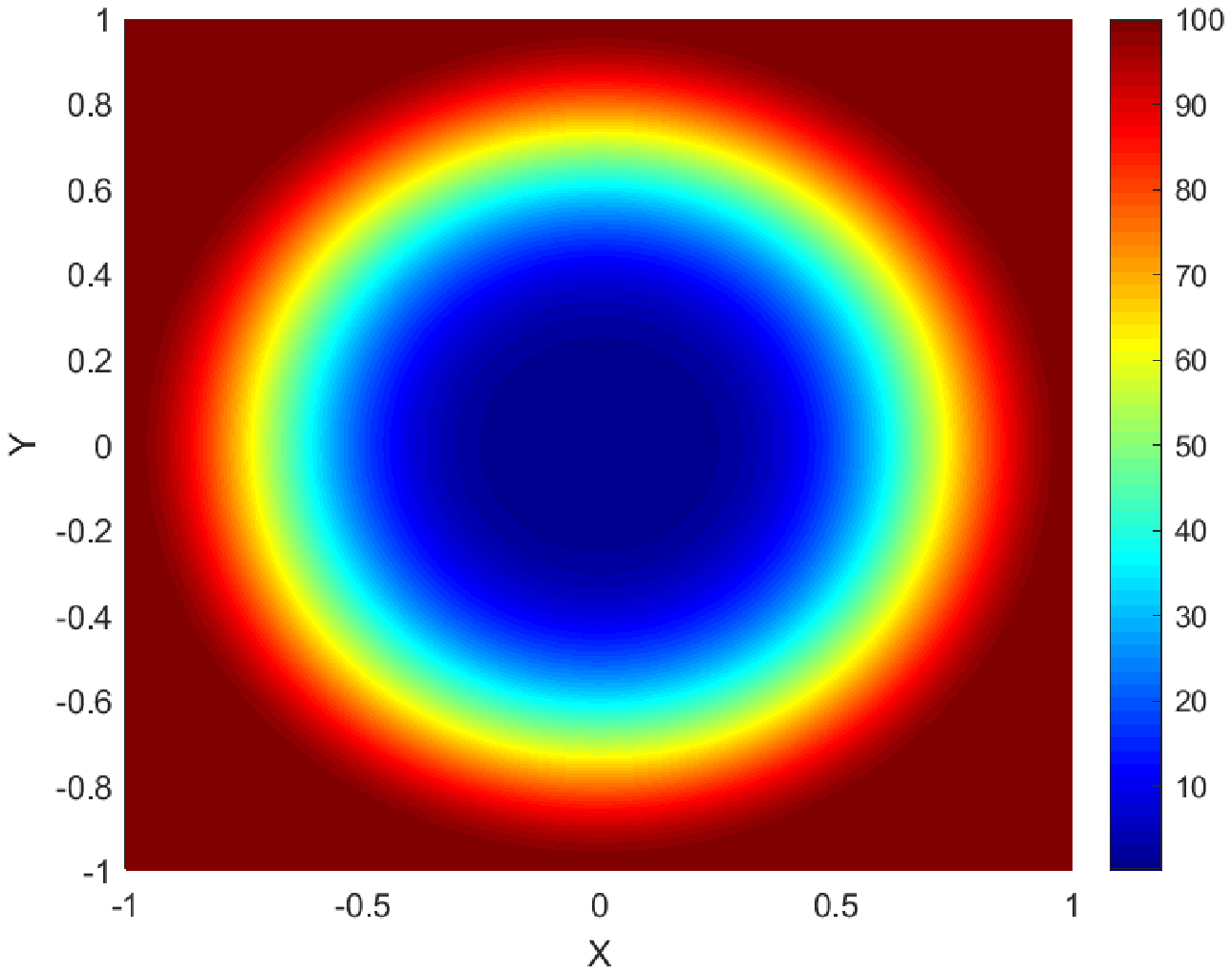}
		\includegraphics[scale=0.3]{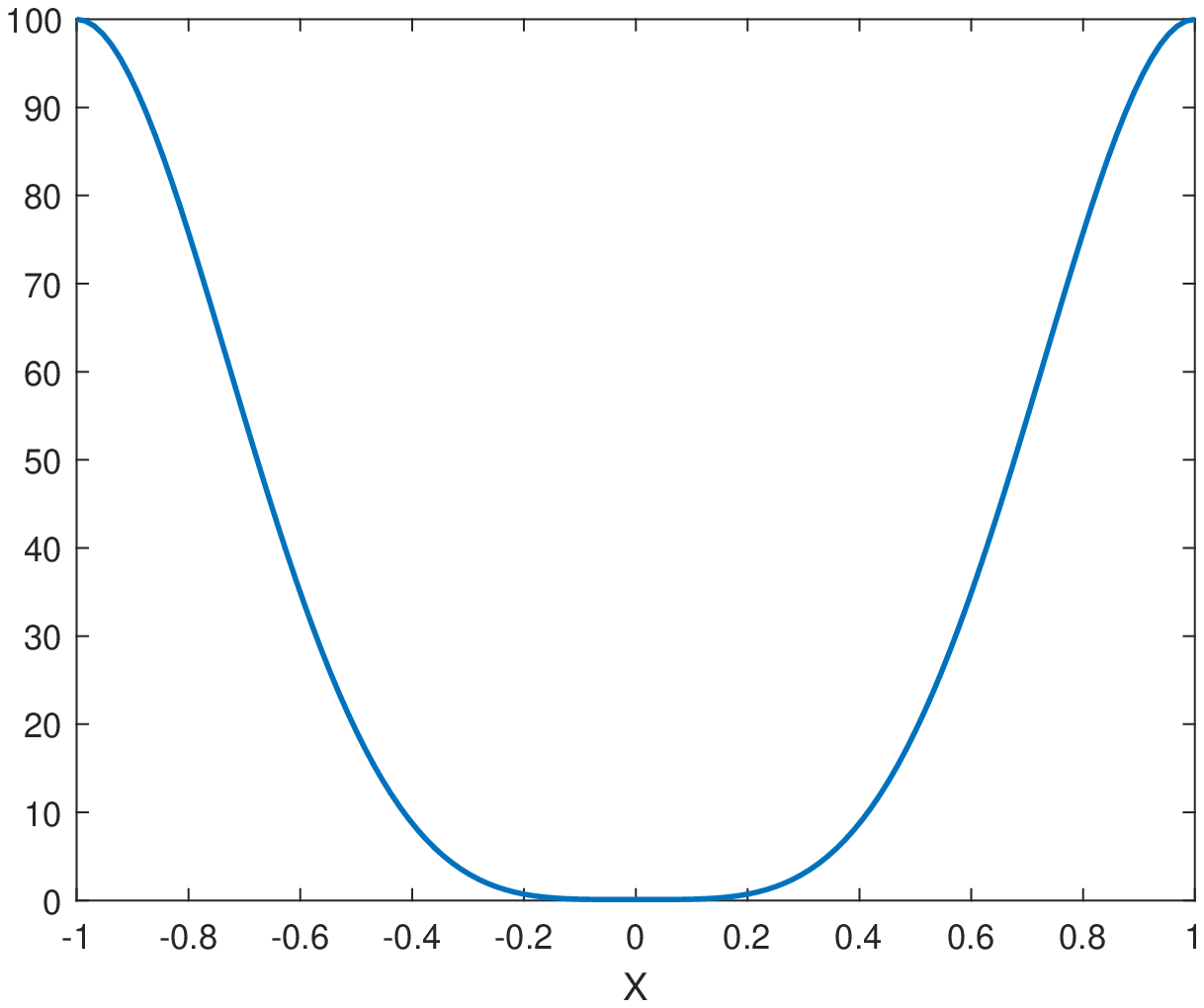}
		\caption{The profile of $\frac{\sigma_S}{\varepsilon}$ (left) and a cutting along $y=0$ (right).}
		\label{two_G_3}
	\end{figure}

\begin{figure}[!ht]
		\centering
		\subfigure[]
		{
			\begin{minipage}{7cm}
				\centering
				\includegraphics[scale=0.36]{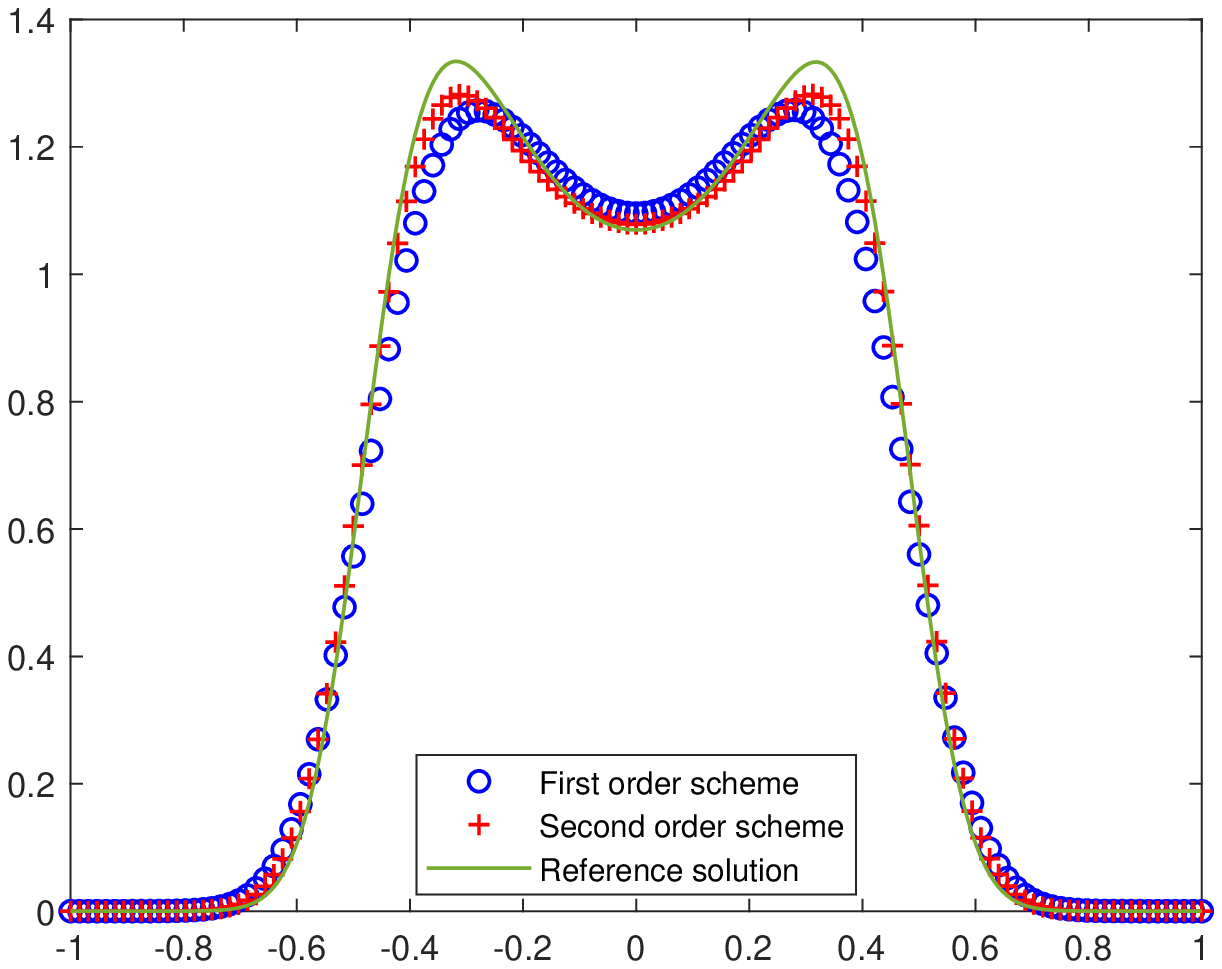}
			\end{minipage}
		}
		\subfigure[]
		{
			\begin{minipage}{7cm}
				\centering
				\includegraphics[scale=0.36]{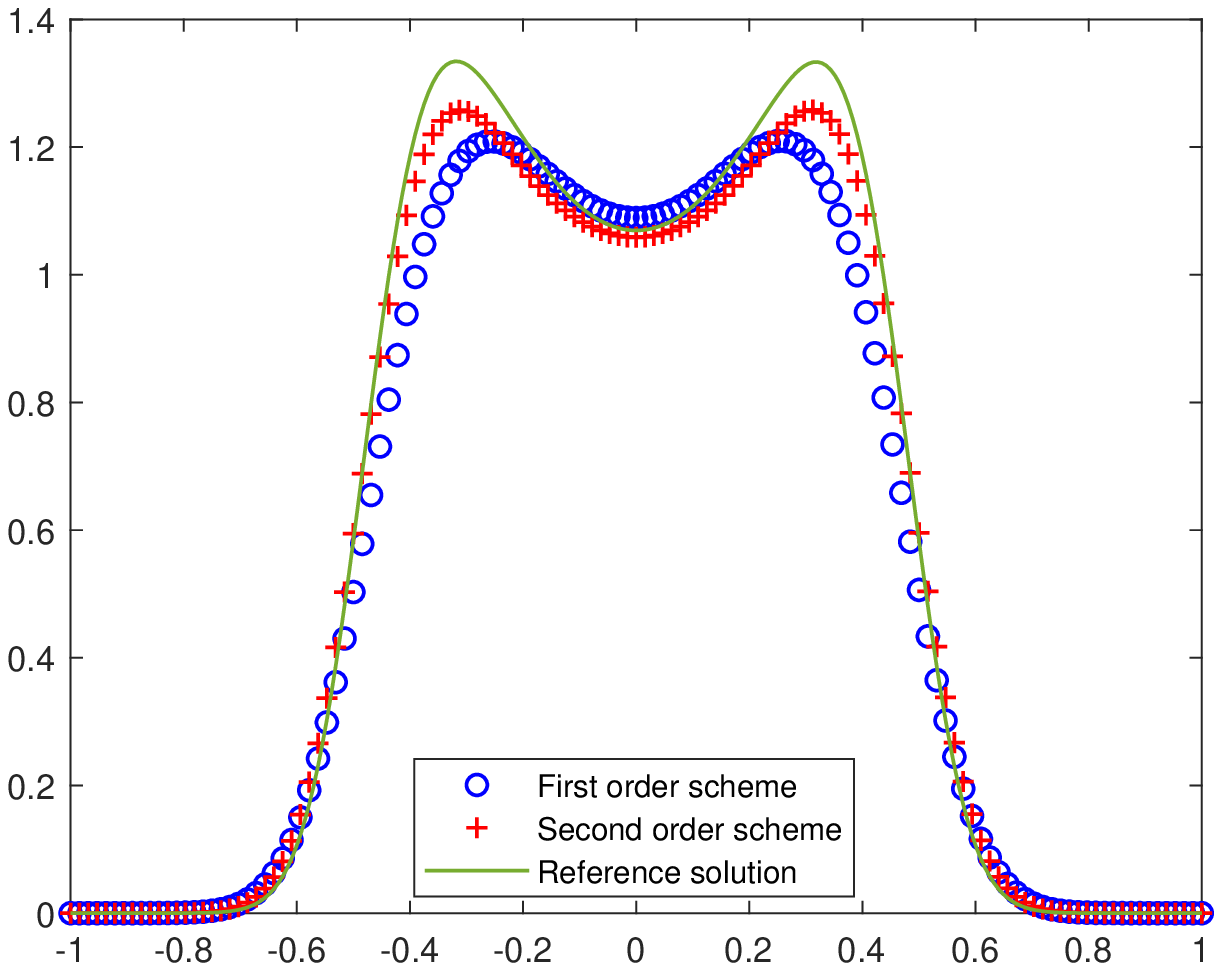}
			\end{minipage}
		}
		\caption{Comparison of first and second order results with a reference solution along $y = 0$. (a): $\Delta t=0.01 \Delta x$, (b): $\Delta t=0.04\Delta x$.   }
		\label{two_G_5}
	\end{figure}

	{{\bf{Case (b): Variable scattering coefficient $\sigma_S$.}}
		Secondly we consider a very challenging case with a variable scattering coefficient $\sigma_S$ \cite{einkemmer2021asymptotic}, which is given by
		\begin{equation}
			\sigma_S(x, y)= \begin{cases}0.999 c^4( c+\sqrt{2} )^2 ( c-\sqrt{2} )^2+0.001, & c=\sqrt{x^2+y^2}<1, \\ 1, & \text { otherwise, }\end{cases}
		\end{equation}
		and we take $\varepsilon=0.01$. We draw the profile of $\frac{\sigma_S}{\varepsilon}$ in Fig. \ref{two_G_3}. We can see that it varies across a large range from $0.1$ to $100$, so that both kinetic and diffusive regimes are included. In Fig. \ref{two_G_5}, we show a comparison of the first and second order results with a reference solution along $y = 0$ at $T=0.006$. Similarly we take $N=128$ along each spatial direction for the first and second order schemes. The reference solution is computed by the first order scheme with a finer mesh $N=400$ and a much smaller time step $\Delta t= 0.0001\Delta x$. Similarly, the first and second results match the reference solution well. A large time step still has very good performance, only with a little larger numerical errors.}
	
\begin{figure}[!ht]
		\centering
		\subfigure[]
		{
			\begin{minipage}{7cm}
				\centering
				\includegraphics[scale=0.36]{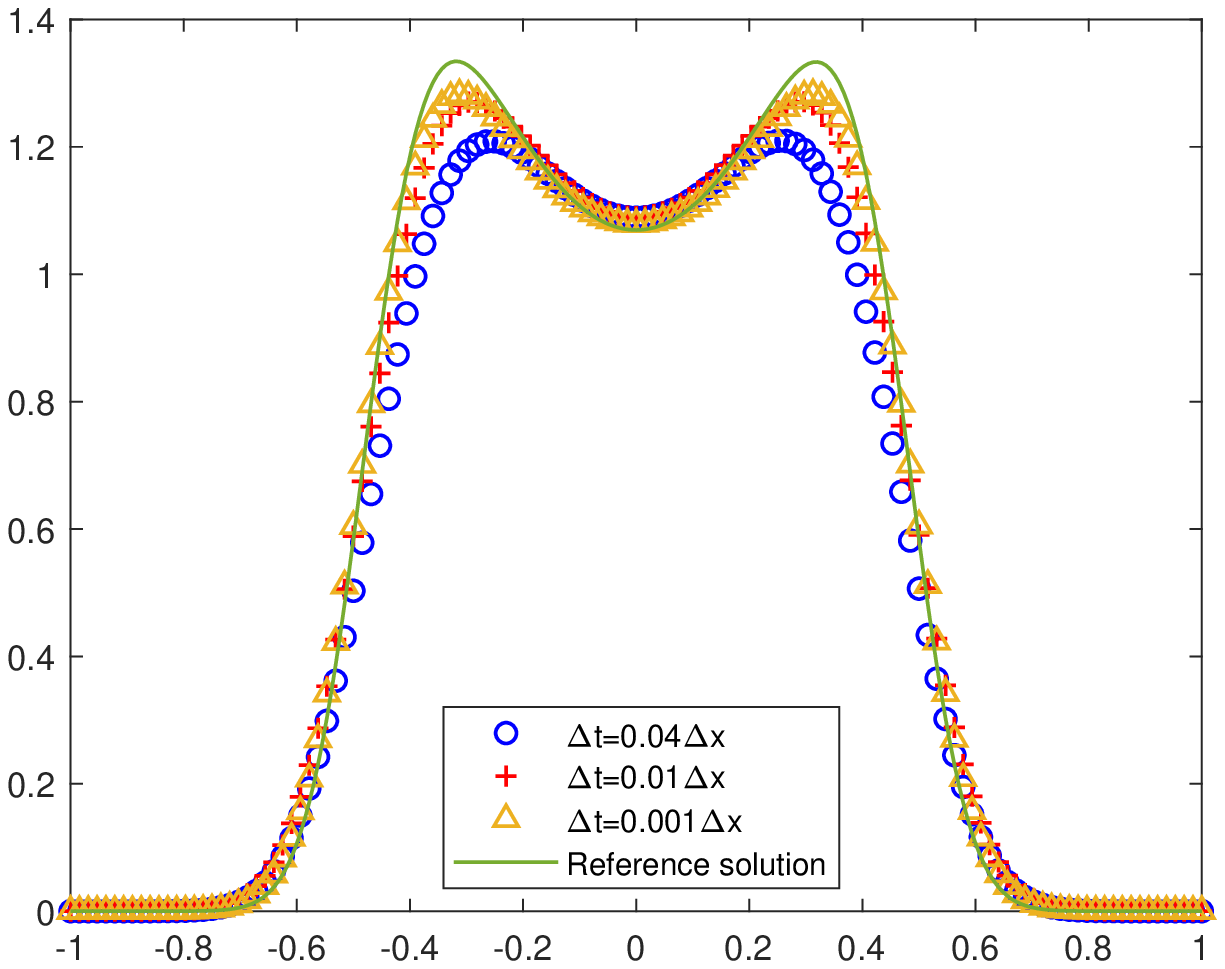}
			\end{minipage}
		}
		\subfigure[]
		{
			\begin{minipage}{7cm}
				\centering
				\includegraphics[scale=0.36]{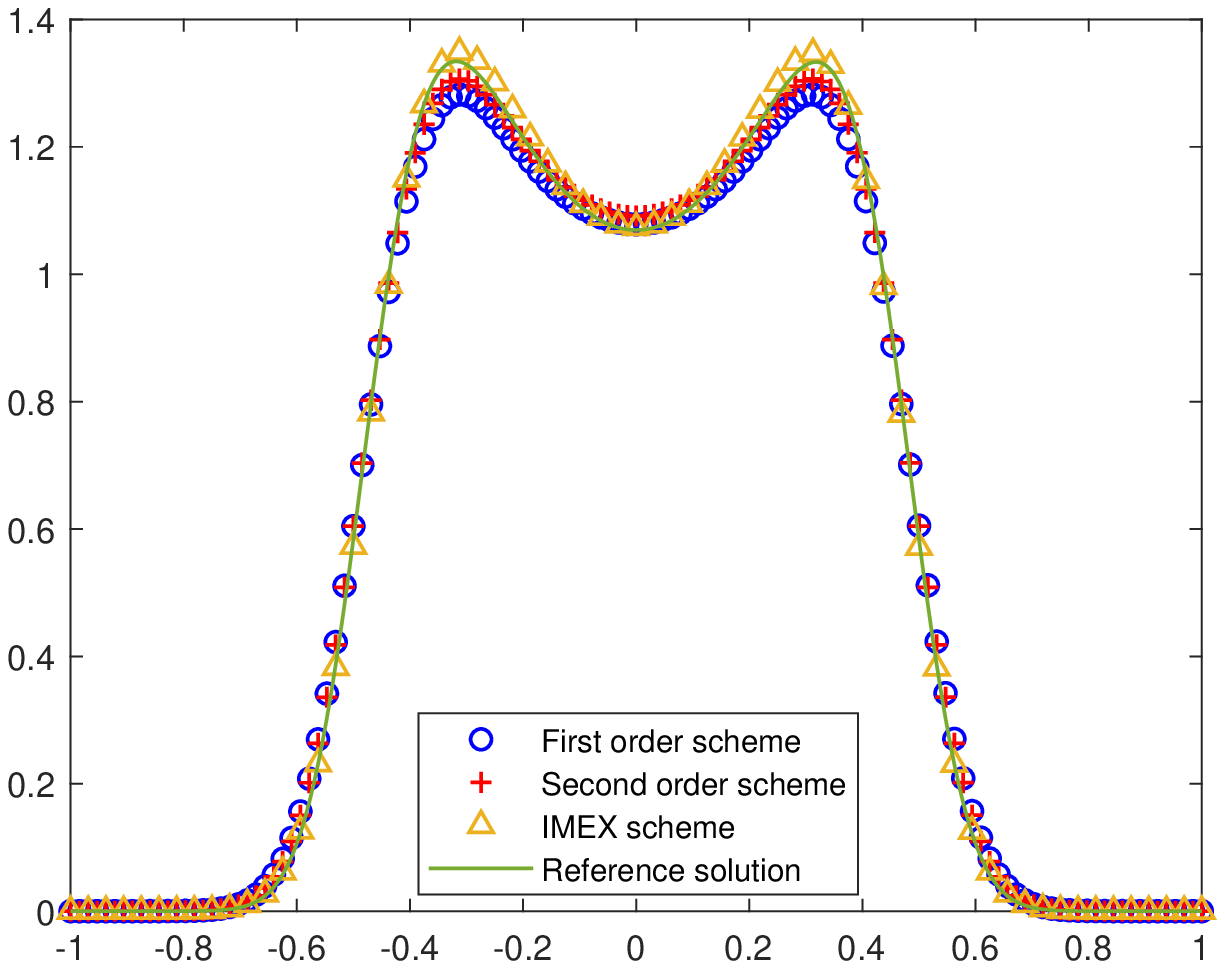}
			\end{minipage}
		}
		\caption{(a): first order scheme with $\Delta t=0.04 \Delta x$, $\Delta t=0.01 \Delta x$, and $\Delta t=0.001 \Delta x$. (b):  comparison of first and second order schemes with the IMEX scheme, $\Delta t= 0.001\Delta x$ for all three schemes. $N=128$ along each spatial direction.  }
		\label{two_G_6}
	\end{figure}

	{We also compare the results and CPU cost of our schemes, with a first order IMEX method  based on a micro-macro decomposition framework as in \cite{lemou2008new}. The IMEX scheme has a time step restriction $\Delta t=0.1 ~ \text{min}(\sigma_S)\Delta x^2+0.1\varepsilon\Delta x$. For this example, it is $\Delta t=0.001\Delta x$, namely in the order of $\mathcal{O}(\eps h)$.
		In Fig. \ref{two_G_6}(a), we show that for our first order scheme, if we take smaller and smaller time steps, such as $\Delta t=0.04\Delta x, \Delta t=0.01\Delta x, \Delta t=0.001\Delta x$, it is getting closer to the reference solution. In Fig. \ref{two_G_6}(b), we compare our schemes with the IMEX scheme using the same mesh size $N=128$ and the same time step $\Delta t= 0.001\Delta x$. It can be observed that all results match the reference solution. However, our schemes allow large time steps, as shown in Fig. \ref{two_G_5}. Finally we compare the CPU cost for these three schemes. We list them in Table. \ref{Table_twod}. With the same time step $\Delta t=0.001\Delta x$, the CPU cost is $6352s$ for our first order scheme, $22249s$ for our second order scheme, while $4323s$ for the IMEX method from \cite{lemou2008new}. However, if we take $\Delta t=0.01\Delta x$ and $\Delta t=0.04\Delta x$, the CPU cost of our first order scheme drops to $621s$ and $157s$, respectively, while for the second order scheme, it becomes $2217s$ and $542s$, respectively. We can see that with allowed large time steps, our schemes are much more efficient, especially the second order scheme has comparable results but much less CPU cost.}
	
	\begin{table}[htbp]
		\scriptsize
		\centering
		\caption{ Comparison of CPU cost (seconds) for three different schemes. \label{Table_twod}}
		\begin{tabular}{|c|c|c|c|}
			\hline	     	
			$\Delta t$  & first order scheme  & second order scheme & the IMEX scheme \\ \hline
			$0.04$  & $157$  & $542$    & --      \\ \hline
			$0.01$  & $621$  & $2217$   & --    \\ \hline
			$0.001$ & $6352$ & $22249$  & $4323$ \\ \hline
		\end{tabular}
	\end{table}

	\section{Conclusions}
	\label{sec_conc}
	\setcounter{equation}{0}
	\setcounter{figure}{0}
	\setcounter{table}{0}
	
	In this paper, we have proposed uniformly unconditionally stable first and second order finite difference schemes for kinetic transport equations in the diffusive scaling. We first derived an approximation model based on a formal solution of the original equation. The approximation error was analyzed. Then, first and second order schemes both in time and in space based on the approximation model were designed. For the developed schemes, their uniformly unconditional stability has been validated by a Fourier analysis, which is consistent with our numerical results. Numerical experiments have demonstrated the first and second orders of accuracy in both the kinetic and diffusive regimes, uniform stability under a large time step condition and good performances of our proposed schemes. {Order reductions are observed in the intermediate regime when the time step size is comparable to the parameter $\eps$, which is due to a modeling error and it is also a typical issue for  multi-scale type problems \cite{xiong2021high}.} The applicability to high-dimensional problems and its higher efficiency as compared to an IMEX scheme based on a micro-macro decomposition have also been demonstrated. The schemes might be extended to other complicated problems, such as radiative transfer equations \cite{li2020unified,xiong2021high}. We will also study higher order model approximations and design corresponding higher order schemes in our future work.
	

	\normalem
	\bibliographystyle{siam}
	\bibliography{referenceXiongTao}

\begin{thebibliography}{10}

\bibitem{sphere}
{\em {SPHERE\_LEBEDEV\_RULE}},
  https://people.sc.fsu.edu/jburkardt/datasets/sphere\_lebedev\_rule/spherelebedevrule.html,
   (accessed 1 August 2022).

\bibitem{albi2020implicit}
{\sc G.~Albi, G.~Dimarco, and L.~Pareschi}, {\em Implicit-explicit multistep
  methods for hyperbolic systems with multiscale relaxation}, SIAM Journal on
  Scientific Computing, 42 (2020), pp.~A2402--A2435.

\bibitem{Bird1994}
{\sc G.~Bird}, {\em Molecular Gas Dynamics and Direct Simulation of Gas Flows},
  Clarendon Press, Oxford, 1994.

\bibitem{boscarino2013implicit}
{\sc S.~Boscarino, L.~Pareschi, and G.~Russo}, {\em Implicit-explicit
  \text{Rung--Kutta} schemes for hyperbolic systems and kinetic equations in
  the diffusion limit}, SIAM Journal on Scientific Computing, 35 (2013),
  pp.~A22--A51.

\bibitem{caflisch1997uniformly}
{\sc R.~Caflisch, S.~Jin, and G.~Russo}, {\em Uniformly accurate schemes for
  hyperbolic systems with relaxation}, SIAM Journal on Numerical Analysis, 34
  (1997), pp.~246--281.

\bibitem{case1967}
{\sc K.~Case and P.~Zweifel}, {\em Linear Transport theory}, Addison-Wesley,
  Reading, MA, 1967.

\bibitem{Chandrasekhar1960}
{\sc S.~Chandrasekhar}, {\em Radiative transfer}, Dover, New York, 1960.

\bibitem{cho2021conservative1}
{\sc S.~Cho, S.~Boscarino, G.~Russo, and S.~Yun}, {\em Conservative
  semi-lagrangian schemes for kinetic equations part i: Reconstruction},
  Journal of Computational Physics, 432 (2021), p.~110159.

\bibitem{crestetto2019asymptotically}
{\sc A.~Crestetto, N.~Crouseilles, G.~Dimarco, and M.~Lemou}, {\em
  Asymptotically complexity diminishing schemes ({ACDS}) for kinetic equations
  in the diffusive scaling}, Journal of Computational Physics, 394 (2019),
  pp.~243--262.

\bibitem{einkemmer2021asymptotic}
{\sc L.~Einkemmer, J.~Hu, and Y.~Wang}, {\em An asymptotic-preserving dynamical
  low-rank method for the multi-scale multi-dimensional linear transport
  equation}, Journal of Computational Physics, 439 (2021), p.~110353.

\bibitem{jang2014analysis}
{\sc J.~Jang, F.~Li, J.~Qiu, and T.~Xiong}, {\em Analysis of asymptotic
  preserving \text{DG-IMEX} schemes for linear kinetic transport equations in a
  diffusive scaling}, SIAM Journal on Numerical Analysis, 52 (2014),
  pp.~2048--2072.

\bibitem{jang2015high}
\leavevmode\vrule height 2pt depth -1.6pt width 23pt, {\em High order
  asymptotic preserving \text{DG-IMEX} schemes for discrete-velocity kinetic
  equations in a diffusive scaling}, Journal of Computational Physics, 281
  (2015), pp.~199--224.

\bibitem{jin1999efficient}
{\sc S.~Jin}, {\em Efficient asymptotic-preserving ({AP}) schemes for some
  multiscale kinetic equations}, SIAM Journal on Scientific Computing, 21
  (1999), pp.~441--454.

\bibitem{jin2010asymptotic}
\leavevmode\vrule height 2pt depth -1.6pt width 23pt, {\em Asymptotic
  preserving ({AP}) schemes for multiscale kinetic and hyperbolic equations: a
  review}, Lecture notes for summer school on methods and models of kinetic
  theory (M $\&$ MKT), Porto Ercole (Grosseto, Italy),  (2010), pp.~177--216.

\bibitem{jin2022asymptotic}
\leavevmode\vrule height 2pt depth -1.6pt width 23pt, {\em
  Asymptotic-preserving schemes for multiscale physical problems}, Acta
  Numerica, 31 (2022), pp.~1--82.

\bibitem{jin1998diffusive}
{\sc S.~Jin, L.~Pareschi, and G.~Toscani}, {\em Diffusive relaxation schemes
  for multiscale discrete-velocity kinetic equations}, SIAM Journal on
  Numerical Analysis, 35 (1998), pp.~2405--2439.

\bibitem{jin2000uniformly}
\leavevmode\vrule height 2pt depth -1.6pt width 23pt, {\em Uniformly accurate
  diffusive relaxation schemes for multiscale transport equations}, SIAM
  Journal on Numerical Analysis, 38 (2000), pp.~913--936.

\bibitem{klar1998asymptotic}
{\sc A.~Klar}, {\em An asymptotic-induced scheme for nonstationary transport
  equations in the diffusive limit}, SIAM Journal on Numerical Analysis, 35
  (1998), pp.~1073--1094.

\bibitem{klar2002uniform}
{\sc A.~Klar and A.~Unterreiter}, {\em Uniform stability of a finite difference
  scheme for transport equations in diffusive regimes}, SIAM Journal on
  Numerical Analysis, 40 (2002), pp.~891--913.

\bibitem{lemou2008new}
{\sc M.~Lemou and L.~Mieussens}, {\em A new asymptotic preserving scheme based
  on micro-macro formulation for linear kinetic equations in the diffusion
  limit}, SIAM Journal on Scientific Computing, 31 (2008), pp.~334--368.

\bibitem{li2017implicit}
{\sc Q.~Li and L.~Wang}, {\em Implicit asymptotic preserving method for linear
  transport equations}, Communications in Computational Physics, 22 (2017),
  pp.~157--181.

\bibitem{li2020unified}
{\sc W.~Li, C.~Liu, Y.~Zhu, J.~Zhang, and K.~Xu}, {\em Unified gas-kinetic
  wave-particle methods iii: Multiscale photon transport}, Journal of
  Computational Physics, 408 (2020), p.~109280.

\bibitem{liu2010analysis}
{\sc J.-G. Liu and L.~Mieussens}, {\em Analysis of an asymptotic preserving
  scheme for linear kinetic equations in the diffusion limit}, SIAM Journal on
  Numerical Analysis, 48 (2010), pp.~1474--1491.

\bibitem{mieussens2013asymptotic}
{\sc L.~Mieussens}, {\em On the asymptotic preserving property of the unified
  gas kinetic scheme for the diffusion limit of linear kinetic models}, Journal
  of Computational Physics, 253 (2013), pp.~138--156.

\bibitem{naldi1998numerical}
{\sc G.~Naldi and L.~Pareschi}, {\em Numerical schemes for kinetic equations in
  diffusive regimes}, Applied Mathematics Letters, 11 (1998), pp.~29--35.

\bibitem{peng2020stability}
{\sc Z.~Peng, Y.~Cheng, J.~Qiu, and F.~Li}, {\em Stability-enhanced \text{AP
  IMEX-LDG} schemes for linear kinetic transport equations under a diffusive
  scaling}, Journal of Computational Physics, 415 (2020), p.~109485.

\bibitem{peng2021stability}
\leavevmode\vrule height 2pt depth -1.6pt width 23pt, {\em Stability-enhanced
  \text{AP IMEX1-LDG} method: Energy-based stability and rigorous ap property},
  SIAM Journal on Numerical Analysis, 59 (2021), pp.~925--954.

\bibitem{peng2021asymptotic}
{\sc Z.~Peng and F.~Li}, {\em Asymptotic preserving \text{IMEX-DG-S} schemes
  for linear kinetic transport equations based on schur complement}, SIAM
  Journal on Scientific Computing, 43 (2021), pp.~A1194--A1220.

\bibitem{pomraning2005equations}
{\sc G.~C. Pomraning}, {\em The equations of radiation hydrodynamics}, Courier
  Corporation, 2005.

\bibitem{shi2020}
{\sc Y.~Shi, P.~Song, and W.~Sun}, {\em An asymptotic preserving unified gas
  kinetic particle method for radiative transfer equations}, Journal of
  Computational Physics, 420 (2020), p.~109687.

\bibitem{sun2015}
{\sc W.~Sun, S.~Jiang, and K.~Xu}, {\em An asymptotic preserving unified gas
  kinetic scheme for gray radiative transfer equations}, Journal of
  Computational Physics, 285 (2015), pp.~265--279.

\bibitem{van1997comparative}
{\sc G.~Van~Albada, B.~Van~Leer, and W.~Roberts}, {\em A comparative study of
  computational methods in cosmic gas dynamics}, in Upwind and high-resolution
  schemes, Springer, 1997, pp.~95--103.

\bibitem{xiong2021high}
{\sc T.~Xiong, W.~Sun, Y.~Shi, and P.~Song}, {\em High order asymptotic
  preserving discontinuous \text{Galerkin} methods for gray radiative transfer
  equations}, Journal of Computational Physics,  (2022), p.~111308.

\end{thebibliography}
\end{document}